\documentclass{article}
\extrafloats{64} 
\PassOptionsToPackage{numbers}{natbib}

\usepackage[preprint]{neurips_2026}

\usepackage[utf8]{inputenc} 
\usepackage[T1]{fontenc}    
\usepackage{hyperref}       
\usepackage{url}            
\usepackage{booktabs}       
\usepackage{amsfonts}       
\usepackage{amsmath,amssymb,array,tabularx,threeparttable,textcomp,longtable} 
\usepackage{nicefrac}       
\usepackage{microtype}      
\usepackage{xcolor}         
\usepackage{changepage}
\usepackage{graphicx}
\usepackage{todonotes}
\usepackage{tikz}
\usepackage{listings}
\usepackage{enumitem}
\usepackage{float}
\usepackage{algorithm}
\usepackage{algorithmic}
\let\originaltable\table
\let\endoriginaltable\endtable
\renewenvironment{table}[1][tbp]
  {\originaltable[#1]\begin{adjustwidth}{-3.2cm}{-3.2cm}}
  {\end{adjustwidth}\endoriginaltable}

\title{LongTails: Automated Research in Optimization}

\hypersetup{
  colorlinks   = true, 
  urlcolor     = blue, 
  linkcolor    = {red!75!black}, 
  citecolor   = {blue!50!black} 
}

\author{%
  Alexander Tyurin \\ AXXX, Moscow, Russia \\ Applied AI Institute, Moscow, Russia
  \And
  QuocViet Pham \\ Applied AI Institute, Moscow, Russia
  \And
  Valentin Khrulkov \\ AXXX, Moscow, Russia
  \And
  Ivan Oseledets \\ AXXX, Moscow, Russia
}

\newcommand{\numbernewresuls}{100}

\begin{document}

\maketitle
\vspace{-0.4cm}
\begin{center}
  \small
  GitHub: \url{https://github.com/k3nfalt/longtails}\\
  Website: \url{https://k3nfalt.github.io/longtails}
\end{center}
\vspace{0.3cm}

\begin{abstract}
We introduce LongTails, a framework for discovering and improving oracle-complexity bounds with formal verification in Lean 4.
Using LongTails, we have already obtained more than \numbernewresuls{} new results across a broad range of optimization problems, including generalized smoothness, heavy-tailed noise, finite-sum and composite optimization, and higher-order and zeroth-order methods. Our goal is to scale discovery and continually develop new methods with non-trivial and valuable theoretical guarantees using modern agentic systems.
The framework is specialized for optimization, where we separate a \emph{Setup}, specifying the problem class, assumptions, available information, and success criterion, from a \emph{Complexity}, specifying a sufficient oracle-call budget. Once a Setup is fixed (e.g., $L$--smooth convex optimization with the goal of finding $x_{\mathrm{out}}$ such that $f(x_{\mathrm{out}}) - f(x^\star) \leq \varepsilon$), agents search for better Complexities (e.g., $\mathcal{O}(\sqrt{L R^2 / \varepsilon})$, where $R \geq \|x_0 - x^\star\|$) by designing algorithms and proving their convergence. The framework then checks the proofs and guarantees that every accepted Complexity can indeed be attained under the given Setup. This turns the search for new methods into an ongoing effort to improve their theoretical guarantees. To the best of our knowledge, LongTails is the first framework to scale automated research in optimization across such a broad range of problems with formally verified theoretical guarantees.
\end{abstract}

\section{Introduction}
The current agentic systems~\citep{openaiChatGPT,anthropicClaudeCode,googleGemini} are becoming incredibly successful in different areas, including software engineering~\citep{yang2024sweagent} and automated scientific experimentation~\citep{lu2024aiscientist}. 
Mathematics is one of the fields affected by these advances. For example, OpenAI reported a proof of finite-time blowup for the three-dimensional Navier--Stokes equations~\citep{openai2026NavierStokes} and released over 700 mathematical manuscripts~\citep{openai2026MathCatalog} that solve challenging mathematical problems.

Mathematical optimization is one of the subfields of mathematics~\citep{nsscNesterov2018,lanbook}, where our goal is to design effective algorithms to solve optimization problems under different setups. Given a setup, we not only want to design a new method but, more importantly, prove that its theoretical convergence rate is sound and ideally better than those of previous methods. 

In this work, we consider the problem of designing an agentic discovery framework specialized for mathematical optimization, online learning, and other numerical methods. 
\begin{quote}
Our goal is to design a universal framework, which would allow us to \textbf{scale discovery and continually develop new methods} with non-trivial and valuable theoretical guarantees.
\end{quote}

\textbf{Contributions.}
In this paper, we present a new framework, LongTails, for automated research in numerical mathematics and optimization. The main novelty is a new formalization in Lean 4~\citep{deMouraUllrich2021Lean4} tailored to optimization. By formalizing Setups and Complexities, which arise naturally in numerical methods, we design the framework to support the discovery and development of new theoretical results at scale. 

In particular, a Setup fixes the problem class, the information available to a method, the oracles it can use, and the required accuracy criterion. Agents can then search for new algorithms and improve their Complexities, while the Setup remains fixed. Thus, LongTails verifies that there exists an algorithm that uses only the allowed information and achieves the required accuracy within the claimed oracle budget (Complexity). This allows us to use the same formal Lean 4 interface and verification procedure across different optimization problems and repeatedly improve the guarantees within each Setup. Therefore, we reduce the problem of searching for new methods and complexities into the problem of improving Complexities under a fixed Setup.

We have already obtained more than \numbernewresuls{} new results across several active areas of optimization, including optimization under generalized smoothness, stochastic optimization with heavy-tailed noise, finite-sum and composite optimization, parameter-free methods, optimization with inexact gradients, second- and higher-order methods, zeroth-order optimization, and monotone inclusions. To the best of our knowledge, this is the first framework capable of producing optimization results at this scale using modern agentic systems. 
See Sections~\ref{gen_inst}, \ref{sec:comparions}, and~\ref{sec:add_comparisons}, and our \href{https://k3nfalt.github.io/longtails/}{website}.

\textbf{Structure.} Section~\ref{sec:user-guide} introduces LongTails, Setups, and complexities through an example. Section~\ref{sec:verification} explains how LongTails certifies that each accepted Complexity is attainable under a given Setup. Section~\ref{gen_inst} highlights selected new theoretical results, with detailed comparisons in Sections~\ref{sec:comparions} and~\ref{sec:add_comparisons}. Finally, Section~\ref{sec:prompting} describes our prompting strategy.

\section{A User's Guide to LongTails}
\label{sec:user-guide}
We now present our new framework developed in Lean 4. The main idea is to consider the structure of optimization and many numerical methods: we have a \emph{Setup}, and our goal is to get the best possible \emph{Complexity}. Let us illustrate this with an example.

\textbf{Example with standard convex and $L$--smooth functions.} \\
Consider the classical problem
\begin{align*}
  \min\limits_{x\in\mathbb{R}^d} f(x),
\end{align*}
where $f$ is differentiable and convex, and its gradient satisfies
$\|\nabla f(x)-\nabla f(y)\|\le L\|x-y\|$ for all $x,y\in\mathbb{R}^d$.
Assume that $L>0$ and that a minimizer $x^\star$ exists with
$\|x_0-x^\star\|\le R$, where $R\ge0$. The method has access to an exact
gradient oracle and must produce $x_{\mathrm{out}}$ such that
$f(x_{\mathrm{out}})-f(x^\star)\le\varepsilon$, for a prescribed
$\varepsilon>0$. These assumptions, the oracle, and the accuracy criterion
form the Setup. The accelerated gradient method~\citep{nesterov1983accelerated} achieves this criterion within
\begin{align}
  \label{eq:agd}
  N = \left\lceil\sqrt{\frac{2LR^2}{\varepsilon}}\right\rceil
\end{align}
gradient queries. This oracle budget is its Complexity.

Motivated by this and many other optimization problems, we develop a framework in which a researcher defines a Setup that describes the problem of interest. Then, potentially with the help of AI, the researcher proves a new Complexity. The framework handles the rest: once a Setup is fixed, our framework ensures that the proved Complexity is indeed correct and verified with Lean 4. 

Here is a simplified presentation of the Setup for the previously considered example with convex and $L$--smooth functions. For readability, we omit the declarations making $E$ a real Hilbert space, their instance registration, and the oracle's random seeds:
\noindent\begin{minipage}{\linewidth}
\begin{lstlisting}[
  basicstyle=\ttfamily\footnotesize,
  sensitive=true,
  morekeywords={structure,where,attribute,instance,extends,abbrev,Type},
  keywordstyle={\color[HTML]{0000FF}},
  emph={NormedAddCommGroup,InnerProductSpace,CompleteSpace},
  emphstyle={\color[HTML]{267F99}},
  columns=fullflexible,
  keepspaces=true,
  breaklines=true,
  showstringspaces=false,
  frame=single,
  frameround=tttt,
  framesep=6pt,
  framerule=0.4pt,
  xleftmargin=6.4pt,
  xrightmargin=6.4pt,
  literate={ℝ}{{$\mathbb{R}$}}1 {ε}{{$\varepsilon$}}1
    {→}{{$\to$}}1 {∧}{{$\land$}}1 {≤}{{$\le$}}1
    {∀}{{$\forall$}}1 {∃}{{$\exists$}}1 {‖}{{$\Vert$}}1
    {ᵐ}{{$^{\mathrm{m}}$}}1 {ω}{{$\omega$}}1 {∂}{{$\partial$}}1
]
structure Public where
  E : Type u
  ...
  L : ℝ
  x0 : E

structure Params extends Public where
  R : ℝ
  ε : ℝ

abbrev spec : Setup where
  Params := Params
  Public := Public
  Instance p := p.E → ℝ
  Query p := p.E
  Reply p := p.E
  Output p := p.E
  publicInput := Params.toPublic
  valid q := 0 < q.L ∧ 0 ≤ q.R ∧ 0 < q.ε
  admissible q f oracle :=
    ConvexOn ℝ Set.univ f ∧
    (∀ k x, HasGradientAt f (oracle.answer k (0, x)) x) ∧
    (∀ k x y,
      ‖oracle.answer k (0, x) - oracle.answer k (0, y)‖ ≤ q.L * ‖x - y‖) ∧
    (∃ xstar, IsMinOn f Set.univ xstar ∧ ‖q.x0 - xstar‖ ≤ q.R)
  successful q f output P :=
    ∃ xstar, IsMinOn f Set.univ xstar ∧
      ∀ᵐ ω ∂(P : MeasureTheory.Measure _), f (output ω) - f xstar ≤ q.ε
\end{lstlisting}
\end{minipage}
\par
Its fields for this example have the following interpretation:
\begin{enumerate}[wide=0pt]
  \item \texttt{Params} stores all parameters of the problem, including the Hilbert space $E$, the smoothness constant $L$, the starting point $x_0$, the initial distance bound $R$, and the target accuracy $\varepsilon$. In contrast, \texttt{Public} contains only the parameters available to the algorithm as input: $E$, $L$, and $x_0$. The algorithm does not have access to the target accuracy $\varepsilon$ or the initial distance bound $R$. 
  For instance, making \texttt{Public} almost empty allows us to build parameter-free methods.
  \item \texttt{Instance} is the type of objectives $f:E\to\mathbb{R}$ that we minimize;
  \texttt{Query}, \texttt{Reply}, and \texttt{Output} are all $E$: a query is
  a point in $E$ that an algorithm sends to the oracle; the oracle's reply is a gradient in $E$; and the output is a candidate minimizer that the algorithm returns at the end.
  \item \texttt{valid} specifies that the Setup assumes $L>0$, $R\ge0$, and $\varepsilon>0$.
  \item \texttt{admissible} is the field specifying the assumptions: it requires convexity of $f$, exact gradient replies from the oracle, the $L$--Lipschitz gradient bound, and the existence of a minimizer within distance $R$ of $x_0$.
  \item \texttt{successful} defines the goal: it requires the algorithm to return $x_{\mathrm{out}}$ such that
  $f(x_{\mathrm{out}})-f(x^\star)\le\varepsilon$ almost surely.
\end{enumerate}
This code formalizes in Lean 4 the assumptions and accuracy criterion typically stated in prose in mathematical papers~\citep{nesterov1983accelerated,nsscNesterov2018}. Once a Setup is fixed, the goal is to establish a Complexity: an explicit bound on the number of oracle calls sufficient for a proposed algorithm to meet the specified accuracy criterion. 
In our framework, omitting type annotations, one should prove a theorem similar to this:
\noindent\begin{minipage}{\linewidth}
\begin{lstlisting}[
  basicstyle=\ttfamily\footnotesize,
  sensitive=true,
  morekeywords={noncomputable,def,let,fun,theorem,by},
  keywordstyle={\color[HTML]{0000FF}},
  emph={Fin},
  emphstyle={\color[HTML]{267F99}},
  columns=fullflexible,
  keepspaces=true,
  breaklines=true,
  showstringspaces=false,
  frame=single,
  frameround=tttt,
  framesep=6pt,
  framerule=0.4pt,
  xleftmargin=6.4pt,
  xrightmargin=6.4pt,
  literate={ε}{{$\varepsilon$}}1 {→}{{$\to$}}1
    {ℕ}{{$\mathbb{N}$}}1 {⌈}{{$\lceil$}}1 {⌉}{{$\rceil$}}1
    {₊}{{$_+$}}1
]
noncomputable def budget := fun q _ => max 1 ⌈Real.sqrt (2 * q.L * q.R ^ 2 / q.ε)⌉₊

theorem result : Convergence spec budget := by ...
\end{lstlisting}
\end{minipage}
\par
The theorem states that there exists an algorithm that, for every valid choice of parameters and every objective and oracle satisfying the Setup's assumptions, returns a point $x_{\mathrm{out}}$ such that $f(x_{\mathrm{out}})-f(x^\star)\le\varepsilon$ almost surely, using at most $\max\left\{1,\left\lceil\sqrt{\frac{2LR^2}{\varepsilon}}\right\rceil\right\}$ gradient oracle calls. This matches the bound in \eqref{eq:agd}.

\textbf{General case.} The same approach extends to a broad range of problems, including nonconvex optimization, stochastic optimization with heavy-tailed noise, finite-sum optimization, and optimization using second-order information, as well as many other settings. The corresponding \href{https://github.com/k3nfalt/longtails/tree/main/LongTails/Setups}{Setups} and \href{https://github.com/k3nfalt/longtails/tree/main/LongTails/Theorems}{Theorems} are available in our repository.

\section{How Are We Sure That the Complexities Are Correct?}
\label{sec:verification}
The main contribution of the framework is to ensure that a proved Complexity is indeed correct and that there exists an algorithm that achieves it under the given Setup. Using Comparator~\citep{comparatorRepository}, we check a proof of the theorem \texttt{result\ :\ Convergence\ spec\ budget} from the preceding example. This statement can be interpreted as ``prove that there exists an algorithm that achieves the required accuracy within the given oracle budget under the specified Setup.''

The central part is how we define ``convergence.'' In \href{https://github.com/k3nfalt/longtails/blob/main/LongTails/Template.lean}{\texttt{Template.lean}}, \texttt{Convergence} takes a Setup and a Complexity (\texttt{budget}) as input. Omitting type annotations and technical details, we can write it schematically as follows:
\noindent\begin{minipage}{\linewidth}
\begin{lstlisting}[
  basicstyle=\ttfamily\footnotesize,
  sensitive=true,
  morekeywords={def,let,fun,Prop},
  keywordstyle={\color[HTML]{0000FF}},
  emphstyle={\color[HTML]{267F99}},
  columns=fullflexible,
  keepspaces=true,
  breaklines=true,
  showstringspaces=false,
  frame=single,
  frameround=tttt,
  framesep=6pt,
  framerule=0.4pt,
  xleftmargin=6.4pt,
  xrightmargin=6.4pt,
  literate={→}{{$\to$}}1 {ℕ}{{$\mathbb{N}$}}1
    {∃}{{$\exists$}}1 {∀}{{$\forall$}}1 {∧}{{$\land$}}1
    {×}{{$\times$}}1 {⟨}{{$\langle$}}1 {⟩}{{$\rangle$}}1
]
def Convergence Setup budget : Prop :=
  ∃ makeAlgorithm,
    ∀ p, Setup.valid p →
      ∀ I oracle, Setup.admissible p I oracle →
        let A := makeAlgorithm (Setup.publicInput p)
        let result := execute oracle A (budget p)
        Setup.successful p I result
\end{lstlisting}
\end{minipage}
\par
This proposition can be read step by step. First, there must exist an algorithm constructor ($\exists$ \texttt{makeAlgorithm}). For every choice of parameters \texttt{p} satisfying \texttt{Setup.valid p}, and every instance \texttt{I} and oracle \texttt{oracle} satisfying \texttt{Setup.admissible p I oracle}, this constructor produces an algorithm \texttt{A} using only the public parameters \texttt{Setup.publicInput p}. Next, \texttt{execute} runs \texttt{A} with \texttt{oracle} within the oracle-call budget \texttt{budget p}. Finally, \texttt{Setup.successful p I result} requires the output \texttt{result} to satisfy the specified accuracy criterion.

\emph{The order of the quantifiers matters.} The algorithm constructor is chosen before the hidden instance and oracle, and receives neither of them as input. It cannot choose a different method by inspecting an unknown objective or its minimizer.

Comparator~\citep{comparatorRepository} checks the theorem \texttt{result\ :\ Convergence\ spec\ budget} against the fixed Setup and audits the budget definition, rejecting unfinished proofs and unpermitted axioms. Thus, we reduce the research problem to agentic optimization of Complexities: agents search for better budgets and corresponding algorithms and proofs, while every accepted result must respect the fixed Setup's assumptions, oracle model, and success criterion.

Schematically, for a fixed Setup, the research objective in optimization is
\[
\begin{tikzpicture}[baseline=(scheme.base)]
\node[draw, rounded corners=5pt, inner sep=8pt] (scheme) {$\displaystyle
\begin{aligned}
&\underset{\texttt{budget}}{\operatorname{minimize}}
&& \texttt{budget} \\
&\text{subject to}
&& \texttt{Convergence}\;\texttt{Setup}\;\texttt{budget}
\quad\text{proved in Lean 4}.
\end{aligned}
$};
\end{tikzpicture}
\]
Here, budget represents the Complexity: a bound on the number of calls to each oracle. The constraint requires an algorithm that achieves this budget under the fixed Setup. The goal of an agent is to minimize budget.

\section{Selected New Theoretical Results}
\label{gen_inst}
Using LongTails, we search for new methods that improve existing Complexities or provide the first guarantees in unexplored Setups.
We use Codex~\citep{openaiCodex} and Claude Code~\citep{anthropicClaudeCode} to consider Setups and produce Complexities. We have proved more than \numbernewresuls{} new results on our \href{https://k3nfalt.github.io/longtails/}{website}. 
We use an additional agentic review layer to produce reviews of the results. We then manually check the Setups, reviews, and Complexities, including whether the Setups match those in prior work or capture the intended mathematical problems. We count only the results for which we agree with the review that the result is new.
The proofs and theorem statements are checked automatically, and the resulting Complexities are converted using a deterministic parser.

In Sections~\ref{sec:comparions} and \ref{sec:add_comparisons}, we present the considered Setups and our results, all formally proved in Lean 4. Below, we highlight several new theoretical convergence guarantees.

\textbf{Generalized Smooth Convex Optimization.} Let us discuss one example in more detail. Consider convex optimization under the generalized smoothness assumption $\|\nabla^2 f(x)\|\le L_0+L_1\|\nabla f(x)\|,$
where $L_0>0$ and $L_1\ge0$ \citep{zhang2020gradientClipping}. Given bounds $R\ge\|x_0-x^\star\|$ and $\Delta\ge f(x_0)-f(x^\star)$, our goal is to find $x_{\mathrm{out}}$ with $f(x_{\mathrm{out}})-f(x^\star)\le\varepsilon$.
For simplicity, we hide logarithmic factors using $\widetilde{\mathcal{O}}$ notation. The previous bounds in~\citep[Thm.~4.3]{cgsTyurin2026v2} and~\citep[Thm.~6.2]{gscVankov2025v3}, respectively, are
\[
\widetilde{\mathcal{O}}\!\left(1+\sqrt{\frac{L_0R^2}{\varepsilon}}+L_1R\right)
\quad\text{and}\quad
\widetilde{\mathcal{O}}\!\left((1+\nu)\left[1+\sqrt{\frac{L_0R^2}{\varepsilon}}+(L_1R)^{2/3}\right]\right).
\]
The first bound counts gradient queries. The second also counts function and gradient queries used by exact line searches, where $\nu$ is the oracle cost of an exact one-dimensional minimization.
Using LongTails, we obtain
\[
\widetilde{\mathcal{O}}\!\left(1+\sqrt{\frac{L_0R^2}{\varepsilon}}+(L_1R)^{2/3}\right)
\]
using only gradient queries. Thus, ignoring logarithmic factors, we improve the $L_1R$ term in the first bound to $(L_1R)^{2/3}$ and remove the exact line-search cost from the second. This answers the question of achieving this accelerated rate using only gradient queries \citep{gscVankov2025v3}. The precise bounds and comparisons are presented in Table~\ref{tab:GeneralizedSmoothAndConvex}. The method is presented in Algorithm~\ref{alg:damped-generalized} in Appendix~\ref{sec:generalized-method}.

The agent came up with a very interesting and non-trivial idea, presented in Algorithm~\ref{alg:damped-generalized}: combine acceleration with a local curvature estimate obtained from the gradient norm. At iteration $k$, the method uses the curvature estimate $G_k$ to choose the acceleration weight $a_k$ and form the extrapolated point $y^k$. It then computes $g^k=\nabla f(y^k)$, sets $M_k=4(L_0+L_1\|g^k\|)$, and takes the gradient step $v^k=y^k-g^k/M_k$. If $M_k>G_k$, the method damps the accelerated update by $\gamma_k=G_k/M_k$ and sets $G_{k+1}=16G_k$. Otherwise, it takes the full update and sets $G_{k+1}=G_k/16$. The second gradient, $h^k=\nabla f(v^k)$, is used to update $z^{k+1}$. The analysis, outlined in Section~\ref{sec:generalized-proof-idea}, separates iterations with low and high local curvature: the former give the usual accelerated term, while the latter are controlled by the decrease of a potential function, yielding the $(L_1R)^{2/3}$ term up to logarithmic factors.

There are many other results.
In Table~\ref{tab:GeneralizedSmoothParameterFree}, we establish
an accelerated, gradient-only guarantee whose updates require
neither the smoothness parameters nor bounds on the initial
distance and objective gap.
In Table~\ref{tab:CompositeGeneralizedSmooth}, we establish an
accelerated complexity bound for constrained composite problems
$f+\psi$, where $f$ is convex and $(L_0,L_1)$-smooth and $\psi$
is convex and admits exact proximal evaluations. 
In Table~\ref{tab:CompositeGeneralizedParameterFree}, we further
obtain an accelerated composite guarantee whose updates use only
gradients and proximal evaluations, with an arbitrary initial trial
step and no smoothness parameters or initial distance and gap bounds.
Compared with~\citep{LiLaiAnHuWen2026v3Composite}, we retain the leading
accelerated complexity while removing these parameters, including
the gradient bound at a minimizer, from the updates; the stopping
budget still depends on them.

\textbf{Generalized Smooth Composite Optimization.} In Table~\ref{tab:SlidingGeneralizedSmooth}, we extend
gradient sliding~\citep{sgsLan2016} to convex objectives $f+h$
with $(L_0,L_1)$-smooth $f$ and bounded subgradients of $h$.
Using exact projections, we obtain $\varepsilon^{-1/2}$ dependence
in gradient calls with polynomial dependence on $L_1R$, improving
the exponential dependence from applying classical gradient sliding
with a localized smoothness bound. The total complexity still
includes $\mathcal{O}(M^2R^2/\varepsilon^2)$ subgradient and
projection calls, where $M$ bounds the subgradients of $h$;
the gradient savings improve the total bound when this term
does not dominate.

\textbf{Generalized Smooth Convex Stochastic Optimization.} For convex $(L_0,L_1)$-smooth objectives, we also obtain
stochastic guarantees.
In Table~\ref{tab:StochasticGeneralizedSmoothAndConvex}, we obtain
accelerated high-probability bounds under unbiased, light-tailed
gradient noise, improving the dependence on $L_1R$ in the regimes
specified in the table.
In Table~\ref{tab:FiniteVarianceGeneralizedSmoothAndConvex},
we obtain accelerated high-probability guarantees with unbiased
stochastic gradients under only a bounded-variance assumption. Up to logarithmic factors, the complexity
is uniformly no worse than the RCSAG bound~\citep{fvLobanov2026v1} on the common
setup, improving the dependence on $L_1R$ from linear to
$(L_1R)^{2/3}$. Our adaptive construction also reduces the
confidence-related logarithmic overhead on the optimization term
in the confidence regimes specified in the table.

\textbf{Smooth Convex Relative-Error Optimization.} In Table~\ref{tab:RelativeErrorSmoothConvex}, we consider
smooth convex optimization with relative gradient errors
$\|g(x)-\nabla f(x)\|\le\alpha\|\nabla f(x)\|$.
For $0\le\alpha\le1/3$, we improve the regularized accelerated bound of~\citep{igsVasinGasnikov2026v1} from
$\mathcal{O}((\sqrt{T}+\alpha T)(1+\log T))$ to
$\mathcal{O}(\sqrt{T}+\alpha^2T)$, where $T=LR^2/\varepsilon>1$.
Our last-iterate guarantee also covers every $\alpha<1.$

\textbf{Strongly Convex Inexact-Gradient Optimization.} In Table~\ref{tab:InexactGradientStronglyConvex}, we consider
smooth strongly convex optimization with inexact gradients
satisfying
$\|g(x)-\nabla f(x)\|\le\alpha\|\nabla f(x)\|+\delta$.
We improve the dependence on the relative error level $\alpha$
and reduce the guaranteed accuracy threshold caused by the
absolute error $\delta$, compared with the bounds and regimes
specified in the table. The guarantee does not require
unbiased or independent errors.

\textbf{BFGS for Linear-Predictor Objectives.} In Table~\ref{tab:LinearPredictorBFGS}, we establish an accelerated-order guarantee for standard BFGS with Armijo--Wolfe line search on objectives $f(x)=h(Ax)$, where $h$ is convex and smooth and satisfies the specified quadratic-growth condition. For fixed remaining parameters, the number of iterations needed to achieve $f(x_k)-f^\star\le\varepsilon$ improves from $\mathcal{O}(\varepsilon^{-1})$ in~\citep{lpbDingYangZhou2026v1} to $\mathcal{O}(\varepsilon^{-1/2})$ on this subclass. This is a BFGS-specific guarantee.

\textbf{Similar Finite-Sum Quasar Optimization.} In Table~\ref{tab:ContinuizedSimilarQuasar}, we consider smooth strongly
quasar-convex finite sums under a similarity assumption on the
component gradients. We obtain improved component-gradient
complexities in high-similarity regimes, with deterministic
oracle budgets and expected last-iterate objective-gap guarantees.

\textbf{Adam under Heavy-Tailed Gradient Noise.} In Table~\ref{tab:HeavyTailedAdam}, we establish a convergence guarantee for coordinate Adam on smooth nonconvex objectives under a bounded coordinatewise $\alpha$-moment noise assumption, $1<\alpha\le2$. A uniformly sampled iterate satisfies $\mathbb{E}\|\nabla f(\widehat x)\|\le\varepsilon$ with $\mathcal{O}(\varepsilon^{-2-\alpha/(\alpha-1)})$ stochastic-gradient calls for fixed remaining parameters. This matches the known accuracy exponent of normalized stochastic gradient descent with momentum~\citep{htaHubler2025}. To the best of our knowledge, this is the first convergence guarantee for coordinate Adam under heavy-tailed gradient noise that achieves the optimal accuracy exponent.

\textbf{Nonconvex Optimization with a Similar Auxiliary Function.} In Table~\ref{tab:AuxiliarySimilarityNonconvex}, we consider a target $f$ with a $\rho$-Lipschitz Hessian and a cheap auxiliary function $h$ with an $L$-Lipschitz gradient, where $\nabla(f-h)$ is $\delta$-Lipschitz. Neither function is required to be convex. With initial objective gap at most $\Delta$, we find $\widehat x$ satisfying $\|\nabla f(\widehat x)\|\le\varepsilon$ using
$\widetilde{\mathcal{O}}\bigl((1+\Delta\sqrt{\rho}\,\varepsilon^{-3/2})\sqrt{1+\delta/\sqrt{\rho\varepsilon}}\bigr)$
expensive value--gradient calls to $f$.
In the small-accuracy regime, the leading term replaces the $\sqrt{L+\delta}$ dependence of restarted AGD~\citep{gtoLiLin2023} with $\sqrt{\delta}$, up to logarithmic factors. This improves the expensive-oracle bound when $L$ is sufficiently larger than $\delta$ to offset those factors. The method also requires $\widetilde{\mathcal{O}}(1+L/\delta)$ cheap value--gradient calls to $h$ per expensive call; the improvement therefore concerns expensive-oracle complexity, with both costs accounted for separately.

\textbf{Nonconvex Optimization with Comparison Oracles.} In Table~\ref{tab:ComparisonHessianLipschitz}, we establish a deterministic guarantee for visiting an $\varepsilon$-stationary point using only function-value comparisons. For functions with an $L$-Lipschitz gradient and a $\rho$-Lipschitz Hessian, the leading complexity is $\widetilde{\mathcal{O}}(n\Delta\sqrt{L}\,\rho^{1/4}\varepsilon^{-7/4})$, where $n$ is the dimension and $\Delta$ bounds the initial objective gap. Compared with~\citep{chltWang2026}, the dimension dependence improves from quadratic to linear, while the accuracy dependence worsens from $\varepsilon^{-3/2}$ to $\varepsilon^{-7/4}$. The total bound improves when $n\gg(L^2/(\rho\varepsilon))^{1/4}$, beyond logarithmic factors.

\textbf{Nonconvex Optimization with Lipschitz Third Derivatives.} In Table~\ref{tab:ThirdDerivativeSmoothNonconvex}, we obtain a dimension-free gradient--Hessian tradeoff for finding $\widehat x$ with $\|\nabla f(\widehat x)\|\le\varepsilon$, using only function values, gradients, and Hessians. For fixed remaining parameters, the Hessian count improves from $\mathcal{O}(d+d^{1/3}\varepsilon^{-4/3})$ in~\citep{tdsDoikovGrapiglia2025v1} to $\mathcal{O}(1+\varepsilon^{-1})$, while the value--gradient count is $\widetilde{\mathcal{O}}(1+\varepsilon^{-5/3})$. Here $d$ is the dimension. If a Hessian costs $B$ value--gradient calls, the total cost improves when $Bd^{1/3}\gg\varepsilon^{-1/3}$, up to logarithmic factors, in the small-accuracy regime. The result trades more value--gradient calls for fewer Hessian calls.

\textbf{Gap-Dependent Convex Stochastic Optimization.} In Table~\ref{tab:GapSmoothStrongGrowth}, we establish
a high-probability accelerated guarantee for convex objectives
under $(H_0,H_1)$ gap-dependent smoothness and a second-moment
strong-growth condition. Compared with RCSAG~\citep{gssgLobanov2026}, the stochastic accuracy term improves from
$\varepsilon^{-1}$ to $\varepsilon^{-1/2}$, up to logarithmic
factors, yielding a better total complexity when that term
dominates.

\textbf{Generalized Lipschitz Monotone Root Finding.} In Table~\ref{tab:GeneralizedLipschitzMonotone}, we improve
the residual complexity for monotone operators under
$(L_0,L_1)$ generalized Lipschitz continuity. For fixed remaining
parameters, the leading accuracy dependence improves from
$\varepsilon^{-2}$ in the extragradient bounds~\citep{glmChoudhuryLoizou2025v1,glmYoonEtAl2026v1} to
$\varepsilon^{-1}$, for the criterion
$\|F(x_{\mathrm{out}})\|\le\varepsilon$. In Table~\ref{tab:StochasticGeneralizedMonotone}, we establish
a high-probability residual guarantee under the same generalized
Lipschitz condition with unbiased, finite-variance stochastic
oracles, without a quasi-sharpness assumption. Up to logarithmic
factors, the oracle complexity is
$\mathcal{O}(1+L_0R/\varepsilon+(L_1R)^2+\sigma^2/\varepsilon^2)$,
matching the polynomial accuracy dependence of the globally Lipschitz result~\citep{sgmJiLanZhu2026v1} while extending to $L_1>0$.

\textbf{Constrained Generalized Monotone Variational Inequalities.} In Table~\ref{tab:ConstrainedGeneralizedMonotone}, we extend
$(L_0,L_1)$-Lipschitz monotone guarantees to constrained variational
inequalities with exact projections. The method returns a feasible
last iterate $x$ with a normal-cone certificate $\xi\in N_X(x)$
satisfying $\|F(x)+\xi\|\le\varepsilon$.
For fixed remaining parameters, the total operator and projection
complexity is $\mathcal{O}(\varepsilon^{-1})$, improving the derived
$\mathcal{O}(\varepsilon^{-2})$ backtracking bound
of~\citep{cgmShen2026v2}. Its constants are polynomial in $L_1R$,
avoiding the potentially exponential dependence from Lipschitz localization. The updates require known $L_0,L_1$, and the operator
must satisfy monotonicity and generalized Lipschitz continuity on
the whole ambient space. On the unconstrained overlap ($X=E$, $G=0$), the constrained bound has an extra
$1+L_1R$ factor in its leading term compared with our unconstrained
result, where $R$ bounds the initial distance to a solution.


\textbf{Exact-Oracle Monotone Root Finding.} In Table~\ref{tab:HolderMonotone}, we establish a residual
guarantee for monotone equations under the mixed continuity condition
$\|F(x)-F(y)\|\le L\|x-y\|+H\|x-y\|^\alpha$,
without requiring $L,H,\alpha$ or the solution-distance bound $R$
in the updates. For fixed remaining parameters, the complexity
is $\widetilde{\mathcal{O}}(\varepsilon^{-1/\alpha})$,
improving the parameter-free pure H\"older bound $\mathcal{O}(\varepsilon^{-2/\alpha})$ of PolyakEG-LS~\citep{hmYoon2026v1}.
The guarantee concerns a stopped iterate with
$\|F(\widehat x)\|\le\varepsilon$; in the Lipschitz case,
it incurs an extra logarithmic factor compared with the known-parameter anchored extragradient bound~\citep{glmYoonRyu2021v2}.

\textbf{Hessian-Lipschitz Stochastic Minimax Optimization.} In Table~\ref{tab:StochasticHessianNcsc}, we consider
nonconvex--strongly-concave minimax optimization with a Lipschitz
Hessian and unbiased, finite-variance stochastic gradients.
For $\Phi(x)=\max_y f(x,y)$, we obtain
$\mathbb{E}\|\nabla\Phi(\widehat x)\|\le\varepsilon$ using only
stochastic gradient calls. Up to logarithmic factors, the leading
accuracy dependence improves from $\varepsilon^{-4}$ in the
gradient-only bound~\citep{shncscWu2026v1} to $\varepsilon^{-7/2}$ for fixed
remaining parameters and sufficiently small $\varepsilon$.
This improvement uses Hessian Lipschitzness and comes with worse
dependence on the condition number $\kappa=L/\mu$.


\textbf{Convex Optimization with Growing Oracle Moments.} In Table~\ref{tab:TerminalClippingConvex}, we establish
a high-probability last-iterate guarantee for projected convex
optimization with unbiased stochastic subgradients satisfying
$\mathbb{E}\|g(x,\xi)\|^2\le B^2\|x-x_0\|^2+G^2$.
The oracle complexity is
$\mathcal{O}(1+(BR^2+G^2/B)^2(1+\log(1/\delta))^2/\varepsilon^2)$.
Compared with the expected last-iterate bound for projected Halpern~\citep{tcAMW2026v3}, converted using Markov's inequality, we replace polynomial dependence on
$1/\delta$ by logarithmic dependence and remove the accuracy
logarithm at fixed confidence. The guarantee requires neither
smoothness nor uniformly bounded subgradients.

\textbf{Smooth Convex Zeroth-Order Optimization.} In Table~\ref{tab:ZerothOrderSmoothAndConvex}, we establish
an accelerated two-point guarantee for smooth convex optimization
with independent, unbiased function-value noise of variance at most
$\sigma^2$. The retuned two-point SGD bound derived
from~\citep{zoscBachPerchet2016} is $\mathcal{O}(dT)$ in the low-noise regime;
we obtain
$\mathcal{O}(d\sqrt{T}+d^2T\sigma/\varepsilon)$ when
$\sigma<\varepsilon/d$, where $T=LR^2/\varepsilon\ge1$.
In particular, when $\sigma\le\varepsilon/(d\sqrt{T})$, the bound
is $\mathcal{O}(d\sqrt{T})$. The guarantee is an expected objective
gap at the last iterate; for $\sigma\ge\varepsilon/d$, the
complexity matches the order of the retuned baseline.

\textbf{Matrix-Generalized Smooth Zeroth-Order Optimization.} In Table~\ref{tab:ZerothOrderMatrixGeneralizedSmooth}, we establish
an accelerated high-probability guarantee for convex optimization
under matrix $(L_0,L_1)$-smoothness using only exact function values.
The method does not require the smoothness matrix $H$; its complexity
uses a supplied trace bound $\operatorname{tr}H\le r$ and has no
dependence on $\lambda_{\min}(H)^{-1}$.
Up to logarithmic factors, the leading accuracy term is
$\sqrt{d(r+1)L_0R^2/\varepsilon}$.
We extend the matrix-smooth RS-NAG result~\citep{zmgsOmiya2026v1} to $L_1>0$
and replace directional-derivative access by function values,
matching its leading rate when $L_1=0$ with additional lower-order
terms specified in the table.

\textbf{Known-Value Local Linear Optimization.} In Table~\ref{tab:KnownValueSmoothConvexLocalLinear}, we establish
an accelerated guarantee for smooth convex optimization over a closed
convex set using exact value--gradient queries and point-only linear
minimization over the feasible set intersected with a ball.
With the exact optimal value known, we obtain a feasible last primal
iterate with objective gap at most $\varepsilon$ using
$\mathcal{O}(1+\sqrt{LR^2/\varepsilon})$ combined oracle calls,
including unsuccessful radius trials. This improves the
$\mathcal{O}(1+LR^2/\varepsilon)$ fixed-radius BroxGD bound
of~\citep[Thm.~3.2(ii)]{kvBrox2026v2} when $LR^2/\varepsilon\ge1$,
but requires the exact optimal value and $L$, which that baseline
does not need. The result permits unbounded feasible sets and has
no additional logarithmic search factor.

\textbf{Smooth Convex Splitting with Proximal Access.} In Table~\ref{tab:DouglasRachfordSmoothConvex}, we establish
an accelerated last-iterate objective-gap guarantee for fast
Douglas--Rachford splitting with a general smooth convex term,
extending the analysis of~\citep{patrinosStellaBemporad2014DRv2} beyond quadratic objectives.
Using only the two proximal oracles at a fixed parameter $\gamma$,
the method achieves
$\mathcal{O}(1+R/\sqrt{\lambda\gamma\varepsilon})$ proximal calls
under $0<\lambda\le1$ and $4\gamma L\le\lambda$.
This extends the scope of the accelerated rate for this method,
with a more restrictive step size than the quadratic result of~\citep{patrinosStellaBemporad2014DRv2}.

\textbf{Adaptive Heavy-Tailed Simplex Optimization.} In Table~\ref{tab:AdaptiveHeavyTailSimplex}, we establish
an accelerated high-probability guarantee for smooth convex
optimization on the simplex under heavy-tailed gradient noise.
The updates require neither the noise moment order $p\in(1,2]$
nor the noise level $\sigma$, and the result permits arbitrary
constrained minimizers. The stopping budget still depends on
$p$ and $\sigma$; the improvement is primarily in adaptivity
and scope, rather than the accuracy exponent.

\textbf{These examples represent only a selection of our new results.} Many additional results and comparisons with prior work are available in Sections~\ref{sec:comparions} and~\ref{sec:add_comparisons} and on our \href{https://k3nfalt.github.io/longtails/}{website}.

\section{Our Prompting Strategy for Agentic Systems}
\label{sec:prompting}
Our research pipeline is relatively simple and could be potentially improved in terms of the strength of the final results and token efficiency. We run two agents: one is a researcher, and the other is a reviewer. We give them a prompt explaining that the researcher has to produce a result, while the reviewer has to verify it and decide whether to accept it. They interact until the reviewer accepts a new research result. See the prompt we use in our \href{https://github.com/k3nfalt/longtails/blob/main/scripts/prompts/research_prompt.md}{research prompt} and \href{https://github.com/k3nfalt/longtails/blob/main/scripts/run_research_prompt.py}{research run script}.

In general, the results come from two types of requests: i) we simply run the \href{https://github.com/k3nfalt/longtails/blob/main/scripts/run_research_prompt.py}{research run script} without additional instructions, which roughly asks the agents to produce ``any new result on a modern optimization topic that is interesting to the community''; ii) some results were obtained with additional requests. For instance, we were aware of a gap in the existing guarantees for $(L_0, L_1)$--smooth convex optimization and asked the agents to find a new method and a new Complexity with uniformly better guarantees. The agents obtained a better Complexity in a single run lasting approximately 2--7 hours.

\section{Related Work}
\label{sec:related-work}
\textbf{Formalized mathematics and optimization.}
Lean 4~\citep{deMouraUllrich2021Lean4} and the community-developed mathematical library mathlib~\citep{mathlib2020library} provide infrastructure for machine-checked mathematics. Within optimization, CvxLean~\citep{bentkamp2023verified} formalizes reductions between optimization problems, including transformations used in disciplined convex programming. Optlib~\citep{optlibRepository} and the work of Li et al.~\citep{li2024formalization} formalize convex analysis and convergence guarantees for first-order methods. These efforts establish that both optimization models and algorithmic guarantees can be verified in Lean. LongTails builds on this direction by organizing the search for new methods and improved oracle-complexity bounds around a fixed Setup that specifies admissible instances, available information, and the success criterion.

\textbf{Language models for formal theorem proving.}
LeanDojo~\citep{yang2023leandojo} provides tools for interacting with Lean and supports retrieval-augmented proof search. DeepSeek-Prover-V2~\citep{ren2025deepseekproverv2} combines subgoal decomposition with reinforcement learning for formal theorem proving. Closer to constructive tasks in applied mathematics, the Construction--Verification framework and AMBER benchmark~\citep{yang2026constructionverification} require agents to construct explicit solutions before proving their correctness, including tasks in optimization and algorithm design. LongTails focuses on the additional research objective of improving a Complexity while preserving the Setup.

\textbf{Automated scientific and algorithmic discovery.}
The AI Scientist~\citep{lu2024aiscientist} automates a research workflow involving idea generation, experiments, and scientific writing. AlphaEvolve has been used to discover mathematical constructions through language-model-guided evolutionary search and automated evaluation; some applications combine this search with formal proof tools~\citep{georgiev2025discovery}. Recent work by Lehre~\citep{lehre2026agentic} is particularly closely related: a staged agent pipeline designs black-box optimization algorithms and proves expected-runtime bounds in Lean using potential functions and drift analysis. AutoOPT~\citep{kim2026autoopt} combines numerical design through branch-and-bound performance estimation, symbolic discovery using large language models, and Lean 4 verification to develop first-order optimization methods. LongTails shares the goal of discovering algorithms with machine-checked efficiency guarantees, while targeting oracle-complexity bounds across deterministic and stochastic numerical optimization and monotone operator problems through a common Setup--Complexity interface.

\textbf{Computer-assisted complexity analysis.}
PEPit~\citep{goujaud2022pepit} automates worst-case analysis of first-order methods by formulating performance estimation problems as semidefinite programs and solving them numerically. This provides a complementary route to identifying tight guarantees and guiding algorithm design. LongTails expresses the guarantee as a Lean proposition about an algorithm executed under an explicit oracle budget, allowing the assumptions and the resulting bound to be checked together. Detailed comparisons with the optimization literature for each Setup are given in Section~\ref{sec:comparions}.

Our contribution lies in combining a reusable formal specification, explicit oracle accounting, and an agentic research workflow that scales across a broad collection of numerical problems. This scalability is central to LongTails: the same framework supports the discovery of new guarantees across different Setups and the repeated improvement of Complexities within each Setup. Using this framework, we have already obtained more than \numbernewresuls{} new theoretical guarantees in optimization. We continue to expand the collection of Setups, improve existing Complexities, and develop the framework.

\section{Conclusions}
In this work, we present LongTails, which, to the best of our knowledge, is the first framework to scale automated research in optimization across a broad range of problem classes while formally verifying the resulting guarantees in Lean 4. We also present at least \numbernewresuls{} new theoretical results across different areas of optimization, improving oracle Complexities in the stated regimes, broadening the assumptions under which guarantees hold, or establishing new guarantees for existing methods.

One limitation of this work is that we do not present each result in a separate paper. Although many of these results arguably warrant individual publications in the optimization community, preparing and reviewing on the order of \numbernewresuls{} separate papers is infeasible. As automated discovery scales, our framework and this paper may motivate discussion of a higher threshold for standalone publications or alternative ways of presenting, evaluating, and disseminating theoretical results.

\bibliographystyle{plain}
\bibliography{additional_citations,results_table}

\appendix
\newpage
\tableofcontents

\newpage
\clearpage
\section{Accelerated Method for Generalized Smooth Convex Optimization}
\label{sec:generalized-method}
We present the method discussed in Section~\ref{gen_inst}. Let $f$ be convex with a differentiable gradient satisfying $\|\nabla^2 f(x)\|\le L_0+L_1\|\nabla f(x)\|$ for all $x$. Assume that $L_0>0$, $L_1\ge0$, and that a minimizer $x^\star$ exists. The inputs $R,\Delta>0$ satisfy $\|x^0-x^\star\|\le R$ and $f(x^0)-f(x^\star)\le\Delta$. The method uses only exact gradients.

\begin{algorithm}[H]
\caption{Damped accelerated gradient method under $(L_0,L_1)$-smoothness}
\label{alg:damped-generalized}
\begin{algorithmic}[1]
\STATE \textbf{Input:} starting point $x^0$, parameters $L_0,L_1,R,\Delta$, accuracy $\varepsilon>0$
\IF{$\Delta\le\varepsilon$}
  \STATE \textbf{return} $x^0$
\ENDIF
\STATE Set $K$ according to~\eqref{eq:generalized-iterations}
\STATE Init $z^0=x^0$, $A_0=R^2/(2\Delta)$, $g^0=\nabla f(x^0)$, $G_0=4(L_0+L_1\|g^0\|)$
\FOR{$k=0,\ldots,K-1$}
  \STATE $a_k=\dfrac{1+\sqrt{1+4G_kA_k}}{2G_k}$, \quad $y^k=\dfrac{A_kx^k+a_kz^k}{A_k+a_k}$
  \STATE If $k>0$, compute $g^k=\nabla f(y^k)$; otherwise, reuse $g^0$ since $y^0=x^0$
  \STATE $M_k=4(L_0+L_1\|g^k\|)$, \quad $v^k=y^k-\dfrac{1}{M_k}g^k$
  \STATE Compute $h^k=\nabla f(v^k)$
  \IF{$M_k\le G_k$}
    \STATE $A_{k+1}=A_k+a_k$, \quad $x^{k+1}=v^k$
    \STATE $z^{k+1}=z^k-a_kh^k$, \quad $G_{k+1}=G_k/16$
  \ELSE
    \STATE $\gamma_k=G_k/M_k$, \quad $b_k=\gamma_ka_k$, \quad $A_{k+1}=A_k+b_k$
    \STATE $x^{k+1}=\dfrac{(1-\gamma_k)A_k}{A_{k+1}}x^k+
      \dfrac{\gamma_k(A_k+a_k)}{A_{k+1}}v^k$
    \STATE $z^{k+1}=z^k-b_kh^k$, \quad $G_{k+1}=16G_k$
  \ENDIF
\ENDFOR
\STATE \textbf{return} $x^K$
\end{algorithmic}
\end{algorithm}

The parameter $G_k>0$ is an estimate of the local curvature used to choose the acceleration weight $a_k$ before querying the gradient at $y^k$. We then compute $M_k=4(L_0+L_1\|\nabla f(y^k)\|)$ from this gradient. If $M_k>G_k$, the estimate is too small, so we damp the update by $\gamma_k=G_k/M_k$ and set $G_{k+1}=16G_k$. Otherwise, we take the full update and set $G_{k+1}=G_k/16$ to try a larger acceleration weight in the next iteration. Thus, each iteration adjusts the curvature estimate without an inner line search.

\textbf{Convergence guarantee.} Let $\ell=\log(4+4\Delta/\varepsilon)$ and choose
\begin{equation}
\label{eq:generalized-iterations}
K=\left\lceil 4096\left(1+\sqrt{\frac{L_0R^2}{\varepsilon}}
+\min\left\{\left\lceil(L_1R)^{2/3}\right\rceil\ell,\;
\ell+(L_1R\ell)^{2/3}\right\}\right)\right\rceil.
\end{equation}
Then Algorithm~\ref{alg:damped-generalized} returns $x^K$ satisfying $f(x^K)-f(x^\star)\le\varepsilon$ using $2K$ gradient queries. The algorithm and this guarantee are formalized in
\href{https://github.com/k3nfalt/longtails/blob/main/LongTails/Algorithms/DampedGeneralizedAccelerated.lean}{\texttt{DampedGeneralizedAccelerated.lean}} and
\href{https://github.com/k3nfalt/longtails/blob/main/LongTails/Theorems/GeneralizedSelfNormalizedConvex.lean}{\texttt{GeneralizedSelfNormalizedConvex.lean}}, respectively.

\clearpage
\subsection{Proof Idea}
\label{sec:generalized-proof-idea}
We asked an agent to explain the proof idea behind Algorithm~\ref{alg:damped-generalized} using the prompt ``Can you create a subsection drafting the proof idea?'' The sketch below was refined through follow-up requests.

\texttt{AI Generated Proof Idea From Here:} Let us explain the main steps of the proof. We consider $\Delta>\varepsilon$, since otherwise the initial point already satisfies the required accuracy.

\textbf{Potential decrease.} Consider
\begin{equation}
\label{eq:generalized-potential}
\mathcal{E}_k=A_k\bigl(f(x^k)-f(x^\star)\bigr)
+\frac12\|z^k-x^\star\|^2.
\end{equation}
The initialization gives $\mathcal{E}_0\le R^2$. Generalized smoothness ensures that
$\|h^k+M_k(v^k-y^k)\|\le(M_k/2)\|v^k-y^k\|$.
Writing $\theta_k=\min\{1,G_k/M_k\}$, the acceleration weight satisfies
$\theta_ka_k^2M_k\le G_ka_k^2=A_k+a_k$.
These two inequalities, together with convexity, imply
\begin{equation}
\label{eq:generalized-potential-decrease}
\mathcal{E}_{k+1}+d_k\le\mathcal{E}_k,
\qquad d_k=\theta_k(A_k+a_k)\frac{3M_k}{8}\|v^k-y^k\|^2.
\end{equation}
Thus, damping preserves potential decrease even when $G_k<M_k$. Summing~\eqref{eq:generalized-potential-decrease} gives $\sum_{k=0}^{K-1}d_k\le R^2$ and $f(x^K)-f(x^\star)\le R^2/A_K$.

\textbf{Controlling the curvature adjustments.} We group consecutive iterations into blocks: the first starts at iteration $0$, and each later iteration with $M_k>G_k$ starts a new block. A block therefore contains an increase of the curvature estimate followed by decreases, except for the first block, where $M_0=G_0$. We classify blocks by the curvature at their starting points and count iterations $1,\ldots,K-1$, leaving out the initial iteration. The updates of $G_k$ and the growth of $A_k$ bound the number of iterations in blocks starting with $M_k\le8L_0$ by $96\sqrt{L_0A_K}$.

For $L_1>0$, let $\mathcal{H}$ contain the starting points of the remaining blocks. At these iterations, $G_k\le M_k$ and $\|v^k-y^k\|\ge1/(8L_1)$. Hence~\eqref{eq:generalized-potential-decrease} gives
\begin{equation}
\label{eq:generalized-high-curvature}
\sum_{k\in\mathcal{H}}A_{k+1}G_k
\le\frac{512L_1^2}{3}\sum_{i=0}^{K-1}d_i
\le\frac{512}{3}(L_1R)^2.
\end{equation}
This limits the cost of blocks with high local curvature. The counting argument combines~\eqref{eq:generalized-high-curvature} with logarithmic growth of the acceleration weights through H\"older's inequality: the potential bound enters with power $1/3$, and the logarithmic term with power $2/3$. The resulting bound is
\begin{equation}
\label{eq:generalized-iteration-count}
K-1\le96\sqrt{L_0A_K}+4D_K+512(L_1R D_K)^{2/3},
\qquad D_K=\log\frac{2A_K}{A_0}.
\end{equation}
The detailed counting argument is formalized in
\href{https://github.com/k3nfalt/longtails/blob/main/LongTails/Proofs/GeneralizedLogCharging.lean}{\texttt{GeneralizedLogCharging.lean}} and
\href{https://github.com/k3nfalt/longtails/blob/main/LongTails/Proofs/GeneralizedSelfNormalizedConvergence.lean}{\texttt{GeneralizedSelfNormalizedConvergence.lean}}.

\textbf{Concluding the bound.} If $f(x^K)-f(x^\star)>\varepsilon$, then~\eqref{eq:generalized-potential}--\eqref{eq:generalized-potential-decrease} imply $A_K<R^2/\varepsilon$ and $D_K\le\ell$. Substituting these into~\eqref{eq:generalized-iteration-count} contradicts the choice of $K$ in~\eqref{eq:generalized-iterations}. Both branches of the minimum suffice: for $L_1>0$ and $\ell>1$,
$4\ell+512(L_1R\ell)^{2/3}\le516\lceil(L_1R)^{2/3}\rceil\ell$.
When $L_1=0$, all blocks have low curvature, and $K-1\le96\sqrt{L_0A_K}$ gives the same conclusion. Finally, each iteration uses two gradient queries.

\clearpage
\section{Comparisons for Results Discussed in the Main Text}
\label{sec:comparions}
In this section, we compare previous approaches with our newly discovered results discussed in Section~\ref{gen_inst}. All tables have the same structure:
\begin{enumerate}
  \item At the top of each table, we briefly introduce the setup.
  \item At the bottom, we define abbreviations and introduce notation for more compact formulas.
  \item In the \textbf{Result} column, we list our results and the previous methods relevant to the setup.
  \item In the \textbf{Oracle complexity} column, we list the oracle complexities of the methods.
  \item In the \textbf{Novelty} column, we explain how the new results improve upon the method in that row. For instance, ``(1): No $p,\sigma$ in updates; smoothness term $\varepsilon^{-1/2}$ up to logs instead of $\varepsilon^{-1}$; removes $G_0$'' in Table~\ref{tab:AdaptiveHeavyTailSimplex} means that our result \texttt{AdaptiveEntropicHeavyTailConvex}, labeled (1), does not require the noise moment order $p$ or noise level $\sigma$ in its updates, improves the accuracy dependence of the smoothness term from $\varepsilon^{-1}$ to $\varepsilon^{-1/2}$ up to logarithmic factors, and removes dependence on the initial gradient norm $G_0$ from the bound compared with result (2).
  \item In the \textbf{Limitations} column, we similarly explain the limitations of the new results compared with the method in that row. For instance, ``simplex rather than general Bregman geometry'' in Table~\ref{tab:AdaptiveHeavyTailSimplex} means that our result (1) is restricted to optimization on the simplex, whereas result (2) applies to more general domains and geometries.
\end{enumerate}
See also other new results on our \href{https://k3nfalt.github.io/longtails/}{website}. The following tables are AI-generated using the script \href{https://github.com/k3nfalt/longtails/blob/main/scripts/run_results_table_prompt.py}{\texttt{run\_results\_table\_prompt.py}}.

\let\comparisonTable\table
\def\table[#1]{\let\table\comparisonTable\comparisonTable[H]%
  \def\endtable{\end{adjustwidth}\csname float@endH\endcsname}}
\begingroup
\begin{table}[tp]
\centering
\caption{Heavy-tailed stochastic convex optimization on the simplex: \href{https://github.com/k3nfalt/longtails/blob/main/LongTails/Setups/AdaptiveHeavyTailSimplex.lean}{AdaptiveHeavyTailSimplex}.}\label{tab:AdaptiveHeavyTailSimplex}
\parbox{\linewidth}{\footnotesize On $\Delta_d=\{x\in\mathbb R^d:x_i\ge0,\ \sum_i x_i=1\}$, $d\ge2$, $f$ is convex, has an ambient gradient at every feasible point, and satisfies $\|\nabla f(x)-\nabla f(y)\|_\infty\le L\|x-y\|_1$ for feasible $x,y$, with $L>0$. A constrained minimizer $x_\star$ exists. The jointly measurable oracle uses IID samples of law $\nu$, with integrable replies and noise moments, $\int g(x,\omega)\,d\nu=\nabla f(x)$ and $\int\|g(x,\omega)-\nabla f(x)\|_\infty^p\,d\nu\le\sigma^p$, $1<p\le2$, $\sigma\ge0$. Target: a measurable feasible output with $P\{f(\mathrm{output})-f(x_\star)>\varepsilon\}\le\delta$, $\varepsilon>0$, $0<\delta<1$. Public inputs are $d,L,\delta$; the stopping budget also uses $p,\sigma,\varepsilon$.}\par\smallskip
\begin{threeparttable}
\footnotesize
\setlength{\tabcolsep}{4pt}\renewcommand{\arraystretch}{1.25}
\begin{tabularx}{\linewidth}{@{}>{\raggedright\arraybackslash}p{0.20\linewidth}>{\raggedright\arraybackslash}p{0.32\linewidth}>{\raggedright\arraybackslash}p{0.25\linewidth}>{\raggedright\arraybackslash}X@{}}
\toprule
\textbf{Result} & \textbf{Oracle complexity} & \textbf{Novelty} & \textbf{Limitations} \\
\midrule
\textbf{(1)} \href{https://github.com/k3nfalt/longtails/blob/main/LongTails/Theorems/AdaptiveEntropicHeavyTailConvex.lean}{Adaptive\allowbreak Entropic\allowbreak Heavy\allowbreak Tail\allowbreak Convex} \textbf{(Ours)} $(\star\star\star)$ \href{https://github.com/k3nfalt/longtails/blob/main/scripts/reviews/AdaptiveEntropicHeavyTailConvex.md}{(review)}
& $\displaystyle\mathcal{O}\!\left(\begin{gathered}\gamma+\sqrt{\frac{LH}{\varepsilon}\left[\gamma+\log\!\left(2+\frac{LH}{\varepsilon}\right)\right]}\\ {}+H\left(\frac{2^{38}\sigma}{\varepsilon}\right)^q\end{gathered}\right)$ & & \\
\midrule[1.2pt]
\textbf{(2)} Clipped-SMD \citep{nguyen2023heavyv2}, Thm.~3.1 (known horizon); entropic specialization with global smoothness and pilot
& $\displaystyle\mathcal{O}\!\left(\begin{gathered}\gamma+\frac{\gamma[LD+\sqrt D(G_0+\sigma)]}{\varepsilon}\\ {}+\gamma\left(\frac{\sigma\sqrt D}{\varepsilon}\right)^q\end{gathered}\right)$
& \textbf{(1)}: No $p,\sigma$ in updates; smoothness term $\varepsilon^{-1/2}$ up to logs instead of $\varepsilon^{-1}$; removes $G_0$. For fixed $p<2,\delta$, stochastic dimension factor $D^{q/2}\to D$.
& \textbf{(1)}: Stationary IID rather than conditional moments; simplex rather than general Bregman geometry. Large $2^{38q}$ penalty: dimensional gain requires roughly $D^{q/2-1}\gg2^{38q}$. \\
\midrule
\textbf{(3)} Clipped-ASMD \citep{nguyen2023heavyv2}, Thm.~4.1; globally smooth special case $\nabla f(x_\star)=0$
& $\displaystyle\mathcal{O}\!\left(1+\gamma\sqrt{\frac{LD}{\varepsilon}}+\gamma\left(\frac{\sigma\sqrt D}{\varepsilon}\right)^q\right)$
& \textbf{(1)}: Scope: arbitrary constrained minimizers. Updates need neither $p$ nor $\sigma$. On the overlap, deterministic term improves when $H[\gamma+\log(2+LH/\varepsilon)]\ll\gamma^2D$; no improved accuracy exponent.
& \textbf{(1)}: At fixed $d,\delta$, deterministic term loses $\sqrt{\log(L/\varepsilon)}$ as $L/\varepsilon\to\infty$; larger noise constant. Stationary IID and simplex restriction; stopping still needs $p,\sigma,\varepsilon$. \\

\bottomrule
\end{tabularx}
\begin{tablenotes}[flushleft]\footnotesize
\item[] \textbf{Abbr:} $q=p/(p-1)$, $D=\log d$, $\gamma=1+\log(1/\delta)$, $H=D+\gamma$; $x_0=(1/d,\ldots,1/d)$ and $G_0=\|\nabla f(x_0)\|_\infty$. Costs count individual stochastic-gradient replies, including the pilot in \textbf{(2)}; prior entropic bounds use entropy radius at most $\log d$. In \textbf{(4)}, $L_2$ and $\sigma_2$ are Euclidean smoothness and noise-moment constants, $R\ge\|x_0-x_\star\|_2$, $\zeta_\star\ge\|\nabla f(x_\star)\|_2$, and $\ell=1+\log\!\bigl((2+\sqrt{L_2R^2/\varepsilon}+\sqrt{R\zeta_\star/\varepsilon}+(\sigma_2R/\varepsilon)^q)/\delta\bigr)$. Prior-work big-$\mathcal{O}$ constants may depend on fixed $p$; the explicit $2^{38q}$ in \textbf{(1)} is retained because it is nonuniform as $p\downarrow1$.
\end{tablenotes}
\end{threeparttable}
\end{table}
\clearpage
\begin{table}[tp]
\centering
\caption{AdaptiveHeavyTailSimplex (continued).}
\begin{threeparttable}
\footnotesize
\setlength{\tabcolsep}{4pt}\renewcommand{\arraystretch}{1.25}
\begin{tabularx}{\linewidth}{@{}>{\raggedright\arraybackslash}p{0.20\linewidth}>{\raggedright\arraybackslash}p{0.32\linewidth}>{\raggedright\arraybackslash}p{0.25\linewidth}>{\raggedright\arraybackslash}X@{}}
\toprule
\textbf{Result} & \textbf{Oracle complexity} & \textbf{Novelty} & \textbf{Limitations} \\
\midrule
\textbf{(4)} Prox-clipped-SSTM-shift \citep{gorbunov2024heavyv2}, Thm.~2.7, $n=1$; Euclidean regularity on $B_{5\sqrt2 R}(x_\star)$
& $\displaystyle\mathcal{O}\!\left(\begin{gathered}\left(1+\sqrt{\frac{L_2R^2}{\varepsilon}}\right)\ell^2\\ {}+\sqrt{\frac{R\zeta_\star}{\varepsilon}}\,\ell^{3/2}+\left(\frac{\sigma_2R}{\varepsilon}\right)^q\ell\end{gathered}\right)$
& \textbf{(1)}: Updates need neither $p$ nor $\sigma$. Rates not comparable on the full setup. On a common regular extension, worst-case $L_2=dL$, $\sigma_2=\sqrt d\sigma$, $R\le\sqrt2$ give polynomial rather than logarithmic dimension factors.
& \textbf{(1)}: Simplex only; \textbf{(4)} covers composite objectives and unbounded domains. Conversion gains need not hold when $L_2=L$, $\sigma_2=\sigma$; \textbf{(1)} retains its noise factor $H$. \\
\midrule
\textbf{(5)} ASMD \citep{liu2026expectationv1}, Thm.~3.3, known-$p,\sigma,T$ tuning; Derived bound (Markov)
& $\displaystyle\mathcal{O}\!\left(\begin{gathered}1+\delta^{-1/p}\sqrt{\frac{LD}{\varepsilon}}\\ {}+\delta^{-2/(p-1)}\left(\frac{\sigma\sqrt D}{\varepsilon}\right)^q\end{gathered}\right)$
& \textbf{(1)}: Logarithmic instead of polynomial confidence cost for the shared probability target; updates need neither $p$ nor $\sigma$. \textbf{(5)} derives from a bound on $\mathbb E[\mathrm{gap}^{p/2}]$.
& \textbf{(1)}: Requires an ambient gradient and the simplex; \textbf{(5)} allows general composite Bregman problems. \\
\bottomrule
\end{tabularx}
\begin{tablenotes}[flushleft]\footnotesize
\item[] \textbf{Abbr:} $q=p/(p-1)$, $D=\log d$, $\gamma=1+\log(1/\delta)$, $H=D+\gamma$; $x_0=(1/d,\ldots,1/d)$ and $G_0=\|\nabla f(x_0)\|_\infty$. Costs count individual stochastic-gradient replies, including the pilot in \textbf{(2)}; prior entropic bounds use entropy radius at most $\log d$. In \textbf{(4)}, $L_2$ and $\sigma_2$ are Euclidean smoothness and noise-moment constants, $R\ge\|x_0-x_\star\|_2$, $\zeta_\star\ge\|\nabla f(x_\star)\|_2$, and $\ell=1+\log\!\bigl((2+\sqrt{L_2R^2/\varepsilon}+\sqrt{R\zeta_\star/\varepsilon}+(\sigma_2R/\varepsilon)^q)/\delta\bigr)$. Prior-work big-$\mathcal{O}$ constants may depend on fixed $p$; the explicit $2^{38q}$ in \textbf{(1)} is retained because it is nonuniform as $p\downarrow1$.
\end{tablenotes}
\end{threeparttable}
\end{table}
\clearpage

\begin{table}[tp]
\centering
\caption{Nonconvex optimization with auxiliary similarity: \href{https://github.com/k3nfalt/longtails/blob/main/LongTails/Setups/AuxiliarySimilarityNonconvex.lean}{AuxiliarySimilarityNonconvex}.}\label{tab:AuxiliarySimilarityNonconvex}
\parbox{\linewidth}{\footnotesize On a finite-dimensional real inner-product space, $f\in C^2$, $h$ is differentiable, $\nabla^2 f$ is $\rho$-Lipschitz, $\nabla h$ is $L$-Lipschitz, and $\nabla(f-h)$ is $\delta$-Lipschitz. Assume $f(x_0)-f(x)\leq\Delta$ for every $x$, $L,\Delta\geq0$, and $\delta,\rho,\varepsilon>0$. Exact oracles return $(f,\nabla f)$ and $(h,\nabla h)$. The target is a returned point $\hat x$ with $\|\nabla f(\hat x)\|\leq\varepsilon$ almost surely; \textbf{(1)} is deterministic.}\par\smallskip
\begin{threeparttable}
\footnotesize
\setlength{\tabcolsep}{4pt}\renewcommand{\arraystretch}{1.25}
\begin{tabularx}{\linewidth}{@{}>{\raggedright\arraybackslash}p{0.20\linewidth}>{\raggedright\arraybackslash}p{0.32\linewidth}>{\raggedright\arraybackslash}p{0.25\linewidth}>{\raggedright\arraybackslash}X@{}}
\toprule
\textbf{Result} & \textbf{Oracle complexity} & \textbf{Novelty} & \textbf{Limitations} \\
\midrule
\textbf{(1)} \href{https://github.com/k3nfalt/longtails/blob/main/LongTails/Theorems/AuxiliarySimilarityAcceleration.lean}{Auxiliary\allowbreak{}Similarity\allowbreak{}Acceleration} \textbf{(Ours)} $(\star\star\star)$ \href{https://github.com/k3nfalt/longtails/blob/main/scripts/reviews/AuxiliarySimilarityAcceleration.md}{(review)} &
$N_f=\mathcal{O}((1+A)\sqrt{1+\delta/s}\,\ell_\delta\ell_L)$;\newline
$N_h=\mathcal{O}(N_f(1+L/\delta)$\newline $\qquad\cdot[\ell_L+\log(1+A)])$ & & \\
\midrule[1.2pt]
\textbf{(2)} AuxMOM \citep{asnChaytiV5}, Sec.~4.2, noiseless &
$N_f=\mathcal{O}(1+\delta\Delta/\varepsilon^2)$;\newline
$N_h=\mathcal{O}((1+L/\delta)N_f)$.
Random output: $\mathbb E\|\nabla f(\hat x)\|^2\leq\varepsilon^2$. &
\textbf{(1)}: fewer $f$ calls when $A\geq1$ and $\delta/s\gg(\ell_\delta\ell_L)^2$; gain $\sqrt{\delta/s}/(\ell_\delta\ell_L)$. For $\rho=\Delta=\delta=1$, $L=\varepsilon^{-1/2}$: $\mathcal{O}(\varepsilon^{-7/4}\log^2(1/\varepsilon))$ versus $\mathcal{O}(\varepsilon^{-2})$ $f$ calls. &
\textbf{(1)} needs Hessian Lipschitzness and $\rho,\Delta$; extra helper logarithms and larger constants. Worse $f$ cost when $\delta\lesssim s$. \\
\midrule
\textbf{(3)} Restarted AGD \citep{asnLiLin2023}, Thm.~1 &
$N_f=\mathcal{O}(1+\sqrt{L+\delta}\,\rho^{1/4}\Delta\varepsilon^{-7/4})$, $N_h=0$; $\varepsilon\lesssim(L+\delta)^2/\rho$. &
\textbf{(1)}: for $A\geq1$, $\delta\geq s$, replaces $\sqrt{L+\delta}$ by $\sqrt\delta\,\ell_\delta\ell_L$; gains when $L\gg\delta(\ell_\delta\ell_L)^2$. The example in \textbf{(2)} also improves \textbf{(3)}'s $\mathcal{O}(\varepsilon^{-2})$ $f$ calls. &
\textbf{(1)} needs helper calls and $\Delta$; loses logarithmic factors when $L\lesssim\delta$. \\
\bottomrule
\end{tabularx}
\begin{tablenotes}[flushleft]\footnotesize
\item[] \textbf{Abbr:} $N_f,N_h$: respective total oracle calls (gradient calls for prior work); $s=\sqrt{\rho\varepsilon}$, $A=\Delta\sqrt\rho\,\varepsilon^{-3/2}$, $\ell_\delta=\log(2+\delta/s)$, $\ell_L=\log(2+(L+\delta)/s)$. In the $N_h$ formulas, $N_f$ denotes its displayed upper-bound expression. $\lesssim$ hides universal constants; $\gg$ denotes asymptotic separation.
\end{tablenotes}
\end{threeparttable}
\end{table}
\clearpage

\begin{table}[tp]
\centering
\caption{Comparison-oracle stationarity under \href{https://github.com/k3nfalt/longtails/blob/main/LongTails/Setups/ComparisonHessianLipschitz.lean}{ComparisonHessianLipschitz}.}\label{tab:ComparisonHessianLipschitz}
\parbox{\linewidth}{\footnotesize On a finite-dimensional real inner-product space, $f\in C^2$ has $L$-Lipschitz gradient and $\rho$-Lipschitz Hessian, with $f(x_0)-f(x)\leq\Delta$ for every $x$; $L,\rho,\varepsilon>0$, $\Delta\geq0$. Queries return exact bits $[f(x)\leq f(y)]$. The target is a list containing $x$ with $\|\nabla f(x)\|\leq\varepsilon$ on every run. The setup restricts list entries to $x_0$ or queried points, not to a particular method.}\par\smallskip
\begin{threeparttable}
\footnotesize
\setlength{\tabcolsep}{4pt}\renewcommand{\arraystretch}{1.25}
\begin{tabularx}{\linewidth}{@{}>{\raggedright\arraybackslash}p{0.20\linewidth}>{\raggedright\arraybackslash}p{0.32\linewidth}>{\raggedright\arraybackslash}p{0.25\linewidth}>{\raggedright\arraybackslash}X@{}}
\toprule
\textbf{Result} & \textbf{Oracle complexity} & \textbf{Novelty} & \textbf{Limitations} \\
\midrule
\textbf{(1)} \href{https://github.com/k3nfalt/longtails/blob/main/LongTails/Theorems/ComparisonTrustAccelerated.lean}{Comparison\allowbreak Trust\allowbreak Accelerated} \textbf{(Ours)} $(\star\star\star)$ \href{https://github.com/k3nfalt/longtails/blob/main/scripts/reviews/ComparisonTrustAccelerated.md}{(review)} & $\mathcal{O}\!\left((1+\Delta\sqrt\rho\,\varepsilon^{-3/2})\right.$\newline $\left.{}\cdot(n+1)(1+q)\ell^3\right)$ comparisons & & \\
\midrule[1.2pt]
\textbf{(2)} \citep[Thm.~8, v1]{chltWang2026} & $\mathcal{O}\!\left((1+\Delta\sqrt\rho\,\varepsilon^{-3/2})\right.$\newline $\left.{}\cdot n^2\operatorname{polylog}(n,L,\rho,\varepsilon^{-1})\right)$ comparisons; success $2/3$. & \textbf{(1)}: leading cost saves $n/(1+q)$ up to logs when $n\gg(1+q)\ell^3$; deterministic instead of success $2/3$. & \textbf{(1)} loses second-order accuracy $\lambda_{\min}(\nabla^2 f)\geq-\sqrt{\rho\varepsilon}$; polynomial $L^{1/2}$ dependence, worse when $n\lesssim q$. \\
\midrule
\textbf{(3)} \citep[Thm.~4, v1]{csZhangLi2024v1} & $\mathcal{O}\!\left((1+L\Delta/\varepsilon^2)n\log(2+n)\right)$ comparisons; selected point succeeds with probability $2/3$. & \textbf{(1)}: leading cost saves $q$ up to logs for $\varepsilon\ll L^2/\rho$. At $L=\rho=\Delta=1$, $n=\varepsilon^{-1/2}$: $\mathcal{O}(\varepsilon^{-9/4}\log^3(1/\varepsilon))$ versus $\mathcal{O}(\varepsilon^{-5/2}\log(1/\varepsilon))$. & \textbf{(1)} needs Lipschitz Hessian and guarantees a containing list, not a selected stationary point. \\
\bottomrule
\end{tabularx}
\begin{tablenotes}[flushleft]\footnotesize
\item[] \textbf{Abbr:} $n=\dim E$, $q=(L^2/(\rho\varepsilon))^{1/4}$, $\ell=\log(2+n+L/\sqrt{\rho\varepsilon}+\sqrt{L\Delta}/\varepsilon)$. $\operatorname{polylog}$ denotes the logarithmic factors suppressed in the cited bound. Leading-cost comparisons assume $\Delta\sqrt\rho\,\varepsilon^{-3/2}\gtrsim1$; the numerical regime is as $\varepsilon\downarrow0$.
\end{tablenotes}
\end{threeparttable}
\end{table}
\clearpage

\begin{table}[tp]
\centering
\caption{Composite optimization with parameter-free updates: \href{https://github.com/k3nfalt/longtails/blob/main/LongTails/Setups/CompositeGeneralizedParameterFree.lean}{CompositeGeneralizedParameterFree}.}\label{tab:CompositeGeneralizedParameterFree}
\parbox{\linewidth}{\footnotesize Minimize $F=f+\psi$ over convex $Q$ in a real Hilbert space. The function $f$ is convex with continuous gradient and satisfies $\|\nabla f(x)-\nabla f(y)\|\le (L_0+L_1M)\|x-y\|$ whenever $\|\nabla f(z)\|\le M$ on $[x,y]$; $\psi$ is finite-valued and convex on $Q$. Exact gradients and minimizers of $\psi(u)+\langle v,u-z\rangle+\|u-z\|^2/(2a)$ over $u\in Q$ are available for every $a>0,z,v$. There is a minimizer $x_\star\in Q$ with $x_0\in Q$, $\|x_0-x_\star\|\le R$, $F(x_0)-F(x_\star)\le\Delta$, and $\|\nabla f(x_\star)\|\le G$. Here $L_0,R,\Delta,h_0,\varepsilon>0$ and $L_1,G\ge0$. Target: deterministic $\widehat x\in Q$ with $F(\widehat x)-F(x_\star)\le\varepsilon$. The constructed method receives only $(x_0,h_0)$; its stopping budget depends on all setup parameters.}\par\smallskip
\begin{threeparttable}
\footnotesize
\setlength{\tabcolsep}{4pt}\renewcommand{\arraystretch}{1.25}
\begin{tabularx}{\linewidth}{@{}>{\raggedright\arraybackslash}p{0.18\linewidth}>{\raggedright\arraybackslash}p{0.37\linewidth}>{\raggedright\arraybackslash}p{0.24\linewidth}>{\raggedright\arraybackslash}X@{}}
\toprule
\textbf{Result} & \textbf{Oracle complexity} & \textbf{Novelty} & \textbf{Limitations} \\
\midrule
\textbf{(1)} \href{https://github.com/k3nfalt/longtails/blob/main/LongTails/Theorems/BacktrackingCompositeGeneralized.lean}{Backtracking\allowbreak Composite\allowbreak Generalized} \textbf{(Ours)} $(\star\star\star)$ \href{https://github.com/k3nfalt/longtails/blob/main/scripts/reviews/BacktrackingCompositeGeneralized.md}{(review)}
& $\displaystyle\mathcal{O}\!\left(\begin{aligned}
&\sqrt{L_\star R^2/\varepsilon}\\
&+(1+s)\log\!\left(2+s+\frac{1}{L_\star h_0}+\frac{L_1^2\Delta}{L_\star}\right)\\
&+\log(2+L_\star h_0+L_1^2\Delta h_0)
\end{aligned}\right)$\par Gradient and prox calls; all rejected trials included.
& & \\
\midrule[1.2pt]
\textbf{(2)} Accelerated proximal gradient \citep[Thm.~B.1]{LiLaiAnHuWen2026v3Composite}\par Twice-differentiable $f$ on $\mathbb R^d$, $\|\nabla^2f\|\le L_0+L_1\|\nabla f\|$; proper closed convex proximable $h=\psi+\iota_Q$.
& $\displaystyle\mathcal{O}\!\left(\begin{aligned}
&\sqrt{L_\star R^2/\varepsilon}\\
&+\log\!\left(2+\frac{\Delta}{L_\star R^2}\right)\\
&+(1+s)\log(2+\Delta/\varepsilon)
\end{aligned}\right)$\par Gradient and prox calls; $\Gamma_0=2\Delta/R^2$.
& \textbf{(1)}: Updates need no $L_0,L_1,G,R,\Delta$; matched leading rate. Hilbert-space, first-order smoothness scope. With $h_0\asymp L_\star^{-1}$, logarithmic overhead improves for $L_1>0$, $\varepsilon\ll\min\{L_\star/L_1^2,\Delta/(1+s)\}$.
& \textbf{(1)}: Roughly threefold call overhead per accepted step; $h_0$-mismatch logs. For $h_0\asymp L_\star^{-1}$, $\varepsilon\gtrsim L_\star/L_1^2$ can lose a logarithm. Finite-valued $\psi$ rather than general extended-valued $h$. \\
\midrule
\textbf{(3)} Accelerated GRAAL \citep[Cor.~3]{BorodichKovalev2025v2Composite}\par Special case: $Q=\mathbb R^d$, $\psi=0$, $G=0$; $C^2$ Hessian condition as in \textbf{(2)}.
& $\displaystyle\mathcal{O}\!\left(\begin{aligned}
&1+\sqrt{L_0R^2/\varepsilon}+s^3\\
&+(1+s^2)\log\!\frac{1}{\eta_0L_0}
\end{aligned}\right)$\par Gradient and value calls; requires $\eta_0L_0e^s\le1$.
& \textbf{(1)}: Composite scope; no values (prox is closed form here); arbitrary $h_0$. On the overlap, $s\gg1$, $h_0\asymp L_0^{-1}$, $\eta_0=e^{-s}/L_0$, $L_1^2\Delta/L_0\lesssim e^s$: $\mathcal{O}(1+\sqrt{L_0R^2/\varepsilon}+s^2)$ versus $\mathcal{O}(1+\sqrt{L_0R^2/\varepsilon}+s^3)$; total gain when $\varepsilon\gg L_0R^2/s^6$.
& \textbf{(1)}: Constant-factor backtracking overhead. \\
\midrule
\textbf{(4)} AC-FGM \citep[Cor.~1 and subsequent remark]{LiLan2024v3Composite}\par Initial-search variant; closed convex $Q\subseteq\mathbb R^d$, closed convex $\psi$; local Lipschitz bound $\bar L$ covering search points and trials.
& $\displaystyle\mathcal{O}\!\left(1+\sqrt{\bar L R^2/\varepsilon}+\log_+\!\frac{\bar L}{\ell}\right)$\par Value, gradient and prox calls, including initialization; initial curvature estimate $0<\ell\le\bar L$.
& \textbf{(1)}: No function values; explicit polynomial parameter dependence. Rates are not comparable: no polynomial bound on $\bar L$ in $L_0,L_1,G,R,\Delta$ is established here.
& \textbf{(1)}: Repeated backtracking instead of initial search only. \\
\bottomrule
\end{tabularx}
\begin{tablenotes}[flushleft]\footnotesize
\item[] \textbf{Abbr:} $L_\star=L_0+L_1G$; $s=L_1R$; $\log_+t=\max\{0,\log t\}$. $\eta_0$ is the initial step in \textbf{(3)}; $\bar L$ bounds the local Lipschitz estimates and initial trials in \textbf{(4)}, and $\ell$ is its initial curvature estimate. $\iota_Q$ is $0$ on $Q$ and $+\infty$ outside. Costs count individual exact oracle calls; $\asymp$ and $\lesssim$ suppress universal factors.
\end{tablenotes}
\end{threeparttable}
\end{table}
\clearpage

\begin{table}[tp]
\centering
\caption{Accelerated composite optimization under generalized smoothness: \href{https://github.com/k3nfalt/longtails/blob/main/LongTails/Setups/CompositeGeneralizedSmooth.lean}{CompositeGeneralizedSmooth}.}\label{tab:CompositeGeneralizedSmooth}
\parbox{\linewidth}{\footnotesize Minimize $F=f+\psi$ over convex $Q$ in a real Hilbert space. The function $f$ is convex and continuously differentiable on the whole space, with $\|\nabla f(x)-\nabla f(y)\|\le (L_0+L_1M)\|x-y\|$ whenever $\|\nabla f\|\le M$ on $[x,y]$; $\psi$ is convex on $Q$, and exact gradient and proximal queries are available. Here $L_0,R,\Delta,\varepsilon>0$, $L_1,G\ge0$, $x_0\in Q$, and some minimizer $x_\star\in Q$ satisfies $\|x_0-x_\star\|\le R$, $F(x_0)-F(x_\star)\le\Delta$, and $\|\nabla f(x_\star)\|\le G$. Target: a feasible output $x$ with $F(x)-F(x_\star)\le\varepsilon$, deterministically.}\par\smallskip
\begin{threeparttable}
\footnotesize
\setlength{\tabcolsep}{4pt}\renewcommand{\arraystretch}{1.25}
\begin{tabularx}{\linewidth}{@{}>{\raggedright\arraybackslash}p{0.18\linewidth}>{\raggedright\arraybackslash}p{0.37\linewidth}>{\raggedright\arraybackslash}p{0.24\linewidth}>{\raggedright\arraybackslash}X@{}}
\toprule
\textbf{Result} & \textbf{Oracle complexity} & \textbf{Novelty} & \textbf{Limitations} \\
\midrule
\textbf{(1)} \href{https://github.com/k3nfalt/longtails/blob/main/LongTails/Theorems/CompositeGeneralizedAccelerated.lean}{Composite\allowbreak Generalized\allowbreak Accelerated} \textbf{(Ours)} $(\star\star\star)$ \href{https://github.com/k3nfalt/longtails/blob/main/scripts/reviews/CompositeGeneralizedAccelerated.md}{(review)}
& $\displaystyle\mathcal{O}\!\left(\begin{gathered}\sqrt{\frac{(L_0+L_1G)R^2}{\varepsilon}}\\ {}+(1+L_1R)\log\!\left(2+(1+L_1R)\sqrt{\frac{\Delta}{L_0R^2}}\right)\end{gathered}\right)$ gradient/prox calls ($3K$ total). & & \\
\midrule[1.2pt]
\textbf{(2)} \citep{cgsTyurin2026v2}, Thm.~4.3. Special case: $Q=\mathbb{R}^d$, $\psi=0$, $G=0$; twice differentiable $f$; $R=\|x_0-x_\star\|$.
& $\displaystyle\mathcal{O}\!\left(\begin{gathered}1+\sqrt{\frac{L_0R^2}{\varepsilon}}\\ {}+L_1R\log_+\!\left(\min\!\left\{\frac{L_1^2\Delta}{L_0},\frac{\Delta}{\varepsilon}\right\}\right)\\ {}+\log_+\!\frac{\Delta}{L_0R^2}\end{gathered}\right)$ gradient calls; $\bar R=R$, $\Gamma_0=2\Delta/R^2$.
& \textbf{(1)}: Scope extends acceleration to constrained/composite objectives with $\nabla f(x_\star)\ne0$. Matches the leading rate on the overlap, up to logarithmic terms.
& \textbf{(1)}: Larger logarithmic term can occur when $L_1>0$, $\varepsilon>L_0/L_1^2$; worse leading constant on the overlap. \\
\midrule
\textbf{(3)} AC-FGM \citep{cgsLiLan2024v3}, Cor.~1 and local-smoothness remark. Finite-dimensional, closed $Q$ and closed convex $\psi$.
& $\displaystyle\mathcal{O}\!\left(1+\sqrt{\frac{\bar L(R^2+d_1^2)}{\varepsilon}}\right)$ gradient, value and prox calls in total; arbitrary initial $\eta_1>0$, fixed $\beta=1/10$.
& \textbf{(1)}: No function values; explicit accelerated composite bound in setup parameters. Rates are not comparable through an unspecified trajectory curvature $\bar L$.
& \textbf{(1)}: Needs $L_0,L_1,R,\Delta$ for updates and $G$ for its stopping budget; \textbf{(3)} has parameter-free updates. \\
\bottomrule
\end{tabularx}
\begin{tablenotes}[flushleft]\footnotesize
\item[] \textbf{Abbr:} Each gradient, value or prox evaluation counts as one call. A prox query minimizes $\|u-z\|^2/(2a)+\langle v,u-z\rangle+\psi(u)$ over $u\in Q$, $a>0$. $\log_+t=\max\{0,\log t\}$ for $t>0$, with $\log_+0=0$. In \textbf{(3)}, $z_0=x_0$, $d_1=\|z_1-z_0\|$, and $\bar L=\sup_{t\ge1}\max\{[4(1-\beta)\eta_1]^{-1},\ell_1,\ldots,\ell_t\}$; $\ell_t$ are the curvature estimates in Eqs.~(2.12)--(2.13), distinct from $L_0,L_1$.
\end{tablenotes}
\end{threeparttable}
\end{table}
\clearpage

\begin{table}[tp]
\centering
\caption{Constrained generalized monotone variational inequalities: \href{https://github.com/k3nfalt/longtails/blob/main/LongTails/Setups/ConstrainedGeneralizedMonotone.lean}{ConstrainedGeneralizedMonotone}.}\label{tab:ConstrainedGeneralizedMonotone}
\parbox{\linewidth}{\footnotesize In a real inner-product space $E$, $X$ is convex with an exact projection oracle. The continuous, monotone operator $F:E\to E$ satisfies $\|F(x)-F(y)\|\le (L_0+L_1\sup_{z\in[x,y]}\|F(z)\|)\|x-y\|$ for all $x,y\in E$. Here $L_0,R,\varepsilon>0$, $L_1,G,G_0\ge0$, $x_0\in X$, and there is a VI solution $x^\star\in X$ with $\|x_0-x^\star\|\le R$, $\|F(x^\star)\|\le G$, and $\|F(x_0)\|\le G_0$. Target: $x\in X$ and some $\xi\in N_X(x)$ with $\|F(x)+\xi\|\le\varepsilon$, using exact $F$ and projection calls.}\par\smallskip
\begin{threeparttable}
\footnotesize
\setlength{\tabcolsep}{4pt}\renewcommand{\arraystretch}{1.25}
\begin{tabularx}{\linewidth}{@{}>{\raggedright\arraybackslash}p{0.18\linewidth}>{\raggedright\arraybackslash}p{0.37\linewidth}>{\raggedright\arraybackslash}p{0.24\linewidth}>{\raggedright\arraybackslash}X@{}}
\toprule
\textbf{Result} & \textbf{Oracle complexity} & \textbf{Novelty} & \textbf{Limitations} \\
\midrule
\textbf{(1)} \href{https://github.com/k3nfalt/longtails/blob/main/LongTails/Theorems/ConstrainedGeneralizedMonotone.lean}{Constrained\allowbreak Generalized\allowbreak Monotone} \textbf{(Ours)} $(\star\star\star)$ \href{https://github.com/k3nfalt/longtails/blob/main/scripts/reviews/ConstrainedGeneralizedMonotone.md}{(review)}
& $\mathcal{O}\!\left(\begin{gathered}\dfrac{(L_0+L_1G)(1+L_1R)R}{\varepsilon}\\{}+(1+L_1R)^2\ell\end{gathered}\right)$ total calls; two $F$ calls and one projection per iteration. & & \\
\midrule[1.2pt]
\textbf{(2)} PF-NE-EG \citep{cgmShen2026v2}, App.~A, Prop.~2, Lemma~10. Locally Lipschitz $F$ on closed convex $X\subseteq\mathbb{R}^d$.
& Derived bound: $\mathcal{O}\!\left(1+\dfrac{R^2}{\eta_*^2\varepsilon^2}+\log\dfrac{\eta_0}{\eta_*}\right)$ total $F$/projection calls, including rejected trials; implicit $\eta_*>0$, fixed backtracking factors.
& \textbf{(1)}: $\varepsilon^{-1}$ instead of $\varepsilon^{-2}$ as $\varepsilon\to0$ with other parameters fixed; explicit polynomial growth constants, with logarithmic $G_0$ dependence.
& \textbf{(1)} needs known $L_0,L_1$ and global generalized smoothness and monotonicity; \textbf{(2)} needs only local Lipschitz continuity and monotonicity on $X$. \\
\midrule
\textbf{(3)} Composite-EAG \citep{cgmCai2024v6}, Thm.~3.2. Special case: globally $L$-Lipschitz $F$ on $\mathbb{R}^d$, $A=N_X$.
& $\mathcal{O}(1+LR/\varepsilon)$ total $F$/projection calls.
& \textbf{(1)}: Scope extends to operators without a finite global Lipschitz constant. Matches known rate when $L_1=0$ and $L=L_0$.
& \textbf{(1)} can lose factors $1+L_1R$ and logarithmic overhead for nearly Lipschitz operators; does not cover general proximal inclusions or negative comonotonicity covered by \textbf{(3)}. \\
\midrule
\textbf{(4)} EG \citep{cgmChoudhury2025v1}, Thm.~3.6. Unconstrained special case: $X=E=\mathbb{R}^d$, $G=0$.
& $\mathcal{O}\!\left(1+\dfrac{L_0^2R^2}{\varepsilon^2}+L_1^2R^2\right)$ $F$ calls.
& \textbf{(1)}: $\varepsilon^{-1}$ instead of $\varepsilon^{-2}$ for fixed remaining parameters and small $\varepsilon$; last instead of best iterate. Scope extends to constraints and $F(x^\star)\ne0$.
& \textbf{(1)} can lose a logarithmic factor when the $L_1^2R^2$ term dominates \textbf{(4)}. \\
\midrule
\textbf{(5)} Repository anchored extragradient \citep{cgmRepositoryAnchored}. Unconstrained special case: $X=E$, $G=0$.
& $\mathcal{O}\!\left(\begin{gathered}1+\dfrac{L_0R}{\varepsilon}\\{}+(1+L_1R)\log\!\left(1+\dfrac{(L_0+L_1G_0)R}{\varepsilon}\right)\end{gathered}\right)$ $F$ calls.
& \textbf{(1)}: Scope extends to constrained VIs with $F(x^\star)\ne0$.
& \textbf{(1)} loses a factor $1+L_1R$ in the leading term on $X=E$, and has larger constants. \\
\bottomrule
\end{tabularx}
\begin{tablenotes}[flushleft]\footnotesize
\item[] \textbf{Abbr:} $N_X(x)=\{\xi:\langle\xi,w-x\rangle\le0\ \forall w\in X\}$; $[x,y]$ is the line segment; $\ell=1+\log(1+(L_0+L_1G_0)R/\varepsilon)$. $L$ is a global Lipschitz constant. $\eta_0$ is the initial step and $\eta_*=\inf_t\eta_t>0$ the accepted-step infimum for \textbf{(2)}. Total calls sum operator evaluations and projections; projection is the identity when $X=E$.
\end{tablenotes}
\end{threeparttable}
\end{table}
\clearpage

\begin{table}[tp]
\centering
\caption{Similar strongly quasar-convex finite sums: \href{https://github.com/k3nfalt/longtails/blob/main/LongTails/Setups/ContinuizedSimilarQuasar.lean}{ContinuizedSimilarQuasar}.}\label{tab:ContinuizedSimilarQuasar}
\parbox{\linewidth}{\footnotesize On a real Hilbert space, $F=n^{-1}\sum_{i=1}^n f_i$, $n\ge1$, has $L$-Lipschitz gradient; each $f_i$ has a gradient. A global minimizer $x_\star$ satisfies $F(x_0)-F(x_\star)\le\Delta$ and $F(x_\star)\ge F(x)+\gamma^{-1}\langle\nabla F(x),x_\star-x\rangle+\frac\mu2\|x-x_\star\|^2$. For all $x,y$, $n^{-1}\sum_i\|[\nabla f_i-\nabla F](x)-[\nabla f_i-\nabla F](y)\|^2\le\delta^2\|x-y\|^2$. Here $0<\mu\le L$, $0<\gamma\le1$, $\delta,\Delta\ge0$, and $\varepsilon>0$. The oracle returns exact component gradients; interpolation is not required. Target: $\int [F(X(\omega))-F(x_\star)]\,dP(\omega)\le\varepsilon$, with integrable gap and $P$ the execution law. Our constructions have fixed budgets and expected final-iterate accuracy; the setup imposes no algorithm restriction.}\par\smallskip
\begin{threeparttable}
\footnotesize
\setlength{\tabcolsep}{4pt}\renewcommand{\arraystretch}{1.25}
\begin{tabularx}{\linewidth}{@{}>{\raggedright\arraybackslash}p{0.20\linewidth}>{\raggedright\arraybackslash}p{0.32\linewidth}>{\raggedright\arraybackslash}p{0.25\linewidth}>{\raggedright\arraybackslash}X@{}}
\toprule
\textbf{Result} & \textbf{Oracle complexity} & \textbf{Novelty} & \textbf{Limitations} \\
\midrule
\textbf{(1)} \href{https://github.com/k3nfalt/longtails/blob/main/LongTails/Theorems/MovementSimilarQuasar.lean}{Movement\allowbreak Similar\allowbreak Quasar} \textbf{(Ours)} $(\star\star\star)$ \href{https://github.com/k3nfalt/longtails/blob/main/scripts/reviews/MovementSimilarQuasar.md}{(review)}
& $\mathcal{O}\!\left(n+\left[n+\frac{\sqrt\kappa}{\gamma}+\frac{(n\delta/\mu)^{2/3}}{\gamma^{4/3}}\right]\ell_\gamma\right)$ component gradients.
& & \\
\midrule
\textbf{(2)} \href{https://github.com/k3nfalt/longtails/blob/main/LongTails/Theorems/ContinuizedSimilarQuasar.lean}{Continuized\allowbreak Similar\allowbreak Quasar} \textbf{(Ours)} $(\star\star\star)$ \href{https://github.com/k3nfalt/longtails/blob/main/scripts/reviews/ContinuizedSimilarQuasar.md}{(review)}
& $\mathcal{O}\!\left(\left[n+\frac{\sqrt\kappa}{\gamma}+\frac{\delta^2}{\gamma^3\mu^2}\right](1+\ell_\gamma)\right)$ component gradients.
& & \\
\midrule[1.2pt]
\textbf{(3)} PAGE \citep{csqLi2021}, Thm.~5, Cor.~6. Derived bound on $\mathbb{R}^d$ via PL constant $\gamma^2\mu$.
& $\mathcal{O}\!\left(n+\left[n+\frac{\sqrt n\sqrt{L^2+\delta^2}}{\gamma^2\mu}\right]\ell_1\right)$ expected component gradients.
& \textbf{(1)}: no worse polynomial dependence; extra accuracy logarithm. For $\gamma=1$, $n=\kappa$, $\varepsilon\le\Delta$: $\delta=L$ gives $\mathcal{O}(\kappa^{4/3}\ell_1)$ instead of $\mathcal{O}(\kappa^{3/2}\ell_1)$; $\delta=0$ gives $\mathcal{O}(\kappa\ell_1)$ for \textbf{(1)}, \textbf{(2)}. Fixed budgets instead of expected counts.
& \textbf{(1)}, \textbf{(2)} require quasar structure; \textbf{(3)} only PL. Extra $\log(1/\gamma)$. At $\gamma=1$, $n=\kappa$, $\delta=L$, $\varepsilon\le\Delta$, \textbf{(2)} costs $\mathcal{O}(\kappa^2\ell_1)$. \\
\midrule
\textbf{(4)} Full-batch acceleration \citep{csqHinder2020}, Thm.~1; $\mathbb{R}^d$.
& $\mathcal{O}\!\left(\frac{n\sqrt\kappa}{\gamma}[1+\log(\kappa/\gamma)](1+\ell_\gamma)\right)$ component values and gradients; deterministic accuracy.
& \textbf{(1)}, \textbf{(2)}: gradients only. At $\delta=0$, save a factor $\asymp\min(n,\sqrt\kappa/\gamma)[1+\log(\kappa/\gamma)]$ for $\varepsilon\le\Delta$. \textbf{(2)} retains this gain if $\delta^2\lesssim\gamma^3\mu^2(n+\sqrt\kappa/\gamma)$.
& \textbf{(1)}, \textbf{(2)} give expected accuracy. \textbf{(1)} loses for $\delta\gg\sqrt{n\gamma}\,\mu^{1/4}L^{3/4}[1+\log(\kappa/\gamma)]^{3/2}$; \textbf{(2)} can also lose at large $\delta$. \\
\midrule
\textbf{(5)} QASVRG \citep{csqFu2023}, Cor.~3.6, Option I. Special case: $\mathbb{R}^d$, interpolation, $L_{\max}$-smooth components.
& $\begin{aligned}
&\mathcal{O}\!\Bigl(\Bigl[n+\min\!\Bigl\{\frac{\kappa_{\max}}{\gamma^2},\frac{n\sqrt{\kappa_{\max}}}{\gamma}\Bigr\}\\
&\quad{}\cdot[1+\log(1/\gamma)]\\
&\quad{}\cdot[1+\log(\kappa_{\max}/\gamma)]\Bigr](1+\ell_\gamma)\Bigr)
\end{aligned}$ component values and gradients.
& \textbf{(1)}, \textbf{(2)}: scope extends to non-interpolating sums; gradients only on the overlap. For \textbf{(2)}, a rate gain requires controlling $\delta^2/(\gamma^3\mu^2)$ relative to $n+\min\{\kappa_{\max}/\gamma^2,n\sqrt{\kappa_{\max}}/\gamma\}$.
& Small similarity is needed for a rate gain. On the overlap, \textbf{(2)} has quadratic $\kappa_{\max}$ dependence at $\delta=L_{\max}$; \textbf{(5)} has at most linear dependence, apart from logarithms. \\
\midrule
\textbf{(6)} AccSVRS \citep{csqLin2023v2}, \S3.3, App.~E. Special case: $\mathbb{R}^d$, $\gamma=1$, strongly convex $F$, $L$-smooth $f_1$, $\mu\le\delta\le L$.
& $\mathcal{O}\!\left(\left[n+n^{3/4}\left(\sqrt{\delta/\mu}+\sqrt{L/\delta}\right)\right.\right.$\newline
$\left.\left.{}+\sqrt{n\kappa}\right]\left[1+\log\frac{\sqrt n\delta+L}{\mu}\right](1+\ell_1)\right)$ expected component gradients, including inner solves.
& \textbf{(1)}: at $\delta=\mu$, $n\le\sqrt\kappa$, $\varepsilon\le\Delta$, $\mathcal{O}(\sqrt\kappa\ell_1)$ instead of $\mathcal{O}(n^{3/4}\sqrt\kappa[1+\log\kappa]\ell_1)$. Scope: \textbf{(1)}, \textbf{(2)} allow nonconvex $F$, $\delta\notin[\mu,L]$, and no $L$-smoothness requirement on $f_1$.
& \textbf{(1)} has worse polynomial cost at $\delta=L$, $\kappa>\sqrt n$: factor $(\kappa^2/n)^{1/12}$, excluding logarithms. Similarity dependence is not optimal on this convex overlap. \\
\bottomrule
\end{tabularx}
\begin{tablenotes}[flushleft]\footnotesize
\item[] \textbf{Abbr:} $\kappa=L/\mu$; $\kappa_{\max}=L_{\max}/\mu$; $L_{\max}$ is a common component smoothness bound. $\ell_a=\log(1+\Delta/(a\varepsilon))$ for $a\in\{1,\gamma\}$. PL means $\|\nabla F(x)\|^2\ge2\gamma^2\mu[F(x)-F(x_\star)]$. Costs count every component query, including full-gradient refreshes, batches, searches, and inner solves. $\lesssim$ suppresses universal constants.
\end{tablenotes}
\end{threeparttable}
\end{table}
\clearpage

\begin{table}[tp]
\centering
\caption{Finite-variance generalized-smooth convex optimization: \href{https://github.com/k3nfalt/longtails/blob/main/LongTails/Setups/FiniteVarianceGeneralizedSmoothAndConvex.lean}{FiniteVarianceGeneralizedSmoothAndConvex}.}\label{tab:FiniteVarianceGeneralizedSmoothAndConvex}
\parbox{\linewidth}{\footnotesize Convex $f$ on a separable real Hilbert space has differentiable gradient and $\|\nabla^2 f(x)\|\le L_0+L_1\|\nabla f(x)\|$. A minimizer $x^\star$ satisfies $\|x_0-x^\star\|\le R$ and $f(x_0)-f(x^\star)\le\Delta$. IID samples with marginal $p$ and joint law $P$ feed a jointly measurable oracle $g$, with integrable $g(x,\cdot)$ and squared error, $\int g(x,\omega)\,dp=\nabla f(x)$ and $\int\|g(x,\omega)-\nabla f(x)\|^2\,dp\le\sigma^2$. Known parameters satisfy $L_0,R,\Delta,\varepsilon>0$, $L_1,\sigma\ge0$, $0<\delta<1$. Target: measurable final output $\hat x$ with $P\{f(\hat x)-f(x^\star)>\varepsilon\}\le\delta$. Costs count all stochastic-gradient samples, including pilots, corrections, and padding.}\par\smallskip
\begin{threeparttable}
\footnotesize
\setlength{\tabcolsep}{4pt}\renewcommand{\arraystretch}{1.25}
\begin{tabularx}{\linewidth}{@{}>{\raggedright\arraybackslash}p{0.20\linewidth}>{\raggedright\arraybackslash}p{0.32\linewidth}>{\raggedright\arraybackslash}p{0.25\linewidth}>{\raggedright\arraybackslash}X@{}}
\toprule
\textbf{Result} & \textbf{Oracle complexity} & \textbf{Novelty} & \textbf{Limitations} \\
\midrule
\textbf{(1)} \href{https://github.com/k3nfalt/longtails/blob/main/LongTails/Theorems/AdaptiveRobustGeneralizedAcceleratedConvex.lean}{Adaptive\allowbreak Robust\allowbreak Generalized\allowbreak Accelerated\allowbreak Convex} \textbf{(Ours)} $(\star\star\star)$ \href{https://github.com/k3nfalt/longtails/blob/main/scripts/reviews/AdaptiveRobustGeneralizedAcceleratedConvex.md}{(review)}
& $\begin{aligned}[t]
\mathcal{O}\Bigl(& (W+v\ell)\Bigl[1+\frac{\log(eM/\delta)}{\log(e+W\sqrt M/(v^2\ell^2))}\Bigr]\\
&+v^2\ell+v\ell\log(eM/\delta)\Bigr).
\end{aligned}$
\par For fixed $c>0$, $\delta\ge(e+W)^{-c}$ implies $\mathcal{O}_c(W+v^2\ell)$.
& & \\
\midrule
\textbf{(2)} \href{https://github.com/k3nfalt/longtails/blob/main/LongTails/Theorems/FiniteVarianceGeneralizedAcceleratedConvex.lean}{Finite\allowbreak Variance\allowbreak Generalized\allowbreak Accelerated\allowbreak Convex} \textbf{(Ours)} $(\star\star)$ \href{https://github.com/k3nfalt/longtails/blob/main/scripts/reviews/FiniteVarianceGeneralizedAcceleratedConvex.md}{(review)}
& $\mathcal{O}\!\left((W+v^2\ell^2)\log\frac{e(1+T+q+\lambda)}{\delta}\right)$ & & \\
\midrule[1.2pt]
\textbf{(3)} RCSAG \citep{fvLobanov2026v1}, Thm.~9, App.~C; derived bound on $\mathbb R^d$, $p=2$, constant variance, $\Delta>\varepsilon$.
& $\mathcal{O}\!\left((\lambda+\sqrt T+qJ+v^2)\log\frac{e(1+\lambda)}{\delta}\right)$
& \textbf{(1)}: removes the $\ell$ factor on $\sqrt T$ when $\delta\ge(e+W)^{-c}$ and optimization dominates.
\par \textbf{(1)}, \textbf{(2)}: for $T=1$, $v\le1$, $\Delta=4\varepsilon$, $\delta=1/2$, $q\to\infty$, cost $\mathcal{O}(q)$ becomes $\mathcal{O}(q^{2/3})$, $\mathcal{O}(q^{2/3}\log q)$, respectively. Extend to separable Hilbert spaces.
& \textbf{(1)}, \textbf{(2)}: require a Hessian bound and finite variance; \textbf{(3)} permits first-order local smoothness, $p<2$, gradient/gap-dependent noise.
\par \textbf{(1)}: worse for $T=v=1$, $q=0$, $\Delta=2\varepsilon$, $\delta\to0$: $\mathcal{O}(\ell^2)$ versus $\mathcal{O}(\ell)$.
\par \textbf{(2)}: extra noise and optimization logarithms. \\
\midrule
\textbf{(4)} Clipped accelerated mirror descent \citep{fvNguyen2023v2}, Thm.~4.1; special case $L_1=0$, $p=2$, Euclidean geometry.
& $\mathcal{O}\!\left(1+(\sqrt T+v^2)\ell\right)$
& \textbf{(1)}: $\mathcal{O}(\sqrt T)$ versus $\mathcal{O}(\sqrt T\log T)$ for $T\to\infty$, $\delta=T^{-1}$, $v^2\ell\le\sqrt T$, using $\Delta\le L_0R^2/2$. No gain at constant $\delta$. \textbf{(1)}, \textbf{(2)} extend to $L_1>0$.
& \textbf{(1)}, \textbf{(2)} need $\Delta$ and differentiable gradients; \textbf{(4)} allows $p<2$ and mirror geometries. \textbf{(1)} has the tiny-$\delta$ loss shown against \textbf{(3)}. \\
\midrule
\textbf{(5)} \citep{fvTyurin2026v2}, Thm.~4.3; special case $\sigma=0$, exact gradients, $\Delta>\varepsilon$.
& $\begin{aligned}[t]\mathcal{O}\Bigl(&1+\sqrt T+q\bigl[1+\log_+\min\{\tfrac{L_1^2\Delta}{L_0},\tfrac{\Delta}{\varepsilon}\}\bigr]\\
&+\log_+\tfrac{\Delta}{L_0R^2}\Bigr)\end{aligned}$
& \textbf{(2)}: at $T=1$, $\Delta=4\varepsilon$, $\delta=1/2$, $q\to\infty$, $\mathcal{O}(q^{2/3}\log q)$ versus $\mathcal{O}(q)$. \textbf{(1)}, \textbf{(2)} extend to finite-variance stochastic gradients.
& \textbf{(2)} pays $\log(e(1+T+q+\lambda)/\delta)$ on $\sqrt T$ when it dominates, even with exact gradients. \\
\midrule
\textbf{(6)} \href{https://github.com/k3nfalt/longtails/blob/main/LongTails/Theorems/FiniteVarianceGeneralizedRadiusAcceleratedConvex.lean}{Finite\allowbreak Variance\allowbreak Generalized\allowbreak Radius\allowbreak Accelerated\allowbreak Convex} \textbf{(Ours)} $(\star\star\star)$ \href{https://github.com/k3nfalt/longtails/blob/main/scripts/reviews/FiniteVarianceGeneralizedRadiusAcceleratedConvex.md}{(review)}
& $\mathcal{O}\!\left((1+q+\sqrt T+v^2)H^3\right)$,
\par $H=\log\!\dfrac{e(2+T+q+v)}{\delta}$.
& & \\
\midrule[1.2pt]
\textbf{(7)} Clip-SSTM \citep{fvRadiusElKadi2026v2}, Thm.~2, Remark~1; derived bound on $\mathbb R^d$, $L_1>0$, $\alpha=2$, supplied gap.
& $\begin{aligned}[t]\mathcal{O}\Bigl(\bigl(& (1+q)[D+\log((1+q)H_D)]\\
&+\sqrt T+v^2\bigr)H_D\Bigr),\end{aligned}$
\par $D=\log(2+\Delta/(L_0R^2))$,
$H_D=\log\!\dfrac{e(2+T+q+v+D)}{\delta}$.
& \textbf{(6)}: curvature $qD\to q$ up to logs; gain when large $D$ dominates logarithmic overhead and unchanged terms. Updates and budget need no $\Delta$.
& \textbf{(6)}: extra logarithmic overhead and larger leading constant; covers only $\alpha=2$, whereas \textbf{(7)} allows $1<\alpha<2$. \\
\midrule
\textbf{(8)} Clip-SSTM \citep{fvRadiusElKadi2026v2}, Thm.~2, Remark~2; derived radius-only bound on $\mathbb R^d$, $L_1>0$, $\alpha=2$.
& $\begin{aligned}[t]\mathcal{O}\Bigl(\bigl(& (1+q)[q+\log((1+q)H)]\\
&+\sqrt T+v^2\bigr)H\Bigr).
\end{aligned}$
& \textbf{(6)}: curvature $q^2\to q$ up to logs. For $T=v=1$, $\delta=1/2$, $q\to\infty$: $\mathcal{O}(q\log^3q)$ versus $\mathcal{O}(q^2\log q)$. No total gain implied when unchanged terms dominate.
& \textbf{(6)}: extra logarithmic overhead and larger leading constant; excludes the $\alpha<2$ noise allowed by \textbf{(8)}. \\
\midrule
\textbf{(9)} Earlier repository bound \citep{fvRadiusEarlierRepository}; comparison with \textbf{(2)} on the identical setup.
& $\mathcal{O}\!\left((W+v^2\ell^2)\log\frac{e(1+T+q+\lambda)}{\delta}\right)$
& \textbf{(6)}: no $\Delta$ needed. For $\lambda\asymp q\to\infty$, $T=v=1$, $\delta=1/2$, admissible cosh instances give $\mathcal{O}(q\log^3q)$ versus $\mathcal{O}(q^{4/3}\log q)$.
& \textbf{(6)}: worse for small known gaps; e.g. fixed $\lambda,T,v,\delta$ gives $\mathcal{O}(q\log^3q)$ versus $\mathcal{O}(q^{2/3}\log q)$. Extra optimization/noise logarithms and larger leading constant. \\
\midrule
\textbf{(10)} Clipped SGD \citep{sgacGaash2025v2}, Thm.~1; derived bound on $\mathbb R^d$, sub-Gaussian noise: $\mathbb E\exp(\|g-\nabla f\|^2/\sigma^2)\le e$.
& $\mathcal{O}\!\left((1+T+q^2)H+v^2H^3\right)$ stochastic-gradient calls, including double sampling.
& \textbf{(6)}: $T\to\sqrt T$ and $q^2\to q$ up to logs on the sub-Gaussian overlap; extends to finite variance. Unchanged noise can mask a total gain.
& \textbf{(6)}: extra logarithmic factors on optimization terms. \\
\midrule
\textbf{(11)} \citep{fvRadiusRudomanov2026v1}, Thm.~3.7; special case $\sigma=0$, $L_1>0$, $\mathbb R^d$.
& $\mathcal{O}\!\left(\sqrt T+(1+q)\log(2+q)\right)$ exact-gradient calls, including line-search queries.
& \textbf{(6)}: extends the radius-only guarantee to finite-variance stochastic gradients; retains linear curvature dependence up to logarithms.
& \textbf{(6)}: at $\sigma=0$, pays $H^3$ on $\sqrt T$ and $(1+q)$, versus no logarithm and $\log(2+q)$, respectively. \\
\bottomrule
\end{tabularx}
\begin{tablenotes}[flushleft]\footnotesize
\item[] \textbf{Abbr:} $T=L_0R^2/\varepsilon$, $q=L_1R$, $v=\sigma R/\varepsilon$, $\ell=\log(e/\delta)$, $\lambda=1+\log(1+\Delta/\varepsilon)$, $W=1+\sqrt T+\lambda+(q\lambda)^{2/3}$, $M=1+T+q+\lambda+v^2\ell^2$; $J=1+\min\{\log(\Delta/\varepsilon),\log_+(L_1^2\Delta/L_0)\}$ in \textbf{(3)}. $\log_+u=\max\{0,\log u\}$, with $\log_+0=0$. The logarithmic fraction in \textbf{(1)} is zero at $v=0$. $c>0$ is fixed; constants in $\mathcal{O}_c$ may depend on $c$. $p$ in prior-work rows denotes the noise moment order.
\end{tablenotes}
\end{threeparttable}
\end{table}
\clearpage

\begin{table}[tp]
\centering
\caption{High-probability optimization under \href{https://github.com/k3nfalt/longtails/blob/main/LongTails/Setups/GapSmoothStrongGrowth.lean}{GapSmoothStrongGrowth}.}\label{tab:GapSmoothStrongGrowth}
\parbox{\linewidth}{\footnotesize Convex differentiable $f$ on a separable real Hilbert space with its Borel $\sigma$-algebra; a minimizer $x^\star$ satisfies $\|x_0-x^\star\|\le R$ and $f(x_0)-f(x^\star)\le\Delta$. For all $x,y$ with $(2\sqrt3+\sqrt6)\sqrt{H_1}\|y-x\|\le1$, $f(y)\le f(x)+\langle\nabla f(x),y-x\rangle+(H_0+H_1[f(x)-f(x^\star)])\|y-x\|^2$. The jointly measurable oracle is integrable, unbiased, and has integrable squared error satisfying $\int\|g(x,\omega)-\nabla f(x)\|^2\,dp(\omega)\le(\rho-1)\|\nabla f(x)\|^2$; samples are IID with marginal probability measure $p$ and joint law $P$. Parameters satisfy $H_0,H_1\ge0$, $R,\Delta>0$, $\rho\ge1$, $0<\varepsilon\le\Delta$, $0<\delta<1$. Target: a measurable output $\widehat x$ with $P(f(\widehat x)-f(x^\star)>\varepsilon)\le\delta$.}\par\smallskip
\begin{threeparttable}
\footnotesize
\setlength{\tabcolsep}{4pt}\renewcommand{\arraystretch}{1.25}
\begin{tabularx}{\linewidth}{@{}>{\raggedright\arraybackslash}p{0.20\linewidth}>{\raggedright\arraybackslash}p{0.32\linewidth}>{\raggedright\arraybackslash}p{0.25\linewidth}>{\raggedright\arraybackslash}X@{}}
\toprule
\textbf{Result} & \textbf{Oracle complexity} & \textbf{Novelty} & \textbf{Limitations} \\
\midrule
\textbf{(1)} \href{https://github.com/k3nfalt/longtails/blob/main/LongTails/Theorems/GapSmoothStrongGrowthHighProbability.lean}{Gap\allowbreak Smooth\allowbreak Strong\allowbreak Growth\allowbreak High\allowbreak Probability} \textbf{(Ours)} $(\star\star\star)$ \href{https://github.com/k3nfalt/longtails/blob/main/scripts/reviews/GapSmoothStrongGrowthHighProbability.md}{(review)}
& $\mathcal{O}\!\left(\rho\ell_\delta\left[(1+R\sqrt{H_1})\ell_\Delta\right.\right.$\newline
$\left.\left.{}+R\sqrt{H_0/\varepsilon}\right]\right)$
& & \\
\midrule[1.2pt]
\textbf{(2)} RCSAG \citep[Thm.~9, v1]{gssgLobanov2026}
& $\mathcal{O}\!\left(\ell_S\left[S+R\sqrt{H_0/\varepsilon}\right.\right.$\newline
\({}+R\sqrt{H_1}J_H\)\newline
$\left.\left.{}+(\rho-1)R^2(H_0/\varepsilon+H_1J_H)\right]\right)$.\par
$\mathbb R^d$; pure second-moment specialization; $Q=B(x_0,R)$.
& \textbf{(1)}: at fixed $\rho=2$, positive $H_0,H_1,R,\Delta,\delta$, $\varepsilon\to0$: $\mathcal{O}(\varepsilon^{-1/2})$ instead of $\mathcal{O}(\ell_S\varepsilon^{-1})$. More generally, stochastic-term gain when $R\sqrt{H_0/\varepsilon}\gg\rho\ell_\delta/(\rho-1)$.\par
$H_1$-term gain if $(\rho-1)R\sqrt{H_1}J_H\gg\rho\ell_\delta\ell_\Delta$.
& \textbf{(1)} excludes additive/gap noise, moments below two, and martingale sampling. Near $\rho=1$, its curvature term uses $\ell_\Delta$ instead of $J_H$; no uniform dominance. \\
\midrule
\textbf{(3)} \citep[App.~I, Thm.~7, v1]{gssgDvinskikh2026}
& $\mathcal{O}\!\left(\rho\left[1+L_1R\right.\right.$\newline
$\left.\left.{}+(1+L_0R^2/\varepsilon)^3\right]\operatorname{polylog}(Z)\right)$.\par
$\mathbb R^d$, $L_0,L_1>0$; convex $(L_0,L_1)$-smooth components, almost-sure strong growth; supplied gap certificate $\Delta$.
& \textbf{(1)}: $\varepsilon^{-3}\to\varepsilon^{-1/2}$ at fixed positive parameters on the overlap $H_0\asymp L_0$, $H_1\asymp L_1^2$. Scope: second moments replace almost-sure strong growth; no componentwise regularity.
& \\
\midrule
\textbf{(4)} Accelerated SGD \citep[Thm.~2]{gssgVaswani2019}
& $\mathcal{O}(1+\rho\sqrt T)$.\par
$H_1=0$, $\mathbb R^d$; expected target $\mathbb E[f(\widehat x)-f(x^\star)]\le\varepsilon$.
& \textbf{(1)}: scope extends to $H_1>0$; matches known $\rho\sqrt T$ accuracy dependence at $H_1=0$, with logarithmic confidence overhead.
& \textbf{(1)} uses batches $\lceil\rho\rceil$; its high-probability budget includes $\ell_\delta$ and $\ell_\Delta$, absent from this expected guarantee. \\
\midrule
\textbf{(5)} Derived bound: Markov conversion of \citep[Thm.~2]{gssgVaswani2019}
& $\mathcal{O}(1+\rho\sqrt{T/\delta})$.\par
$H_1=0$, $\mathbb R^d$; high-probability target.
& \textbf{(1)}: if $\Delta\le LR^2$ and $T\ge1$, costs $\mathcal{O}(\rho\sqrt T\,\ell_\delta)$ instead; confidence gain $1/(\sqrt\delta\,\ell_\delta)$.
& \textbf{(1)} uses batches $\lceil\rho\rceil$. \\
\midrule
\textbf{(6)} Derived bound: AGNES + proxBoost \citep[Thm.~15, v3]{gssgGupta2025}\citep[Thm.~7]{gssgDavis2020}
& $\mathcal{O}\!\left(\rho(1+\sqrt T)\ell_T[1+\log(\ell_T/\delta)]\right)$.\par
$H_1=0$, $\mathbb R^d$; high-probability target.
& \textbf{(1)}: for $\Delta\le LR^2$, $T\ge1$, removes the $\ell_T$ overhead, up to nested confidence logarithms. Scope extends to $H_1>0$; smooth-case acceleration is already known.
& \textbf{(1)} uses batches $\lceil\rho\rceil$. \\
\bottomrule
\end{tabularx}
\begin{tablenotes}[flushleft]\footnotesize
\item[] \textbf{Abbr:} Costs count individual stochastic-gradient samples, including batches and trials. $\ell_\delta=1+\log(1/\delta)$, $\ell_\Delta=\log(1+\Delta/\varepsilon)$, $S=\lceil\log_2(\Delta/\varepsilon)\rceil$, $\ell_S=1+\log((1+S)/\delta)$. For $H_0,H_1>0$, $J_H=\min\{S,\lceil\max\{0,\log_2(H_1\Delta/H_0)\}\rceil\}$; $J_H=0$ if $H_1=0$ and $J_H=S$ if $H_0=0<H_1$. $L=2H_0$, $T=LR^2/\varepsilon$, $\ell_T=1+\log(1+T)$. $L_0,L_1$ are component generalized-smoothness constants: $\|\nabla f_\omega(y)-\nabla f_\omega(x)\|\le(L_0+L_1\|\nabla f_\omega(x)\|)\|y-x\|$ for $\|y-x\|\le1/L_1$; almost-sure strong growth means $\|\nabla f_\omega(x)\|^2\le\rho\|\nabla f(x)\|^2$. $Z=2+L_1R+L_1^2\Delta/L_0+L_0R^2/\varepsilon+1/\delta$; $\operatorname{polylog}(Z)$ denotes the cited bound's unspecified fixed polynomial in $1+\log Z$.
\end{tablenotes}
\end{threeparttable}
\end{table}
\clearpage

\begin{table}[tp]
\centering
\caption{Exact monotone root finding under \href{https://github.com/k3nfalt/longtails/blob/main/LongTails/Setups/GeneralizedLipschitzMonotone.lean}{GeneralizedLipschitzMonotone}.}\label{tab:GeneralizedLipschitzMonotone}
\parbox{\linewidth}{\footnotesize A continuous monotone map $F:E\to E$ on a real inner-product space satisfies
$\|F(x)-F(y)\|\le (L_0+L_1\max_{0\le t\le1}\|F(x+t(y-x))\|)\|x-y\|$.
Here $L_0>0$, $L_1,G,R\ge0$, $\varepsilon>0$, $\|F(x_0)\|\le G$, and some root $x_\star$ satisfies $\|x_0-x_\star\|\le R$. The oracle returns exact $F(x)$; the target is $\|F(x)\|\le\varepsilon$ deterministically. The setup imposes no algorithm restriction.}\par\smallskip
\begin{threeparttable}
\footnotesize
\setlength{\tabcolsep}{4pt}\renewcommand{\arraystretch}{1.25}
\begin{tabularx}{\linewidth}{@{}>{\raggedright\arraybackslash}p{0.20\linewidth}>{\raggedright\arraybackslash}p{0.32\linewidth}>{\raggedright\arraybackslash}p{0.25\linewidth}>{\raggedright\arraybackslash}X@{}}
\toprule
\textbf{Result} & \textbf{Oracle complexity} & \textbf{Novelty} & \textbf{Limitations} \\
\midrule
\textbf{(1)} \href{https://github.com/k3nfalt/longtails/blob/main/LongTails/Theorems/AnchoredGeneralizedMonotone.lean}{Anchored\allowbreak Generalized\allowbreak Monotone} \textbf{(Ours)} $(\star\star\star)$ \href{https://github.com/k3nfalt/longtails/blob/main/scripts/reviews/AnchoredGeneralizedMonotone.md}{(review)}
& $\mathcal{O}\!\left(\begin{aligned}&1+\frac{L_0R}{\varepsilon}+(1+L_1R)\\[-2pt]&\quad\cdot\log\!\left(1+\frac{(L_0+L_1G)R}{\varepsilon}\right)\end{aligned}\right)$.
Also $\mathcal{O}(1+L_0R/\varepsilon+L_1^2R^2)$ using the intrinsic initial-residual bound.
& & \\
\midrule[1.2pt]
\textbf{(2)} EG \citep[Thm.~3.6]{glmChoudhuryLoizou2025v1}
& $\mathcal{O}(1+L_0^2R^2/\varepsilon^2+L_1^2R^2)$.
$E=\mathbb{R}^d$.
& \textbf{(1)}: $L_0R/\varepsilon$ instead of $(L_0R/\varepsilon)^2$; speedup $\Theta(L_0R/\varepsilon)$ if $\varepsilon\le\min\{L_0R,L_0/(L_1^2R)\}$. Scope: arbitrary real inner-product spaces. For $L_0=L_1=1$, $R=10$, $\varepsilon=10^{-3}$: $\mathcal{O}(10^4)$ versus $\mathcal{O}(10^8)$.
& \\
\midrule
\textbf{(3)} PolyakEG \citep[Thm.~3.4, Eq.~(19)]{glmYoonEtAl2026v1}
& $\mathcal{O}(1+L_0^2R^2/\varepsilon^2+L_1^2R^2)$.
Known $L_0,L_1$; $E=\mathbb{R}^d$.
& \textbf{(1)}: same leading-term gain as against \textbf{(2)}. If $L_1>0$ and $L_0/(L_1^2R)\le\varepsilon\le L_0/L_1$, speedup $\Omega((L_0/(L_1\varepsilon))^2)$. For $\varepsilon\ge L_0/L_1$, moderate $G$ gives $L_1R\log(2+L_1R+L_0R/\varepsilon)$ instead of $L_1^2R^2$; otherwise no strict gain is assured.
& \\
\midrule
\textbf{(4)} PolyakEG-LS \citep[Thm.~3.4, Prop.~3.5(a)]{glmYoonEtAl2026v1}
& Derived bound:
$\mathcal{O}\!\left(1+B+\frac{(\bar\lambda^0)^2R^2}{\varepsilon^2}+(\bar\lambda^1)^2R^2\right)$,
including all trials. $E=\mathbb{R}^d$; fixed $A\in(0,1)$, $\beta>2$, positive initial estimates $\lambda^0,\lambda^1$.
& \textbf{(1)}: $\varepsilon^{-1}$ instead of $\varepsilon^{-2}$ for fixed remaining parameters; no backtracking calls. Last iterate instead of minimum residual among observed iterates.
& \textbf{(1)} needs $L_0,L_1$ in its updates; \textbf{(4)} does not. \\
\midrule
\textbf{(5)} FEG / EAG \citep[Thm.~4.1, $\rho=0$]{glmLeeKim2021v3}\citep[Thms.~1--2, Sec.~2.3]{glmYoonRyu2021v2}
& $\mathcal{O}(1+L_0R/\varepsilon)$.
Special case: $L_1=0$, $E=\mathbb{R}^d$.
& \textbf{(1)} matches the known rate on this overlap; scope extends to generalized Lipschitz continuity with $L_1>0$.
& \textbf{(1)} has a roughly fivefold larger leading cost constant than FEG when $L_1=0$. \\
\bottomrule
\end{tabularx}
\begin{tablenotes}[flushleft]\footnotesize
\item[] \textbf{Abbr:} Costs count individual $F$-evaluations. $R$: initial root-distance bound; $G$: initial residual bound; $\varepsilon$: residual tolerance. The intrinsic bound is $\|F(x_0)\|\le L_0R e^{L_1R}$, so $G$ may be replaced by $\min\{G,L_0R e^{L_1R}\}$. Moderate $G$ means $\log(1+L_1G/L_0)=\mathcal{O}(\log(2+L_1R))$. For \textbf{(4)}, $A$ is the fixed acceptance tolerance, $\beta$ the backtracking factor, $\lambda^i$ the initial estimates,
$B=2\max\{0,\lceil\log(L_0/\lambda^0)/\log(\beta-1)\rceil,\lceil\log(L_1/\lambda^1)/\log(\beta-1)\rceil\}$,
and $\bar\lambda^i=\beta^{B/2}\lambda^i$; omit the $L_1$ logarithm when $L_1=0$. In \textbf{(2)}, interpret $L_0/(L_1^2R)=+\infty$ if $L_1R=0$.
\end{tablenotes}
\end{threeparttable}
\end{table}
\clearpage

\begin{table}[tp]
\centering
\caption{Oracle-complexity upper bounds for \href{https://github.com/k3nfalt/longtails/blob/main/LongTails/Setups/GeneralizedSmoothAndConvex.lean}{GeneralizedSmoothAndConvex}.}\label{tab:GeneralizedSmoothAndConvex}
\parbox{\linewidth}{\footnotesize Convex $f:E\to\mathbb R$ on a real Hilbert space has a differentiable gradient satisfying $\|\nabla^2 f(x)\|\le L_0+L_1\|\nabla f(x)\|$ everywhere. A minimizer $x_\star$ satisfies $\|x_0-x_\star\|\le R$ and $f(x_0)-f(x_\star)\le\Delta$, with $L_0,R,\Delta,\varepsilon>0$ and $L_1\ge0$. Exact gradients are available; $L_0,L_1,R,\Delta$ are public inputs. The target is $f(x_{\rm out})-f(x_\star)\le\varepsilon$ almost surely; our constructions guarantee it for every seed. The setup imposes no algorithm restriction.}\par\smallskip
\begin{threeparttable}
\footnotesize
\setlength{\tabcolsep}{4pt}\renewcommand{\arraystretch}{1.25}
\begin{tabularx}{\linewidth}{@{}>{\raggedright\arraybackslash}p{0.20\linewidth}>{\raggedright\arraybackslash}p{0.32\linewidth}>{\raggedright\arraybackslash}p{0.25\linewidth}>{\raggedright\arraybackslash}X@{}}
\toprule
\textbf{Result} & \textbf{Oracle complexity} & \textbf{Novelty} & \textbf{Limitations} \\
\midrule
\textbf{(1)} \href{https://github.com/k3nfalt/longtails/blob/main/LongTails/Theorems/GeneralizedSelfNormalizedConvex.lean}{Generalized\allowbreak Self\allowbreak Normalized\allowbreak Convex} \textbf{(Ours)} $(\star\star\star)$ \href{https://github.com/k3nfalt/longtails/blob/main/scripts/reviews/GeneralizedSelfNormalizedConvex.md}{(review)}
& $\begin{aligned}[t]\mathcal{O}\Bigl(&1+\sqrt{L_0R^2/\varepsilon}\\ &+\min\{\lceil q^{2/3}\rceil\ell,\\ &\qquad\ell+(q\ell)^{2/3}\}\Bigr)\end{aligned}$ gradients
& & \\
\midrule
\textbf{(2)} \href{https://github.com/k3nfalt/longtails/blob/main/LongTails/Theorems/GeneralizedAcceleratedConvex.lean}{Generalized\allowbreak Accelerated\allowbreak Convex} \textbf{(Ours)} $(\star\star\star)$ \href{https://github.com/k3nfalt/longtails/blob/main/scripts/reviews/GeneralizedAcceleratedConvex.md}{(review)}
& $\mathcal{O}\!\left(1+\sqrt{L_0R^2/\varepsilon}+\lceil q^{2/3}\rceil\ell\right)$ gradients
& & \\
\midrule[1.2pt]
\textbf{(3)} AGD \citep[Theorem~4.3, v2]{cgsTyurin2026v2}\par $E=\mathbb R^d$; derived bound
& $\begin{aligned}[t]\mathcal{O}\Bigl(&1+\sqrt{L_0R^2/\varepsilon}\\ &+q\log\!\left(1+\min\!\left\{\frac{L_1^2\Delta}{L_0},\frac{\Delta}{\varepsilon}\right\}\right)\\ &+\log\!\left(1+\frac{\Delta}{L_0R^2}\right)\Bigr)\end{aligned}$ gradients.\par $\bar R=R$, $\Gamma_0=2\Delta/R^2$, using a flat-coordinate extension.
& \textbf{(1)}: for $q\to\infty$, $\Delta=L_0R^2$, $\varepsilon=L_0R^2q^{-4/3}$, total $\mathcal{O}(q\log q)\to\mathcal{O}((q\log q)^{2/3})$.\par \textbf{(2)}: $\mathcal{O}(q^{2/3}\log q)$ in that regime. Matches the leading $\varepsilon^{-1/2}$ rate.
& \textbf{(2)}: startup can be worse when $\varepsilon<L_0/L_1^2$ and $q^{2/3}\ell>\sqrt{L_0R^2/\varepsilon}$; its logarithm does not saturate.\par \textbf{(1)}, \textbf{(2)}: two gradients per iteration instead of one; restricted to full domains. \\
\midrule
\textbf{(4)} AGMsDR \citep[Theorem~6.2, v3]{gscVankov2025v3}\par $E=\mathbb R^d$, $f\in C^2$
& $\begin{aligned}[t]\mathcal{O}\Bigl((1+\nu)\bigl[&1+\sqrt{L_0R^2/\varepsilon}\\ &+\lceil q^{2/3}\rceil\ell\bigr]\Bigr)\end{aligned}$ gradient/value calls, including exact segment searches; $\nu$ has no parameter bound.
& \textbf{(1)}: startup $q^{2/3}\ell\to(q\ell)^{2/3}$ for $q\ge\sqrt\ell$.\par \textbf{(2)}: matches iteration rate.\par \textbf{(1)}, \textbf{(2)}: bounded gradient-only cost; no values or exact search. Total-call rates are not comparable with unbounded $\nu$.
& \textbf{(1)}, \textbf{(2)}: need $R,\Delta$ in addition. \\
\midrule
\textbf{(5)} Accelerated GRAAL \citep[Corollary~3, v2]{gscBorodichKovalev2025v2}\par $E=\mathbb R^d$, $f\in C^2$, $L_1>0$
& $\begin{aligned}[t]\mathcal{O}\Bigl(&1+\sqrt{L_0R^2/\varepsilon}+q^3\\ &+(1+q^2)\log\frac1{\eta_0L_0}\Bigr)\end{aligned}$ gradient/value calls, for $0<\eta_0\le e^{-q}/L_0$.
& \textbf{(1)}, \textbf{(2)}: no values. For $q\to\infty$, $\eta_0=e^{-q}/L_0$, $\Delta=L_0R^2$, $\varepsilon=L_0R^2q^{-4}$, total $\mathcal{O}(q^3)\to\mathcal{O}(q^2)$.\par \textbf{(2)}: startup $q^3\to q^{2/3}\ell$ when $\ell\ll q^{7/3}$.
& \textbf{(1)}, \textbf{(2)}: need $L_0,L_1,R,\Delta$; \textbf{(5)} updates need only an admissible $\eta_0$. Our startup depends on accuracy. \\
\bottomrule
\end{tabularx}
\begin{tablenotes}[flushleft]\footnotesize
\item[] \textbf{Abbr:} $q=L_1R$; $\ell=1+\log(1+\Delta/\varepsilon)$; $\nabla^2 f=D\nabla f$; $\eta_0$ is the initial stepsize; $\nu$ is the maximum gradient/value call cost of one exact segment minimization. Costs count all gradient evaluations, including trials, and values where used.
\end{tablenotes}
\end{threeparttable}
\end{table}
\clearpage

\begin{table}[tp]
\centering
\caption{Comparison for \href{https://github.com/k3nfalt/longtails/blob/main/LongTails/Setups/GeneralizedSmoothParameterFree.lean}{GeneralizedSmoothParameterFree}.}\label{tab:GeneralizedSmoothParameterFree}
\parbox{\linewidth}{\footnotesize Convex $f:E\to\mathbb R$ on a real Hilbert space, with differentiable gradient and $\|\nabla^2 f(x)\|\le L_0+L_1\|\nabla f(x)\|$ everywhere. A minimizer $x_\star$ satisfies $\|x_0-x_\star\|\le R$ and $f(x_0)-f(x_\star)\le\Delta$, where $L_0,R,\Delta,\varepsilon,h_0>0$ and $L_1\ge0$. Exact gradients are available. Target: $f(\hat x)-f(x_\star)\le\varepsilon$. In \textbf{(1)}, updates receive only $x_0,h_0$; the sufficient stopping budget depends on the problem parameters and $\varepsilon$.}\par\smallskip
\begin{threeparttable}
\footnotesize
\setlength{\tabcolsep}{4pt}\renewcommand{\arraystretch}{1.25}
\begin{tabularx}{\linewidth}{@{}>{\raggedright\arraybackslash}p{0.20\linewidth}>{\raggedright\arraybackslash}p{0.32\linewidth}>{\raggedright\arraybackslash}p{0.27\linewidth}>{\raggedright\arraybackslash}X@{}}
\toprule
\textbf{Result} & \textbf{Oracle complexity} & \textbf{Novelty} & \textbf{Limitations} \\
\midrule
\textbf{(1)} \href{https://github.com/k3nfalt/longtails/blob/main/LongTails/Theorems/BacktrackingGeneralizedAccelerated.lean}{Backtracking\allowbreak{}Generalized\allowbreak{}Accelerated} \textbf{(Ours)} $(\star\star\star)$ \href{https://github.com/k3nfalt/longtails/blob/main/scripts/reviews/BacktrackingGeneralizedAccelerated.md}{(review)}
& $\begin{aligned}[t]
\mathcal{O}\Bigl(&\sqrt{L_0R^2/\varepsilon}+(1+s)\log\Bigl(2+s\\
&+\frac{1}{L_0h_0}+L_0h_0\\
&+\frac{L_1^2\Delta}{L_0}+L_1^2\Delta h_0\Bigr)\Bigr)
\end{aligned}$ gradient calls, including two per trial and all rejected trials.
& & \\
\midrule[1.2pt]
\textbf{(2)} Accelerated GRAAL \citep[Cor.~3, v2]{gspfBorodichKovalev2025v2}. Restriction: $E=\mathbb R^d$, $f\in C^2$, $L_1>0$.
& $\begin{aligned}[t]
\mathcal{O}\Bigl(&1+\sqrt{L_0R^2/\varepsilon}+s^3\\
&+(1+s^2)\log\frac{1}{L_0\eta_0}\Bigr)
\end{aligned}$ gradient and value calls; requires $L_0\eta_0e^s\le1$.
& \textbf{(1)}: for $s\gg1$, $h_0=\eta_0$, $\Delta=f(x_0)-f(x_\star)$, proved start-up improves by a factor of order $s$; total gain when $\varepsilon\gtrsim L_0R^2/s^6$. Allows any $h_0>0$; removes values. Matches order for $s\lesssim1$ or accuracy-dominated costs.\par
Example: $f(x)=e^x-x$, $L_0=L_1=1$, $x_0=R\gg1$, $h_0=\eta_0=e^{-R}$: \textbf{(1)} $\mathcal{O}(R/\sqrt\varepsilon+R^2)$ versus \textbf{(2)} $\mathcal{O}(R/\sqrt\varepsilon+R^3)$.
& \textbf{(1)}: backtracking incurs constant-factor overhead. \\
\midrule
\textbf{(3)} AGD \citep[Thm.~4.3, v2]{cgsTyurin2026v2}. Restriction: $E=\mathbb R^d$; $R=\|x_0-x_\star\|$ for the displayed tuning.
& $\begin{aligned}[t]
\mathcal{O}\Bigl(&1+\sqrt{L_0R^2/\varepsilon}\\
&+s\log\Bigl(2+\min\Bigl\{\frac{L_1^2\Delta}{L_0},\frac{\Delta}{\varepsilon}\Bigr\}\Bigr)\\
&+\log\Bigl(2+\frac{\Delta}{L_0R^2}\Bigr)\Bigr)
\end{aligned}$ gradient calls; supplied $L_0,L_1$, $\bar R=R$, $\Gamma_0=2\Delta/R^2$.
& \textbf{(1)}: removes $L_0,L_1,R,\Delta$ from updates; matches the leading $\varepsilon^{-1/2}$ term. No established rate improvement.
& \textbf{(1)}: possible $(1+s)\log(L_0h_0+1/(L_0h_0))$ and $(1+s)\log(2+s)$ overhead; loses the shortcut when $\Delta/\varepsilon<L_1^2\Delta/L_0$. Backtracking constant factor. \\
\midrule
\textbf{(4)} AGMsDR \citep[Thm.~6.2, v3]{gspfVankovEtAl2025v3}. Restriction: $E=\mathbb R^d$, $f\in C^2$.
& $\begin{aligned}[t]
\mathcal{O}\Bigl((\nu+1)\Bigl[&1+\sqrt{L_0R^2/\varepsilon}\\
&+\lceil s^{2/3}\rceil\bigl(1+\log(\Delta/\varepsilon)\bigr)\Bigr]\Bigr)
\end{aligned}$ gradient/value calls; $0<\varepsilon\le\Delta=f(x_0)-f(x_\star)$. Requires $L_0,L_1$ and exact one-dimensional minimization.
& \textbf{(1)}: no $L_0,L_1$, values, or exact line search; start-up logarithm is independent of $\varepsilon$. Total rates are not comparable: $\nu$ has no parameter bound.
& \textbf{(1)} loses the $s^{2/3}$ start-up power, ignoring the unbounded search factor $\nu$. \\
\bottomrule
\end{tabularx}
\begin{tablenotes}[flushleft]\footnotesize
\item[] \textbf{Abbr:} $s=L_1R$; $h_0$ and $\eta_0$ are initial step sizes; $\bar R$ is a supplied distance bound; $\Gamma_0$ is AGD's initial potential parameter; $\nu$ bounds the gradient/value calls per exact one-dimensional search. $\log$ denotes the natural logarithm.
\end{tablenotes}
\end{threeparttable}
\end{table}
\clearpage

\begin{table}[tp]
\centering
\caption{Heavy-tailed Adam under the \href{https://github.com/k3nfalt/longtails/blob/main/LongTails/Setups/HeavyTailedAdam.lean}{HeavyTailedAdam setup}.}\label{tab:HeavyTailedAdam}
\parbox{\linewidth}{\footnotesize $f:\mathbb{R}^d\to\mathbb{R}$ is differentiable, $\nabla f$ is $L$-Lipschitz, and $f(x_0)-\inf f\le\Delta$, with $L>0$, $\Delta\ge0$. A jointly measurable, integrable oracle uses IID samples of law $\mu$, with $\int g(x,\omega)\,d\mu=\nabla f(x)$ and $\int\sum_i|g_i(x,\omega)-\partial_i f(x)|^\alpha\,d\mu\le\sigma^\alpha$, $1<\alpha\le2$, $\sigma\ge0$. The setup restricts algorithms to coordinate Adam. \textbf{(1)} uses no bias correction, $\beta_1=0.9$, and damping outside the square root; its uniform pre-update output satisfies $\int\|\nabla f(x_I)\|_2\,dP\le\varepsilon$, $\varepsilon>0$, under execution law $P$.}\par\smallskip
\begin{threeparttable}
\footnotesize
\setlength{\tabcolsep}{4pt}\renewcommand{\arraystretch}{1.25}
\begin{tabularx}{\linewidth}{@{}>{\raggedright\arraybackslash}p{0.18\linewidth}>{\raggedright\arraybackslash}p{0.30\linewidth}>{\raggedright\arraybackslash}p{0.23\linewidth}>{\raggedright\arraybackslash}X@{}}
\toprule
\textbf{Result} & \textbf{Oracle complexity} & \textbf{Novelty} & \textbf{Limitations} \\
\midrule
\textbf{(1)} \href{https://github.com/k3nfalt/longtails/blob/main/LongTails/Theorems/AdamHeavyTailedFixedMomentum.lean}{Adam\allowbreak Heavy\allowbreak Tailed\allowbreak Fixed\allowbreak Momentum} \textbf{(Ours)} $(\star\star\star)$ \href{https://github.com/k3nfalt/longtails/blob/main/scripts/reviews/AdamHeavyTailedFixedMomentum.md}{(review)}
& $\mathcal{O}\!\left(1+\dfrac{L\Delta}{\varepsilon^2}\left[1+c_\alpha\left(\dfrac{\sigma}{\varepsilon}\right)^r\right]\right)$ & & \\
\midrule[1.2pt]
\textbf{(2)} NSGD-M \citep[Cor.~19, arXiv v3]{htaHubler2025}
& $\mathcal{O}\!\left(1+\dfrac{L\Delta}{\varepsilon^2}\left[1+\left(\dfrac{\sigma_2}{\varepsilon}\right)^r\right]\right)$
& \textbf{(1)}: Matches known accuracy exponent with fixed $\beta_1$ and no normalization; method-specific advance, no problem-class rate improvement. No tuning-knowledge gain.
& \textbf{(1)}: Larger constants; narrower noise model. Euclidean-scale conversion multiplies its noise term by $d^{(2-\alpha)/(2(\alpha-1))}$. \\
\midrule
\textbf{(3)} SGD \citep[Thm.~4.1, constant steps]{htaLiu2026}
& $\mathcal{O}\!\left(1+\dfrac{L\Delta}{\varepsilon^2}+\left(\dfrac{\sigma_2\sqrt{L\Delta}}{\varepsilon^2}\right)^r\right)$\par $\varepsilon\le\sqrt{2L\Delta}$.
& \textbf{(1)}: For fixed $\alpha<2$, noise-term gain scales as $(\sqrt{L\Delta}/\varepsilon)^{(2-\alpha)/(\alpha-1)}$ at equal scales; after Euclidean conversion, $(\sqrt{L\Delta}/(\sqrt d\,\varepsilon))^{(2-\alpha)/(\alpha-1)}$ when $\sqrt{L\Delta}>\sqrt d\,\varepsilon$. No exponent gain at $\alpha=2$.
& \textbf{(1)}: Coordinate-noise penalty and large constants; loses the stronger guarantee $\mathbb{E}[(K^{-1}\sum_{k=1}^K\|\nabla f(x_k)\|_2^2)^{\alpha/2}]$ $\le\varepsilon^\alpha$. Gain concerns the shared expected-norm target. \\
\midrule
\textbf{(4)} Vector Adam \citep[Thm.~6.4, v2]{htaPang2026}\par $G$-Lipschitz differentiable objectives; $\alpha>4/3$.
& $\mathcal{O}\!\left(1+\left(\dfrac{B}{\varepsilon}\right)^{\frac{8\alpha-4}{3\alpha-4}}\left[1+\log\dfrac{G+\sigma_2}{\varepsilon}\right]\right.$\par $\left.{}+\dfrac{\Delta\sqrt\lambda\,B^{\frac{5\alpha}{3\alpha-4}}}{\varepsilon^{\frac{5\alpha}{3\alpha-4}+3/2}}\right)$\par Expected relaxed $(\lambda,\varepsilon)$-stationarity; $0<\varepsilon<G+\sigma_2$, fixed $\alpha$ and memory coupling.
& \textbf{(1)}: Ordinary expected gradient norm; no bounded-gradient assumption; includes $1<\alpha\le4/3$. Different criteria and classes: rates not comparable.
& \textbf{(1)}: Requires smoothness. Accuracy-dependent damping replaces a fixed stabilizer; $1-\beta_2\le10^{-12}$. \\
\bottomrule
\end{tabularx}
\begin{tablenotes}[flushleft]\footnotesize
\item[] \textbf{Abbr:} Costs count individual stochastic-gradient calls. $r=\alpha/(\alpha-1)$; $c_\alpha=10^{5\alpha/(\alpha-1)}$; $K$ is the iteration count. $\sigma_2$ bounds $(\mathbb{E}\|g-\nabla f\|_2^\alpha)^{1/\alpha}$: the setup implies $\sigma_2\le\sigma$, while Euclidean noise permits $\sigma=d^{1/\alpha-1/2}\sigma_2$. Rate comparisons hold at fixed $\alpha$; prior-work constants may depend on $\alpha$. In \textbf{(4)}, $G$ bounds $\|\nabla f\|_2$, $\lambda>0$ is the relaxation parameter, $B=CG+\sigma_2$, and $C$ is the regret constant in Thm.~6.4 (bounded for fixed memory coupling and sufficiently small accuracy). Relaxed stationarity uses $\inf_{\mathbb{E}Y=x}\{\|\mathbb{E}\nabla f(Y)\|_2+\lambda\mathbb{E}\|Y-x\|_2^2\}$. In \textbf{(1)}, $\beta_2=1-\rho$, $\tau=\varepsilon/(10^5\sqrt\rho)$, and $\rho=(10^{12}+(10^5\sigma/\varepsilon)^r)^{-1}$; tuning uses $L,\Delta,\sigma,\alpha,\varepsilon$.
\end{tablenotes}
\end{threeparttable}
\end{table}
\clearpage

\begin{table}[tp]
\centering
\caption{Inexact-gradient strongly convex optimization: \href{https://github.com/k3nfalt/longtails/blob/main/LongTails/Setups/InexactGradientStronglyConvex.lean}{InexactGradientStronglyConvex}.}\label{tab:InexactGradientStronglyConvex}
\parbox{\linewidth}{\footnotesize On a real Hilbert space, $f$ is differentiable with global quadratic upper and lower models of curvature $L$ and $\mu$, where $0<\mu\le L$. A minimizer $x^\star$ satisfies $\|x_0-x^\star\|\le R$. Every oracle reply obeys $\|g-\nabla f(x)\|\le\alpha\|\nabla f(x)\|+\delta$, with $0\le\alpha<1$ and $\delta,R\ge0$; errors may depend on the call and seed. The target is $f(x_N)-f(x^\star)\le\varepsilon$ almost surely, with $\varepsilon>0$ and $136\delta^2\le\mu(1-\alpha)^2\varepsilon$. The constructed method uses $L,\mu,\alpha,x_0$; the setup imposes no algorithm restriction.}\par\smallskip
\begin{threeparttable}
\footnotesize
\setlength{\tabcolsep}{4pt}\renewcommand{\arraystretch}{1.25}
\begin{tabularx}{\linewidth}{@{}>{\raggedright\arraybackslash}p{0.20\linewidth}>{\raggedright\arraybackslash}p{0.32\linewidth}>{\raggedright\arraybackslash}p{0.25\linewidth}>{\raggedright\arraybackslash}X@{}}
\toprule
\textbf{Result} & \textbf{Oracle complexity} & \textbf{Novelty} & \textbf{Limitations} \\
\midrule
\textbf{(1)} \href{https://github.com/k3nfalt/longtails/blob/main/LongTails/Theorems/AcceleratedInexactStronglyConvex.lean}{Accelerated\allowbreak Inexact\allowbreak Strongly\allowbreak Convex} \textbf{(Ours)} $(\star\star\star)$ \href{https://github.com/k3nfalt/longtails/blob/main/scripts/reviews/AcceleratedInexactStronglyConvex.md}{(review)}
& $\mathcal{O}\!\left(1+\left[\sqrt{\frac{\kappa}{1-\alpha}}+\frac{\alpha^2\kappa}{1-\alpha}\right]\ell\right)$ & & \\
\midrule[1.2pt]
\textbf{(2)} RE-AGM \citep[Thm.~6.1]{igsVasinGasnikov2026v1}
& $\mathcal{O}\!\left(1+(\sqrt\kappa+\alpha\kappa)\ell\right)$; $\alpha\le1/3$, $\varepsilon\ge C\min\{\sqrt\kappa,\alpha^{-1}\}\delta^2/\mu$. Euclidean space.
& \textbf{(1)}: for $\kappa^{-1/2}\le\alpha\le1/3$, gain $\min\{\alpha\sqrt\kappa,\alpha^{-1}\}$. At $\alpha=\kappa^{-1/4}$, $\mathcal{O}(1+\kappa^{1/2}\ell)$ instead of $\mathcal{O}(1+\kappa^{3/4}\ell)$. Certified floor drops to $\mathcal{O}(\delta^2/\mu)$; extends to $\alpha<1$.
& \\
\midrule
\textbf{(3)} RE-AGM \citep[Thm.~4]{igsVasinEtAl2025v3}
& $\mathcal{O}(1+\sqrt\kappa\,\ell)$; $\delta=0$, $\alpha<c_0\kappa^{-1/2}$. Euclidean space.
& \textbf{(1)}: matches this rate; retains it up to $\alpha=\mathcal{O}(\kappa^{-1/4})$ bounded away from $1$, instead of $\mathcal{O}(\kappa^{-1/2})$. Scope: permits absolute error subject to its floor.
& \\
\midrule
\textbf{(4)} RE-AGM \citep[Thm.~5]{igsVasinEtAl2025v3}
& $\mathcal{O}(1+\alpha\kappa\ell)$; $\delta=0$, $\kappa^{-1/2}/3\le\alpha\le1/3$. Euclidean space.
& \textbf{(1)}: gain $\min\{\alpha\sqrt\kappa,\alpha^{-1}\}$ for $\kappa^{-1/2}\le\alpha\le1/3$; reaches $\kappa^{1/4}$ at $\alpha=\kappa^{-1/4}$. Scope: adds absolute error subject to its floor.
& \\
\midrule
\textbf{(5)} GD \citep[Thm.~5.2]{igsVasinGasnikov2026v1}
& $\mathcal{O}\!\left(1+\frac{\kappa(1+\alpha)}{(1-\alpha)^3}\ell\right)$; $\varepsilon\ge C\frac{1+\alpha}{(1-\alpha)^3}\frac{\delta^2}{\mu}$. Euclidean smooth PL functions; displayed initialization on the common strongly convex setup.
& \textbf{(1)}: for large $\kappa$, $\alpha\le\kappa^{-1/4}$ and common admissible accuracy, $\mathcal{O}(1+\sqrt\kappa\,\ell)$ instead of $\mathcal{O}(1+\kappa\ell)$.
& \textbf{(1)} requires strong convexity instead of the PL inequality. \\
\midrule
\textbf{(6)} Fast gradient method \citep[Thm.~3.3]{igsGaucherEtAl2026v1}
& $\mathcal{O}(1+\sqrt\kappa\,\ell)$; $\alpha=0$, $\varepsilon\ge C\sqrt\kappa\,\delta^2/\mu$. Euclidean space; additive-error specialization.
& \textbf{(1)}: at $\alpha=0$, certified floor $\delta^2/\mu$ instead of $\sqrt\kappa\,\delta^2/\mu$ in order; certifies acceleration below \textbf{(6)}'s threshold when $\delta^2/\mu\lesssim\varepsilon\ll\sqrt\kappa\,\delta^2/\mu$. Scope: adds relative error.
& \\
\bottomrule
\end{tabularx}
\begin{tablenotes}[flushleft]\footnotesize
\item[] \textbf{Abbr:} $\kappa=L/\mu$, $\ell=\log(1+LR^2/\varepsilon)$; costs count inexact-gradient calls. $C$ denotes a sufficiently large universal constant in sufficient accuracy conditions; $\alpha^{-1}=+\infty$ at $\alpha=0$. $c_0=(\sqrt2-1)/(18\sqrt2)$. PL: $\|\nabla f(x)\|^2\ge2\mu(f(x)-f(x^\star))$.
\end{tablenotes}
\end{threeparttable}
\end{table}
\clearpage

\begin{table}[tp]
\centering
\caption{Known-value convex optimization with exact local linear minimization: \href{https://github.com/k3nfalt/longtails/blob/main/LongTails/Setups/KnownValueSmoothConvexLocalLinear.lean}{KnownValueSmoothConvexLocalLinear}.}\label{tab:KnownValueSmoothConvexLocalLinear}
\parbox{\linewidth}{\footnotesize A real Hilbert space, a closed convex set $C$, and $x_0\in C$. The differentiable objective satisfies $0\le f(y)-f(x)-\langle\nabla f(x),y-x\rangle\le (L/2)\|y-x\|^2$ for $x,y\in C$. Known inputs are $L>0$, $R\ge0$, $\varepsilon>0$, and $f_\star=f(x_\star)$ for a minimizer $x_\star\in C$ with $\|x_0-x_\star\|\le R$. Exact oracles return $(f(x),\nabla f(x))$ and a minimizer of $\langle g,v\rangle$ over $v\in C\cap B(x,r)$. Target: feasible $y$ with $f(y)-f(x_\star)\le\varepsilon$; \textbf{(1)} guarantees its last primal iterate for every seed and oracle tie-breaking choice.}\par\smallskip
\begin{threeparttable}
\footnotesize
\setlength{\tabcolsep}{4pt}\renewcommand{\arraystretch}{1.25}
\begin{tabularx}{\linewidth}{@{}>{\raggedright\arraybackslash}p{0.20\linewidth}>{\raggedright\arraybackslash}p{0.32\linewidth}>{\raggedright\arraybackslash}p{0.25\linewidth}>{\raggedright\arraybackslash}X@{}}
\toprule
\textbf{Result} & \textbf{Oracle complexity} & \textbf{Novelty} & \textbf{Limitations} \\
\midrule
\textbf{(1)} \href{https://github.com/k3nfalt/longtails/blob/main/LongTails/Theorems/SplitBroxKnownValue.lean}{Split\allowbreak Brox\allowbreak Known\allowbreak Value} \textbf{(Ours)} $(\star\star\star)$ \href{https://github.com/k3nfalt/longtails/blob/main/scripts/reviews/SplitBroxKnownValue.md}{(review)}
& $\mathcal{O}(1+\sqrt{LR^2/\varepsilon})$ total value-gradient and local-LMO calls, including all trials.
& & \\
\midrule[1.2pt]
\textbf{(2)} BroxGD \citep[Thm.~3.2(ii), Cor.~D.2]{kvBrox2026v2}
& $\mathcal{O}(1+LR^2/\varepsilon)$ total value-gradient and local-LMO calls. Finite-dimensional; $R>0$; $L$-Lipschitz gradient on an open convex neighborhood of $C$.
& \textbf{(1)}: factor $\sqrt{LR^2/\varepsilon}$ fewer calls when $LR^2/\varepsilon\ge1$; last instead of best iterate. Scope: quadratic models only on $C$ suffice.
& \textbf{(1)} needs exact $f_\star$ and $L$ additionally; \textbf{(2)} uses only $R,K$ to set radii. \\
\midrule
\textbf{(3)} Polyak-radius BroxGD \citep[Thm.~3.5, Eq.~(24)]{kvBrox2026v2}
& $\mathcal{O}\!\left(1+\dfrac{(\|\nabla f(x_\star)\|+LR)^2R^2}{\varepsilon^2}\right)$ total value-gradient and local-LMO calls. Finite-dimensional; smooth on an open neighborhood of $C$; known $f_\star$.
& \textbf{(1)}: $\varepsilon^{-1/2}$ instead of $\varepsilon^{-2}$ in this general guarantee, without $\|\nabla f(x_\star)\|$ dependence; last instead of averaged iterate.
& \textbf{(1)} needs $L,R$ additionally. \\
\midrule
\textbf{(4)} APMM \citep[Alg.~1, Thm.~2.1]{kvAPMM2026v3}; $\rho=1$, $X_k=C$, no functional constraints
& $\mathcal{O}(1+\sqrt{LR^2/\varepsilon})$ total value-gradient evaluations and exact projections onto $C\cap\{\ell_f\le f_\star\}$. Finite-dimensional smooth case.
& \textbf{(1)} matches the accelerated order using point-only local LMOs, with no logarithmic search overhead. Different oracles prevent a direct cost comparison.
& \textbf{(1)} needs $L,R$ additionally; \textbf{(4)} is parameter-free apart from $f_\star,\varepsilon$. \\
\midrule
\textbf{(5)} CGS \citep[Cor.~3, 2014 preprint]{kvCGS2016}; special case: bounded $C\subseteq\mathbb{R}^d$ of diameter $D$
& $\mathcal{O}(1+LD^2/\varepsilon)$ total gradient and global-LMO calls; gradient calls alone: $\mathcal{O}(1+\sqrt{LD^2/\varepsilon})$.
& \textbf{(1)}: for $R\le D$ and $LD^2/\varepsilon\ge1$, at least a factor $\sqrt{LD^2/\varepsilon}$ fewer linear-oracle calls. Scope: unbounded $C$ allowed.
& \textbf{(1)} needs $f_\star$ and the stronger local LMO; \textbf{(5)} needs only global linear minimization over $C$. \\
\bottomrule
\end{tabularx}
\begin{tablenotes}[flushleft]\footnotesize
\item[] \textbf{Abbr:} LMO: linear minimization oracle; local LMO: one point minimizing a linear function on $C\cap B(x,r)$; $B(x,r)$: closed ball; $D=\sup_{x,y\in C}\|x-y\|$; $K$: iteration horizon; $\ell_f(v;z)=f(z)+\langle\nabla f(z),v-z\rangle$. Costs count each oracle invocation once.
\end{tablenotes}
\end{threeparttable}
\end{table}
\clearpage

\begin{table}[tp]
\centering
\caption{BFGS for linear-predictor objectives: \href{https://github.com/k3nfalt/longtails/blob/main/LongTails/Setups/LinearPredictorBFGS.lean}{LinearPredictorBFGS}.}\label{tab:LinearPredictorBFGS}
\parbox{\linewidth}{\footnotesize $f(x)=h(Ax)$ on $\mathbb{R}^d$, with linear $A:\mathbb{R}^d\to\mathbb{R}^m$, convex $L_h$-smooth $h$, $L\ge L_h\|A\|^2$, $L>0$, and $0<\mu\le L_h$. A minimizer $x_\star$ satisfies $\|x_0-x_\star\|\le R$ and $\frac{\mu}{2}\|A(x-x_\star)\|^2\le f(x)-f(x_\star)$ whenever $f(x)\le f(x_0)$. The setup restricts algorithms to BFGS with positive scalar initialization; \textbf{(1)} uses $B_0=LI$. The oracle returns an Armijo--weak-Wolfe step and its endpoint gradient, with $0<c_1<c_2<1$. Target: $f(x_K)-f(x_\star)\le\varepsilon$, $\varepsilon>0$, for every admissible oracle history. Bundled line-search calls include initialization; internal search trials have no stated gradient/value-call budget.}\par\smallskip
\begin{threeparttable}
\footnotesize
\setlength{\tabcolsep}{4pt}\renewcommand{\arraystretch}{1.25}
\begin{tabularx}{\linewidth}{@{}>{\raggedright\arraybackslash}p{0.20\linewidth}>{\raggedright\arraybackslash}p{0.32\linewidth}>{\raggedright\arraybackslash}p{0.25\linewidth}>{\raggedright\arraybackslash}X@{}}
\toprule
\textbf{Result} & \textbf{Oracle complexity} & \textbf{Novelty} & \textbf{Limitations} \\
\midrule
\textbf{(1)} \href{https://github.com/k3nfalt/longtails/blob/main/LongTails/Theorems/LinearPredictorBFGS.lean}{Linear\allowbreak{}Predictor\allowbreak{}BFGS} \textbf{(Ours)} $(\star\star\star)$ \href{https://github.com/k3nfalt/longtails/blob/main/scripts/reviews/LinearPredictorBFGS.md}{(review)}
& $\mathcal{O}\!\left(1+\kappa^{3/2}\gamma^{-3/2}\sqrt{LR^2/\varepsilon}\right)$ line-search calls.
& & \\
\midrule[1.2pt]
\textbf{(2)} BFGS \citep[Thm.~3.1, Eq.~(11)]{lpbDingYangZhou2026v1}
& $\mathcal{O}\!\left(1+\frac{LW_0^2}{\gamma\varepsilon}\right)$ line-search calls; $B_0=LI$, $W_0<\infty$.
& \textbf{(1)}: Distance-only BFGS bound, independent of design conditioning; $\varepsilon^{-1/2}$ instead of $\varepsilon^{-1}$. For $W_0\asymp R$, total gain when $LR^2/\varepsilon\gg\kappa^3/\gamma$. No uniform $W_0$ bound in setup parameters: $h(z)=\|z\|^2/2$, $A=\operatorname{diag}(1,\delta)$, $x_0=(R,0)$ give $W_0\ge R/(2\delta)$.
& \textbf{(1)} needs composite growth; \textbf{(2)} covers general convex smooth objectives with finite $W_0$. For $W_0\asymp R$, \textbf{(1)} loses at coarse accuracy and has worse $\gamma$ dependence. \\
\midrule
\textbf{(3)} AGD/FISTA \citep[Thm.~4.4]{ifsBeckTeboulle2009}, fixed step, $\psi=0$
& $\mathcal{O}\!\left(1+\sqrt{LR^2/\varepsilon}\right)$ gradient calls, available through $v=0$ queries; applies throughout the target objective class.
& \textbf{(1)} matches the $\varepsilon^{-1/2}$ order for BFGS, with extra parameter factors; no problem-class complexity improvement.
& \textbf{(1)} pays $\kappa^{3/2}\gamma^{-3/2}$ and a larger leading constant; needs composite growth and line-search access instead of gradients alone. \\
\bottomrule
\end{tabularx}
\begin{tablenotes}[flushleft]\footnotesize
\item[] \textbf{Abbr:} $\kappa=L_h/\mu$; $\gamma=c_1(1-c_2)$; $f_\star=f(x_\star)$; $X_\star=\operatorname{argmin} f$; $W_0=\sup_{x\notin X_\star,\,f(x)\le f(x_0)}(f(x)-f_\star)/\|\nabla f(x)\|$; $\psi$ is the FISTA nonsmooth term; $\delta>0$ is a design singular value; $\asymp$ means equality up to universal factors; $\gg$ means asymptotically larger.
\end{tablenotes}
\end{threeparttable}
\end{table}

\begin{table}[tp]
\centering
\caption{Relative gradient error: \href{https://github.com/k3nfalt/longtails/blob/main/LongTails/Setups/RelativeErrorSmoothConvex.lean}{RelativeErrorSmoothConvex}.}\label{tab:RelativeErrorSmoothConvex}
\parbox{\linewidth}{\footnotesize Convex differentiable $f$ on a real Hilbert space satisfies $f(y)\le f(x)+\langle\nabla f(x),y-x\rangle+\frac L2\|y-x\|^2$, with $L>0$, a minimizer $x^\star$, and $\|x_0-x^\star\|\le R$, $R\ge0$. Every oracle reply satisfies $\|g-\nabla f(x)\|\le\alpha\|\nabla f(x)\|$, $0\le\alpha<1$, for every call and seed. Target: $f(x_N)-f(x^\star)\le\varepsilon$, $\varepsilon>0$, almost surely under the execution measure (indeed, for every seed). Parameters $L,\alpha,R,\varepsilon$ are supplied.}\par\smallskip
\begin{threeparttable}
\footnotesize
\setlength{\tabcolsep}{4pt}\renewcommand{\arraystretch}{1.25}
\begin{tabularx}{\linewidth}{@{}>{\raggedright\arraybackslash}p{0.20\linewidth}>{\raggedright\arraybackslash}p{0.32\linewidth}>{\raggedright\arraybackslash}p{0.25\linewidth}>{\raggedright\arraybackslash}X@{}}
\toprule
\textbf{Result} & \textbf{Oracle complexity} & \textbf{Novelty} & \textbf{Limitations} \\
\midrule
\textbf{(1)} \href{https://github.com/k3nfalt/longtails/blob/main/LongTails/Theorems/AcceleratedRelativeErrorConvex.lean}{Accelerated\allowbreak Relative\allowbreak Error\allowbreak Convex} \textbf{(Ours)} $(\star\star\star)$ \href{https://github.com/k3nfalt/longtails/blob/main/scripts/reviews/AcceleratedRelativeErrorConvex.md}{(review)}
& $\mathcal{O}\!\left(1+\sqrt{\frac{T}{1-\alpha}}+\frac{\alpha^2T}{1-\alpha}\right)$ & & \\
\midrule[1.2pt]
\textbf{(2)} Regularized RE-AGM \citep[Thm.~7.2]{igsVasinGasnikov2026v1}
& $\mathcal{O}\!\left((\sqrt T+\alpha T)(1+\log T)\right)$.
$T>1$, $\alpha\le1/3$; Euclidean space.
& \textbf{(1)}: $\alpha^2T$ instead of $\alpha T$; removes logarithm on $T>1$, $\alpha\le1/3$. At $\alpha=T^{-1/4}$, $T\ge81$: $\mathcal{O}(T^{3/4}(1+\log T))\to\mathcal{O}(\sqrt T)$. Scope: every $\alpha<1$.
& \\
\midrule
\textbf{(3)} Regularized GD \citep[Thm.~7.1]{igsVasinGasnikov2026v1}
& $\mathcal{O}\!\left(\frac{(1+\alpha)^2}{(1-\alpha)^6}T(1+\log T)\right)$.
$T>1$, $\alpha<1/2$; Euclidean space.
& \textbf{(1)}: $\mathcal{O}(\sqrt T+\alpha^2T)$ on $T>1$, $\alpha<1/2$; logarithmic saving throughout. For $\alpha\le T^{-1/4}$: $\mathcal{O}(T(1+\log T))\to\mathcal{O}(\sqrt T)$. Scope: every $\alpha<1$.
& \\
\midrule
\textbf{(4)} ISTM, $p=2$ \citep[Thm.~3.1]{rescKornilov2023v1}
& $\mathcal{O}(1+\sqrt T)$ when $\alpha\le c\min\{1,T^{-1/2}\}$; Euclidean space. Certified accuracy floor of order $\alpha^2LR^2$.
& \textbf{(1)}: matches $\mathcal{O}(\sqrt T)$ for $T>1$, $\alpha\lesssim T^{-1/2}$; extends acceleration to $\alpha\lesssim T^{-1/4}$ with $\alpha\le1/3$. Arbitrary accuracy for every $\alpha<1$.
& \\
\bottomrule
\end{tabularx}
\begin{tablenotes}[flushleft]\footnotesize
\item[] \textbf{Abbr:} $T=LR^2/\varepsilon$; complexity counts relative-error gradient calls; $c>0$ is a sufficiently small universal constant; $\lesssim$ suppresses a universal constant. GD: gradient descent; RE-AGM: relative-error accelerated gradient method; ISTM: intermediate similar triangle method.
\end{tablenotes}
\end{threeparttable}
\end{table}
\clearpage

\begin{table}[tp]
\centering
\caption{Gradient sliding for \href{https://github.com/k3nfalt/longtails/blob/main/LongTails/Setups/SlidingGeneralizedSmooth.lean}{SlidingGeneralizedSmooth}.}\label{tab:SlidingGeneralizedSmooth}
\parbox{\linewidth}{\footnotesize Minimize $F=f+h$ over convex $Q$ in a real Hilbert space. The convex differentiable $f$ has continuous gradient and satisfies $\|\nabla f(x)-\nabla f(y)\|\le [L_0+L_1\sup_{z\in[x,y]}\|\nabla f(z)\|]\|x-y\|$ on the whole space. Oracles supply exact gradients, exact projections onto $Q$, and possibly adversarial vectors $s(u)$ with $\|s(u)\|\le M$ and $h(v)\ge h(u)+\langle s(u),v-u\rangle$ for $u,v\in Q$. There is $x^\star\in\arg\min_Q F$ with $x_0\in Q$, $\|x_0-x^\star\|\le R$, $F(x_0)-F(x^\star)\le\Delta$, and $\|\nabla f(x^\star)\|\le G$. Here $L_0,R,\Delta,\varepsilon>0$ and $L_1,M,G\ge0$. Target: a feasible output $x$ with $F(x)-F(x^\star)\le\varepsilon$ almost surely.}\par\smallskip
\begin{threeparttable}
\footnotesize
\setlength{\tabcolsep}{4pt}\renewcommand{\arraystretch}{1.25}
\begin{tabularx}{\linewidth}{@{}>{\raggedright\arraybackslash}p{0.20\linewidth}>{\raggedright\arraybackslash}p{0.32\linewidth}>{\raggedright\arraybackslash}p{0.25\linewidth}>{\raggedright\arraybackslash}X@{}}
\toprule
\textbf{Result} & \textbf{Oracle complexity} & \textbf{Novelty} & \textbf{Limitations} \\
\midrule
\textbf{(1)} \href{https://github.com/k3nfalt/longtails/blob/main/LongTails/Theorems/GeneralizedGradientSliding.lean}{Generalized\allowbreak{}Gradient\allowbreak{}Sliding} \textbf{(Ours)} $(\star\star\star)$ \href{https://github.com/k3nfalt/longtails/blob/main/scripts/reviews/GeneralizedGradientSliding.md}{(review)}
& Total: $\mathcal{O}\!\left(\begin{gathered}\sqrt{AR^2/\varepsilon}+(1+L_1R)\ell\\{}+M^2R^2/\varepsilon^2\end{gathered}\right)$.
\par Gradients: $\mathcal{O}\!\left(\sqrt{\frac{AR^2}{\varepsilon}}+(1+L_1R)\ell\right)$.\par $G$ is used only in the gradient-count cap.
& & \\
\midrule[1.2pt]
\textbf{(2)} Gradient sliding \citep{sgsLan2016}, Cor.~1(a), arXiv v2. Derived bound: finite-dimensional localization to $X=Q\cap\overline B(x_0,2R)$.
& Total: $\mathcal{O}\!\left(1+\sqrt{\frac{Ae^{3L_1R}R^2}{\varepsilon}}+\frac{M^2R^2}{\varepsilon^2}\right)$.
\par Gradients: $\mathcal{O}\!\left(1+\sqrt{\frac{Ae^{3L_1R}R^2}{\varepsilon}}\right)$.
\par Uses $L=Ae^{3L_1R}$ and exact $\Pi_X$; at $L_1=0$, use $X=Q$, $L=L_0$.
& \textbf{(1)}: gradient saving $e^{3L_1R/2}$ when $L_1R\gg1$ and $\sqrt{AR^2/\varepsilon}\gg(1+L_1R)\ell$; total saving additionally requires $M^2R^2/\varepsilon^2\ll\sqrt{AR^2/\varepsilon}$. Only $\Pi_Q$ needed. Matches leading rate for $L_1R\lesssim1$; scope extends global smoothness.
& \textbf{(1)}: extra $\ell$ term at $L_1=0$; no Bregman geometry or additional simple prox term. \\
\midrule
\textbf{(3)} Normalized subgradient \citep{sgsGrimmer2019}, Thm.~2, arXiv v3. Derived bound, $E=\mathbb R^d$; best iterate selected using $F$ values.
& Total and gradients:
$\mathcal{O}\!\left(\begin{gathered}1+(G+M)^2R^2/\varepsilon^2\\{}+AR^2/\varepsilon+L_1^2R^2\end{gathered}\right)$.
\par Includes $\nabla f$, $h'$ and $\Pi_Q$ calls, plus $F$ values for selection; $\Pi_Q$ is unnecessary when $Q=E$.
& \textbf{(1)}: gradients $\varepsilon^{-1/2}$ instead of $\varepsilon^{-2}$ for fixed parameters, $G+M>0$; no $F$ values. Total gains can occur for constrained $G\gg M$. For $Q=E$, take $G=M$: matched total order unless $\varepsilon A\gg M^2$ or $\varepsilon L_1\gg M$; initialization can mask gains.
& \textbf{(1)}: updates additionally need $L_0,L_1,M,\Delta,\varepsilon$; \textbf{(3)} needs only $R,T$. Initialization term $(1+L_1R)\ell$ exceeds $L_1^2R^2$ when $\ell\gg L_1R\ge1$. \\
\bottomrule
\end{tabularx}
\begin{tablenotes}[flushleft]\footnotesize
\item[] \textbf{Abbr:} $A=L_0+L_1G$; $\ell=1+\log\!\left(1+(1+L_1R)\sqrt{\Delta/(L_0R^2)}\right)$; $\Pi_C$ is exact Euclidean projection onto $C$; $h'$ denotes a returned subgradient; $T$ is the iteration horizon in \textbf{(3)}. Total cost counts all primitive oracle calls.
\end{tablenotes}
\end{threeparttable}
\end{table}
\clearpage

\begin{table}[tp]
\centering
\caption{Stochastic generalized monotone root finding: \href{https://github.com/k3nfalt/longtails/blob/main/LongTails/Setups/StochasticGeneralizedMonotone.lean}{StochasticGeneralizedMonotone}.}\label{tab:StochasticGeneralizedMonotone}
\parbox{\linewidth}{\footnotesize On a finite-dimensional real Hilbert space, $F$ is continuous and monotone, $F(x^\star)=0$, $\|x_0-x^\star\|\le R$, and $\|F(x)-F(y)\|\le(L_0+L_1\sup_{z\in[x,y]}\|F(z)\|)\|x-y\|$. The jointly measurable oracle satisfies $\int g(x,\omega)\,dp=F(x)$ and $\int\|g(x,\omega)-F(x)\|^2\,dp\le\sigma^2$, with integrable values and squared errors, and IID samples of law $p$ under $P$. Here $L_0,L_1,R,\sigma\ge0$, $\varepsilon>0$, $0<\delta<1$. Target: measurable output $X$ with $P\{\|F(X)\|>\varepsilon\}\le\delta$.}\par\smallskip
\begin{threeparttable}
\footnotesize
\setlength{\tabcolsep}{4pt}\renewcommand{\arraystretch}{1.25}
\begin{tabularx}{\linewidth}{@{}>{\raggedright\arraybackslash}p{0.20\linewidth}>{\raggedright\arraybackslash}p{0.32\linewidth}>{\raggedright\arraybackslash}p{0.25\linewidth}>{\raggedright\arraybackslash}X@{}}
\toprule
\textbf{Result} & \textbf{Oracle complexity} & \textbf{Novelty} & \textbf{Limitations} \\
\midrule
\textbf{(1)} \href{https://github.com/k3nfalt/longtails/blob/main/LongTails/Theorems/StochasticGeneralizedMonotone.lean}{Stochastic\allowbreak Generalized\allowbreak Monotone} \textbf{(Ours)} $(\star\star\star)$ \href{https://github.com/k3nfalt/longtails/blob/main/scripts/reviews/StochasticGeneralizedMonotone.md}{(review)}
& $\mathcal{O}\!\left(\left[1+\frac{L_0R}{\varepsilon}+(L_1R)^2+\frac{\sigma^2}{\varepsilon^2}\right]\right.$\newline
$\left.{}\cdot\log^3 Q\,\log(Q/\delta)\right)$ stochastic calls, including batches and repetitions.
& & \\
\midrule[1.2pt]
\textbf{(2)} \citep{sgmJiLanZhu2026v1}, Thm.~2.2; globally Lipschitz special case $L_1=0$, $L_0>0$
& $\mathcal{O}\!\left(\frac{L_0R}{\varepsilon}h+\frac{\sigma^2}{\varepsilon^2}\ell^3h^2\right)$ stochastic calls; $\varepsilon\le L_0R$; $\bigl(\mathbb E\|F(X)\|^2\bigr)^{1/2}\le\varepsilon$ on the unconstrained restriction.
& \textbf{(1)}: Scope extends to $L_1>0$ with polynomial dependence and additive variance cost. At $L_1=0$, matches polynomial accuracy dependence; gives high-probability accuracy under finite variance.
& \textbf{(1)}: Extra logarithmic overhead at $L_1=0$, including confidence dependence; no constrained-VI guarantee. \\
\midrule
\textbf{(3)} \citep{sgmVankov2025v2}, Thm.~4.4, $p=2$; quasi-sharp, $0<\alpha\le\tfrac12$
& $\mathcal{O}\!\left(1+\frac{dR}{a\sqrt\tau}+\frac{\sigma^2+K_1\sigma^{2\alpha}}{a^2\tau}\right)$ stochastic calls for $\mathbb E\,\mathrm{dist}^2(X,X^\star)\le\tau$; two calls per iteration.
& \textbf{(1)}: No quasi-sharpness assumption; permits $\alpha=1$, including monotone operators with $\langle F(x),x-x^\star\rangle=0$ near roots. Residual and distance rates are not comparable.
& \textbf{(1)}: Does not cover constrained or structured non-monotone problems covered by \textbf{(3)}. \\
\midrule
\textbf{(4)} Anchored extragradient \citep{cgmRepositoryAnchored}; exact-oracle special case $\sigma=0$, $L_0>0$
& $\mathcal{O}\!\left(1+\frac{L_0R}{\varepsilon}+(1+L_1R)\right.$\newline
$\left.{}\cdot\log\!\left(1+\frac{(L_0+L_1G)R}{\varepsilon}\right)\right)$ exact calls; deterministic residual accuracy.
& \textbf{(1)}: Scope extends to noisy finite-variance oracles; no initial-residual bound $G$ needed for its budget.
& \textbf{(1)}: At $\sigma=0$, quadratic rather than linear $L_1R$ term, with additional logarithms; loses anytime last-iterate guarantee. \\
\bottomrule
\end{tabularx}
\begin{tablenotes}[flushleft]\footnotesize
\item[] \textbf{Abbr:} $Q=2+L_0R/\varepsilon+L_1R$, $\ell=1+\log(L_0R/\varepsilon)$ and $h=1+\log\ell$ in \textbf{(2)}; $G\ge\|F(x_0)\|$ in \textbf{(4)}. In \textbf{(3)}, $X^\star$ is the solution set, $\tau>0$ is squared-distance accuracy, $R^2\ge\mathbb E\,\mathrm{dist}^2(x_0,X^\star)$, and quasi-sharpness means $\langle F(x),x-x^\star\rangle\ge\mu\,\mathrm{dist}^2(x,X^\star)$, $\mu>0$. The exponent $\alpha$ replaces $1$ on the segment operator norm in generalized Lipschitz continuity; its constants are denoted $\bar L_0,\bar L_1$. Set $c_\alpha=2^{\alpha^2/(1-\alpha)}$, $K_0=\bar L_0(c_\alpha+1)$, $K_1=\bar L_1c_\alpha3^\alpha$, $K_2=\bar L_1^{1/(1-\alpha)}c_\alpha3^\alpha(1-\alpha)^{\alpha/(1-\alpha)}$, $a=\mu\min\{1,[2(C_F+\sigma)]^{-1}\}$, $d=\max\{4\mu,2\sqrt3(K_0+K_1+K_2)\}$, where $C_F\ge\sup_k\mathbb E\|F(x_k)\|$.
\end{tablenotes}
\end{threeparttable}
\end{table}
\clearpage

\begin{table}[tp]
\centering
\caption{Stochastic minimax optimization with Lipschitz Hessians: \href{https://github.com/k3nfalt/longtails/blob/main/LongTails/Setups/StochasticHessianNcsc.lean}{StochasticHessianNcsc}.}\label{tab:StochasticHessianNcsc}
\parbox{\linewidth}{\footnotesize On finite-dimensional real inner-product spaces, $f(x,y)$ is jointly $L$-smooth, $\mu$-strongly concave in $y$, and has $\rho$-Lipschitz Hessian $\nabla^2 f$. Let $\Phi(x)=\max_y f(x,y)$, $\Phi(x_0)-\inf_x\Phi(x)\le\Delta$, and $\|\nabla_y f(x_0,y_0)\|\le G$. Here $0<\mu\le L$, $\rho,\sigma,\varepsilon>0$, and $\Delta,G\ge0$. Each call returns a jointly measurable, unbiased joint stochastic gradient with total variance at most $\sigma^2$, using a fresh IID sample. The target is $\int\|\nabla\Phi(\widehat x(\omega))\|\,d\mathbb P(\omega)\le\varepsilon$, where $\mathbb P$ includes the sampling and algorithm randomness.}\par\smallskip
\begin{threeparttable}
\footnotesize
\setlength{\tabcolsep}{4pt}\renewcommand{\arraystretch}{1.25}
\begin{tabularx}{\linewidth}{@{}>{\raggedright\arraybackslash}p{0.20\linewidth}>{\raggedright\arraybackslash}p{0.32\linewidth}>{\raggedright\arraybackslash}p{0.25\linewidth}>{\raggedright\arraybackslash}X@{}}
\toprule
\textbf{Result} & \textbf{Oracle complexity} & \textbf{Novelty} & \textbf{Limitations} \\
\midrule
\textbf{(1)} \href{https://github.com/k3nfalt/longtails/blob/main/LongTails/Theorems/HessianNcscImplicitTransport.lean}{Hessian\allowbreak Ncsc\allowbreak Implicit\allowbreak Transport} \textbf{(Ours)} $(\star\star\star)$ \href{https://github.com/k3nfalt/longtails/blob/main/scripts/reviews/HessianNcscImplicitTransport.md}{(review)}
& $\begin{aligned}[t]
\mathcal{O}\!\Bigl(&K^2(1+\sigma^2/\varepsilon^2)\\
&\cdot(1+\rho\varepsilon K/L^2)\\
&\cdot\bigl[K\Delta\sqrt Q/\varepsilon+1+\sigma/\varepsilon\\
&\qquad+\log(1+G^2Q/\mu^2)\bigr]\Bigr)
\end{aligned}$ joint-gradient calls, including initialization.
& & \\
\midrule[1.2pt]
\textbf{(2)} \citep[Thm.~3.4, v1]{shncscWu2026v1}
& $\begin{aligned}[t]
\mathcal{O}\!\Bigl(\Lambda^{10}\bigl[&\sqrt\kappa L\Delta/\varepsilon^2\\
&+\kappa L\Delta\sigma^2/\varepsilon^4\\
&+\kappa^2\sigma^2/\varepsilon^2\\
&+\sqrt\kappa\bigl(1+\log(1+\\
&\quad\kappa^2(G^2+L\Delta)/\varepsilon^2)\bigr)\bigr]\Bigr)
\end{aligned}$ joint-gradient calls; success probability $2/3$. Applies for $\kappa\ge\kappa_0$, $L\Delta\ge C_0\varepsilon^2$ (universal thresholds).
& \textbf{(1)}: $\varepsilon^{-7/2}$ instead of $\varepsilon^{-4}$ for fixed positive parameters as $\varepsilon\to0$. More generally, gain up to logarithms when $\sigma\ge\varepsilon$ and
$\varepsilon\ll\min\{L^2/(\rho\kappa^5),\sigma/\kappa^2,L\Delta/(\kappa\sigma)\}$.
Expected stationarity implies the shared probability target by tuning to $\varepsilon/3$.
& \textbf{(1)} needs a known Hessian-Lipschitz bound; leading noise dependence $\kappa^{7/2}$ instead of $\kappa$. Low-noise gap term worse by $\kappa^{5/2}$ when $\rho\varepsilon K\le L^2$; larger proved constants. \\
\midrule
\textbf{(3)} HCMM-1 \citep[Thm.~1, v3]{shncscCai2026v3}
& $\mathcal{O}(\varepsilon^{-3})$ stochastic-gradient and stochastic Hessian--vector-product calls in total, for fixed remaining parameters and initialization; their dependence is unspecified. Requires bounded gradients and unbiased, bounded-variance stochastic Hessians.
& \textbf{(1)}: no stochastic Hessian oracle or bounded-gradient assumption. Rates are not comparable under the supplied gradient-only oracle.
& \textbf{(1)} has exponent $7/2$ instead of $3$ with other parameters fixed, under the respective oracle models. \\
\midrule
\textbf{(4)} NIGT \citep[Thm.~3]{shncscCutkosky2020}
& Derived bound:
$\begin{aligned}[t]
\mathcal{O}\!\Bigl(&1+L\Delta/\varepsilon^2\\
&+\Delta\sigma^2\sqrt\rho/\varepsilon^{7/2}\\
&+\frac{\sigma^{13/3}}{\Delta^{4/3}\rho^{2/3}\varepsilon^{7/3}}\Bigr)
\end{aligned}$ gradient calls, $\Delta>0$. Special case $f(x,y)=F(x)-\frac\mu2\|y\|^2$; shift $\inf F$ to zero.
& \textbf{(1)} extends single-level optimization to coupled minimax problems. Matches the leading $\varepsilon^{-7/2}$ rate on the overlap ($\mu=L$, fixed positive parameters).
& \textbf{(1)} retains extra $L\Delta\sigma/\varepsilon^3$-type terms even at $\mu=L$. \\
\bottomrule
\end{tabularx}
\begin{tablenotes}[flushleft]\footnotesize
\item[] \textbf{Abbr:} $\kappa=L/\mu$, $K=1+\kappa$, $Q=\rho K/\varepsilon+L^2/(\varepsilon^2+\sigma^2)$;
$\Lambda=1+\log(1+\kappa)+\log(1+L\Delta/\varepsilon^2)+\log(1+\kappa^2\sigma^2/\varepsilon^2)$.
IID: independent, identically distributed. One joint-gradient call returns the $x$ and $y$ components together.
\end{tablenotes}
\end{threeparttable}
\end{table}
\clearpage

\begin{table}[tp]
\centering
\caption{High-probability convex optimization under distance-growing second moments: \href{https://github.com/k3nfalt/longtails/blob/main/LongTails/Setups/TerminalClippingConvex.lean}{TerminalClippingConvex}.}\label{tab:TerminalClippingConvex}
\parbox{\linewidth}{\footnotesize Measurable $f$ is convex on closed convex $X\subseteq\mathbb R^d$; $x_0,x^\star\in X$, $x^\star\in\arg\min_X f$, and $\|x_0-x^\star\|\le R$. IID samples have law $\nu$ and joint law $P$. The jointly measurable oracle satisfies $\int g(x,\omega)\,d\nu(\omega)\in\partial_X f(x)$ and $\int\|g(x,\omega)\|^2\,d\nu(\omega)\le B^2\|x-x_0\|^2+G^2$ on $X$; exact measurable projection onto $X$ is available. Here $B,G,\varepsilon>0$, $R\ge0$, $0<\delta<1$. The setup restricts algorithms to projected clipped subgradient updates with the current iterate as output. Target: $P\{f(x_T)-f(x^\star)>\varepsilon\}\le\delta$.}\par\smallskip
\begin{threeparttable}
\footnotesize
\setlength{\tabcolsep}{4pt}\renewcommand{\arraystretch}{1.25}
\begin{tabularx}{\linewidth}{@{}>{\raggedright\arraybackslash}p{0.20\linewidth}>{\raggedright\arraybackslash}p{0.32\linewidth}>{\raggedright\arraybackslash}p{0.23\linewidth}>{\raggedright\arraybackslash}X@{}}
\toprule
\textbf{Result} & \textbf{Oracle complexity} & \textbf{Novelty} & \textbf{Limitations} \\
\midrule
\textbf{(1)} \href{https://github.com/k3nfalt/longtails/blob/main/LongTails/Theorems/ClippedBgLastIterate.lean}{Clipped\allowbreak Bg\allowbreak Last\allowbreak Iterate} \textbf{(Ours)} $(\star\star\star)$ \href{https://github.com/k3nfalt/longtails/blob/main/scripts/reviews/ClippedBgLastIterate.md}{(review)}
& $\mathcal{O}\!\left(1+\dfrac{(BR^2+G^2/B)^2\ell^2}{\varepsilon^2}\right)$ & & \\
\midrule[1.2pt]
\textbf{(2)} Projected Halpern \citep[Thm.~3.6]{tcAMW2026v3}; derived bound via Markov
& $\mathcal{O}\!\left(1+\dfrac{A^2}{\varepsilon^2\delta^2}\left[1+\log\!\left(1+\dfrac{A}{\varepsilon\delta}\right)\right]^2\right)$
& \textbf{(1)}: On the same moment class, logarithmic instead of polynomial confidence cost; removes the accuracy logarithm at fixed confidence. Comparison is specifically with Markov conversion.
& \textbf{(1)} needs a fixed horizon instead of an anytime schedule; substantially larger proved numerical constants. \\
\midrule
\textbf{(3)} Clipped SGD \citep[Thm.~3]{tcLZ2023v3}, $p=2$
& $\mathcal{O}\!\left(1+\dfrac{(\beta+R^2/\beta)^2(M^2+\sigma^2)\ell^2}{\varepsilon^2}\right)$\par Special case: $M$-bounded subgradients, uniform variance $\sigma^2$; $\beta>0$, known $\sigma$, fixed horizon.
& \textbf{(1)}: Matches this rate for $G^2=M^2+\sigma^2$, $B=G/\beta$, with last iterate instead of averaging. Scope: allows non-Lipschitz objectives and growing moments.
& \textbf{(1)}: The $G^4/(B^2\varepsilon^2)$ term diverges as $B\downarrow0$; bounded-moment tuning $B=G/R$ needs $R>0$ known. \\
\midrule
\textbf{(4)} Distance-adaptive clipped SGD \citep[Thm.~6]{tcLZ2023v3}, $p=2$
& $\mathcal{O}\!\left(h+\dfrac{aM\ell}{\varepsilon}+z^2[1+\log(1+z)]^2\right)$\par Special case: $M$-bounded subgradients, known uniform variance bound $\sigma^2$; $r>0$, $w_t=1+\log^2t$.
& \textbf{(1)}: Last iterate instead of a selected weighted average; extends to non-Lipschitz objectives with growing moments.
& \textbf{(1)} is not distance-adaptive: $(\beta+R^2/\beta)^2$ mismatch when $B=G/\beta$. Loses anytime tuning. \\

\bottomrule
\end{tabularx}
\begin{tablenotes}[flushleft]\footnotesize
\item[] \textbf{Abbr:} $\partial_X f(x)=\{v:f(y)\ge f(x)+\langle v,y-x\rangle\ \forall y\in X\}$; $\ell=1+\log(1/\delta)$; $A=BR^2+G^2/B$. Costs count subgradient calls plus projection calls (one of each per update); on $X=\mathbb R^d$ projection is the identity. $M$ bounds deterministic subgradients; $\sigma$ is the central-noise scale; $L$ is gradient Lipschitzness; $p,\alpha$ are central moment orders. In \textbf{(4)}, $h=1+\log(1+R/r)$, $a=(r+R)h$, $z=a\sqrt{M^2+\sigma^2\ell}/\varepsilon$. In \textbf{(6)}, $D=\max\{R\sigma,L^{1/5}R^{6/5}\sigma^{4/5},L^{1/3}R^{4/3}\sigma^{2/3}\}$ and $\ell_s=1+\log((1+LR^2/\varepsilon+D^3/\varepsilon^3)/\delta)$.
\end{tablenotes}
\end{threeparttable}
\end{table}
\clearpage
\begin{table}[tp]
\centering
\caption{TerminalClippingConvex (continued).}
\begin{threeparttable}
\footnotesize
\setlength{\tabcolsep}{4pt}\renewcommand{\arraystretch}{1.25}
\begin{tabularx}{\linewidth}{@{}>{\raggedright\arraybackslash}p{0.20\linewidth}>{\raggedright\arraybackslash}p{0.32\linewidth}>{\raggedright\arraybackslash}p{0.23\linewidth}>{\raggedright\arraybackslash}X@{}}
\toprule
\textbf{Result} & \textbf{Oracle complexity} & \textbf{Novelty} & \textbf{Limitations} \\
\midrule
\textbf{(5)} SGD \citep[Thm.~3.5]{tcLZ2026v4}
& $\mathcal{O}\!\left(1+\dfrac{(M^2+\sigma^2\ell)R^2}{\varepsilon^2}\right)$\par Special case: $M$-Lipschitz objective, sub-Gaussian noise of scale $\sigma$; $R$-tuned.
& \textbf{(1)}: Extends the high-probability last-iterate guarantee to finite second moments and non-Lipschitz objectives with growing moments.
& On the overlap, $G^2=M^2+\sigma^2$, $B=G/R$, $R>0$: \textbf{(1)} costs $\mathcal{O}\!\left(1+\dfrac{(M^2+\sigma^2)R^2\ell^2}{\varepsilon^2}\right)$, losing up to $\ell^2$. \\
\midrule
\textbf{(6)} Clipped SGD \citep[Thm.~1, Cor.~1]{tcCPPG2026}, $\alpha=2$
& $\mathcal{O}\!\left(1+\dfrac{LR^2\ell_s^2}{\varepsilon}+\dfrac{D^3\ell_s^6}{\varepsilon^3}\right)$\par Special case: $X=\mathbb R^d$, $L$-smooth, uniform variance $\sigma^2$; $L,R,\sigma>0$.
& \textbf{(1)}: With $B=\sqrt2L$, $G^2=2L^2R^2+\sigma^2$, costs $\mathcal{O}\!\left(1+\dfrac{(LR^2+\sigma^2/L)^2\ell^2}{\varepsilon^2}\right)$; improves the accuracy dependence when $\varepsilon\le\min\{\sigma^3/(L^2R),L^2R^3/\sigma\}$. Scope: no smoothness needed. Updates can use supplied $G^2\ge2\|\nabla f(x_0)\|^2+\sigma^2$ instead of $R$.
& \textbf{(1)} loses anytime tuning and coverage of $\alpha<2$; worse deterministic term: $L^2R^4/\varepsilon^2$ instead of $LR^2/\varepsilon$, apart from logarithms. \\
\bottomrule
\end{tabularx}
\begin{tablenotes}[flushleft]\footnotesize
\item[] \textbf{Abbr:} $\partial_X f(x)=\{v:f(y)\ge f(x)+\langle v,y-x\rangle\ \forall y\in X\}$; $\ell=1+\log(1/\delta)$; $A=BR^2+G^2/B$. Costs count subgradient calls plus projection calls (one of each per update); on $X=\mathbb R^d$ projection is the identity. $M$ bounds deterministic subgradients; $\sigma$ is the central-noise scale; $L$ is gradient Lipschitzness; $p,\alpha$ are central moment orders. In \textbf{(4)}, $h=1+\log(1+R/r)$, $a=(r+R)h$, $z=a\sqrt{M^2+\sigma^2\ell}/\varepsilon$. In \textbf{(6)}, $D=\max\{R\sigma,L^{1/5}R^{6/5}\sigma^{4/5},L^{1/3}R^{4/3}\sigma^{2/3}\}$ and $\ell_s=1+\log((1+LR^2/\varepsilon+D^3/\varepsilon^3)/\delta)$.
\end{tablenotes}
\end{threeparttable}
\end{table}
\clearpage

\begin{table}[tp]
\centering
\caption{Third-derivative-smooth nonconvex optimization: \href{https://github.com/k3nfalt/longtails/blob/main/LongTails/Setups/ThirdDerivativeSmoothNonconvex.lean}{ThirdDerivativeSmoothNonconvex}.}\label{tab:ThirdDerivativeSmoothNonconvex}
\parbox{\linewidth}{\footnotesize On a finite-dimensional real inner-product space, $f\in C^3$ has $L_3$-Lipschitz third derivative, $f(x_0)-\inf f\le\Delta$, and $\|\nabla^2f(x_0)\|\le H_0$, where $L_3,\Delta,\varepsilon>0$ and $H_0\ge0$. Exact value/gradient and full-Hessian oracles are available. The target is $\|\nabla f(x_{\rm out})\|\le\varepsilon$ almost surely; \textbf{(1)} guarantees it for every seed. No global gradient- or Hessian-Lipschitz bound or algorithm restriction is imposed.}\par\smallskip
\begin{threeparttable}
\footnotesize
\setlength{\tabcolsep}{4pt}\renewcommand{\arraystretch}{1.25}
\begin{tabularx}{\linewidth}{@{}>{\raggedright\arraybackslash}p{0.20\linewidth}>{\raggedright\arraybackslash}p{0.32\linewidth}>{\raggedright\arraybackslash}p{0.25\linewidth}>{\raggedright\arraybackslash}X@{}}
\toprule
\textbf{Result} & \textbf{Oracle complexity} & \textbf{Novelty} & \textbf{Limitations} \\
\midrule
\textbf{(1)} \href{https://github.com/k3nfalt/longtails/blob/main/LongTails/Theorems/LazyQuarticHessianTradeoff.lean}{Lazy\allowbreak Quartic\allowbreak Hessian\allowbreak Tradeoff} \textbf{(Ours)} $(\star\star\star)$ \href{https://github.com/k3nfalt/longtails/blob/main/scripts/reviews/LazyQuarticHessianTradeoff.md}{(review)}
& $\mathcal{O}\!\left((1+S^{5/4})\ell+m(1+S^{3/4})\right)$ total; $\mathcal{O}(1+S^{3/4})$ Hessian calls.
& & \\
\midrule[1.2pt]
\textbf{(2)} \citep{tdsDoikovGrapiglia2025v1}, Thm.~3.5, Cor.~3.6 ($p=3$)
& $\mathcal{O}\!\left(m\left[d^{1/3}\Delta\bar L^{1/3}\varepsilon^{-4/3}\right.\right.$\newline
$\left.\left.{}+d\{1+\log_+(L_3/L_{\rm init})\}\right]\right)$ total, counting trials and finite differences; $L_{\rm init}>0$.
& \textbf{(1)}: Dimension-free; for $S\ge1$, $L_{\rm init}\asymp L_3$, Hessians $S^{3/4}$ instead of $d^{1/3}S+d$; total gain when $md^{1/3}\gg S^{1/4}\ell$.
\newline Example: $m=d$, $H_0\le\beta$, $L_{\rm init}\asymp L_3$: $\mathcal{O}(d^{4/3}S+d^2)$ becomes $\mathcal{O}(S^{5/4}[1+\log S]+dS^{3/4})$.
& \textbf{(1)} needs $L_3,\Delta,H_0$; \textbf{(2)} adapts to $L_3$. For $m=d=1$, $L_{\rm init}\asymp L_3$, $S\ge1$: $\mathcal{O}(S)$ becomes $\mathcal{O}(S^{5/4}\ell)$. \\
\midrule
\textbf{(3)} \citep{tdsHinder2018v2}, first-order variant after (8); globally $L_1$-Lipschitz gradient
& $\mathcal{O}\!\left(d(1+S)\left[1+\log_+\!\left(\frac{L_1^3}{L_3\varepsilon^2}\right)\right]\right)$ expected value/gradient calls; no Hessians.
& \textbf{(1)}: Scope extends beyond global $L_1$; dimension-free. On the overlap, for $S\ge1$, gains when $d\gg S^{1/4}$ and $m\ll dS^{1/4}$, up to logarithms.
& \textbf{(1)} needs Hessians; larger cost when $d\lesssim S^{1/4}$, up to logarithms. \\
\midrule
\textbf{(4)} \citep{tdsCarmonEtAl2017}, Thm.~2; globally $L_1$-Lipschitz gradient
& $\mathcal{O}\!\left(\frac{\Delta L_1^{1/2}L_3^{1/6}}{\varepsilon^{5/3}}\left[1+\log_+\!\left(\frac{L_1\Delta}{\varepsilon^2}\right)\right]\right)$ value/gradient calls; no Hessians.
\newline Valid for $\varepsilon^{2/3}\le\delta^{1/2}L_3^{1/6}$ and $8L_3^{1/3}\varepsilon^{2/3}\le L_1$.
& \textbf{(1)}: Scope extends beyond global $L_1$, including quartic least squares. On the overlap, replaces $\sqrt{L_1}$ by $(L_3\Delta)^{1/4}$; total gain for $L_1\gg\beta$ and $m\ll\Delta^{1/4}L_1^{1/2}L_3^{-1/12}\varepsilon^{-2/3}$, up to logarithms.
& \textbf{(1)} needs Hessians and $H_0$; larger gradient prefactor when $L_1<\beta$ within \textbf{(4)}'s admissible accuracy range. \\
\bottomrule
\end{tabularx}
\begin{tablenotes}[flushleft]\footnotesize
\item[] \textbf{Abbr:} $d=\dim E$; $m\ge1$ is the cost of one full Hessian in value/gradient-call units; total cost is $G+mH$, with $G,H$ the respective call counts. $S=\Delta L_3^{1/3}\varepsilon^{-4/3}$; $\beta=\sqrt{L_3\Delta}$; $\ell=1+\log(1+H_0\Delta/\varepsilon^2+L_3\Delta^3/\varepsilon^4)$. $L_1$ is a global gradient-Lipschitz constant; $L_{\rm init}$ is \textbf{(2)}'s initial smoothness estimate and $\bar L=\max\{L_{\rm init},L_3\}$. $\delta=f(x_0)-\inf f\le\Delta$; $\log_+t=\max\{0,\log t\}$; $\asymp$ denotes equality up to universal factors.
\end{tablenotes}
\end{threeparttable}
\end{table}
\clearpage

\begin{table}[tp]
\centering
\caption{Exact-value optimization under \href{https://github.com/k3nfalt/longtails/blob/main/LongTails/Setups/ZerothOrderMatrixGeneralizedSmooth.lean}{ZerothOrderMatrixGeneralizedSmooth}.}\label{tab:ZerothOrderMatrixGeneralizedSmooth}
\parbox{\linewidth}{\footnotesize Convex $f:\mathbb{R}^d\to\mathbb{R}$ has a differentiable gradient and a minimizer $x^\star$, with $\|x_0-x^\star\|\le R$ and $f(x_0)-f(x^\star)\le\Delta$. An unknown symmetric matrix $0\prec H\preceq I$ satisfies $\operatorname{tr}H\le r$ and $\nabla^2f(x)\preceq(L_0+L_1\|\nabla f(x)\|_{H^{-1}})H$. Here $d\ge1$, $L_0,r,R,\Delta,\varepsilon>0$, $L_1\ge0$, and $0<\delta<1$. Access is to exact function values; the target is $\Pr\{f(\hat x)-f(x^\star)>\varepsilon\}\le\delta$.}\par\smallskip
\begin{threeparttable}
\footnotesize
\setlength{\tabcolsep}{4pt}\renewcommand{\arraystretch}{1.25}
\begin{tabularx}{\linewidth}{@{}>{\raggedright\arraybackslash}p{0.20\linewidth}>{\raggedright\arraybackslash}p{0.32\linewidth}>{\raggedright\arraybackslash}p{0.25\linewidth}>{\raggedright\arraybackslash}X@{}}
\toprule
\textbf{Result} & \textbf{Oracle complexity} & \textbf{Novelty} & \textbf{Limitations} \\
\midrule
\textbf{(1)} \href{https://github.com/k3nfalt/longtails/blob/main/LongTails/Theorems/ZerothOrderMatrixGeneralizedAccelerated.lean}{Zeroth\allowbreak Order\allowbreak Matrix\allowbreak Generalized\allowbreak Accelerated} \textbf{(Ours)} $(\star\star\star)$ \href{https://github.com/k3nfalt/longtails/blob/main/scripts/reviews/ZerothOrderMatrixGeneralizedAccelerated.md}{(review)}
& $\mathcal{O}\!\left(D\left[\sqrt{\frac{d(r+1)L_0R^2}{\varepsilon}}\right.\right.$\newline
$\left.\left.{}+(r+1+L_1R\sqrt{d(r+1)})\ell\right]\right)$ function values, including independent runs and output selection.
& & \\
\midrule[1.2pt]
\textbf{(2)} RS-NAG \citep{zmgsOmiya2026v1}, Thm.~3.3, Prop.~5.2. Derived bound: independent runs, best-value selection.
& $\mathcal{O}\!\left(D\left[1+\sqrt{\frac{d(r+1)L_0R^2}{\varepsilon}}\right]\right)$ directional-derivative and selection-value calls. Special case $L_1=0$, $d\ge2$; same probability target.
& \textbf{(1)}: Scope extends to $L_1>0$. Values replace directional derivatives; at $L_1=0$, matches the leading rate.
& \textbf{(1)}: Adds $\mathcal{O}(D(r+1)\ell)$ at $L_1=0$. \\
\midrule
\textbf{(3)} AGD \citep{cgsTyurin2026v2}, Thm.~4.3; Euclidean reduction.
& $\mathcal{O}\!\left(1+\sqrt{\frac{L_0R^2}{\varepsilon}}\right.$\newline
$\left.{}+L_E R\log_+\min\!\left\{\frac{L_E^2\Delta}{L_0},\frac{\Delta}{\varepsilon}\right\}\right.$\newline
$\left.{}+\log_+\frac{\Delta}{L_0R^2}\right)$ exact-gradient calls. Here $R=\|x_0-x^\star\|$, $\varepsilon<\Delta$; deterministic accuracy.
& \textbf{(1)}: Values suffice; no $\mu_H^{-1/2}$ dependence from $L_E$. No established function-value cost comparison with \textbf{(3)}.
& \textbf{(1)}: Probabilistic instead of deterministic accuracy; needs $r$ additionally. \\
\midrule
\textbf{(4)} Uniform coordinate acceleration \citep{zmgsLobanov2026v1}, Thm.~6.
& $\mathcal{O}\!\left(d\widetilde R\sqrt{\frac{H_0}{\varepsilon}} +(1+d\widetilde R\sqrt{H_1})\ell_F\right)$ partial-derivative, value, and exact segment-solve calls. Expected gap $\le\varepsilon$; gap local model (Assump.~2), known $f^\star,\widetilde R<\infty$, $\varepsilon\le F_0$.
& \textbf{(1)}: Values replace derivatives and exact segment solves; $R,\Delta$ replace $f^\star,\widetilde R$. Function-value rates are not comparable: segment-solve cost is unbounded here.
& \\
\midrule
\textbf{(5)} Importance-sampled coordinate acceleration \citep{zmgsLobanov2026v1}, Thm.~8.
& $\mathcal{O}\!\left(\frac{\widetilde R\sum_i\sqrt{H_{0,i}}}{\sqrt{\varepsilon}}\right.$\newline
$\left.{}+(1+\widetilde R\sum_i\sqrt{H_{1,i}})\ell_F\right)$ partial-derivative, value, and exact segment-solve calls. Expected gap $\le\varepsilon$; coordinate local model (Assump.~7), known $f^\star,\widetilde R<\infty$, $\varepsilon\le F_0$.
& \textbf{(1)}: No coordinate-curvature inputs; $H$ remains unknown. Values replace derivatives and exact segment solves. Function-value rates are not comparable without segment-solve costs.
& \\
\bottomrule
\end{tabularx}
\begin{tablenotes}[flushleft]\footnotesize
\item[] \textbf{Abbr:} $D=1+\log(1/\delta)$; $\ell=\log(e+\Delta/\varepsilon)$; $\log_+t=\max\{0,\log t\}$, with $\log_+0=0$; $\mu_H=\lambda_{\min}(H)$; $L_E=L_1/\sqrt{\mu_H}$ (or $L_E=L_1$ when $H=I$); $\|g\|_{H^{-1}}=\sqrt{g^\top H^{-1}g}$; $f^\star=f(x^\star)$; $F_0=f(x_0)-f^\star$; $\ell_F=1+\log(F_0/\varepsilon)$; $\widetilde R=\sup_{f(x)\le f(x_0)}\|x-x^\star\|$. $H_0,H_1$ and $H_{0,i},H_{1,i}$ are the scalar and coordinate gap-local-model coefficients in \textbf{(4)} and \textbf{(5)}; they are not the matrix $H$.
\end{tablenotes}
\end{threeparttable}
\end{table}
\clearpage

\begin{table}[tp]
\centering
\caption{Fast Douglas--Rachford splitting under the \href{https://github.com/k3nfalt/longtails/blob/main/LongTails/Setups/DouglasRachfordSmoothConvex.lean}{DouglasRachfordSmoothConvex setup}.}\label{tab:DouglasRachfordSmoothConvex}
\parbox{\linewidth}{\footnotesize On a real Hilbert space, $f$ is differentiable and convex with the global quadratic upper bound of parameter $L>0$; $C$ is convex and $g$ is convex on $C$. Exact proximal minimizers of $f$ and of $g$ over $C$ are available at fixed $\gamma$. A minimizer $x^\star\in C$ exists with $\|w_1-x^\star-\gamma\nabla f(x^\star)\|\le R$. Assume $R,\varepsilon,\gamma>0$, $0<\lambda\le1$, and $4\gamma L\le\lambda$. The constructed fast Douglas--Rachford method guarantees $y_N\in C$ and $F(y_N)-F(x^\star)\le\varepsilon$ deterministically for its last iterate.}\par\smallskip
\begin{threeparttable}
\footnotesize
\setlength{\tabcolsep}{4pt}\renewcommand{\arraystretch}{1.25}
\begin{tabularx}{\linewidth}{@{}>{\raggedright\arraybackslash}p{0.20\linewidth}>{\raggedright\arraybackslash}p{0.32\linewidth}>{\raggedright\arraybackslash}p{0.25\linewidth}>{\raggedright\arraybackslash}X@{}}
\toprule
\textbf{Result} & \textbf{Oracle complexity} & \textbf{Novelty} & \textbf{Limitations} \\
\midrule
\textbf{(1)} \href{https://github.com/k3nfalt/longtails/blob/main/LongTails/Theorems/AcceleratedDouglasRachfordConvex.lean}{Accelerated\allowbreak Douglas\allowbreak Rachford\allowbreak Convex} \textbf{(Ours)} $(\star\star)$ \href{https://github.com/k3nfalt/longtails/blob/main/scripts/reviews/AcceleratedDouglasRachfordConvex.md}{(review)}
& $\mathcal{O}\!\left(1+\dfrac{R}{\sqrt{\lambda\gamma\varepsilon}}\right)$ proximal calls.
At $\lambda=1$, $\gamma=1/(4L)$: $\mathcal{O}(1+\sqrt{LR^2/\varepsilon})$.
& & \\
\midrule[1.2pt]
\textbf{(2)} Fast Douglas--Rachford \citep[Theorem~5]{patrinosStellaBemporad2014DRv2}
& $\mathcal{O}\!\left(1+\dfrac{R}{\sqrt{\lambda\gamma\varepsilon}}\right)$ proximal calls.
Special case: finite-dimensional quadratic $f$, closed proper $g+\iota_C$; $0<\gamma<1/L$, $\lambda=(1-\gamma L)/(1+\gamma L)$.
& \textbf{(1)}: Scope: extends the same method's accelerated last-iterate guarantee to nonquadratic smooth convex $f$. Matches the rate on the quadratic overlap; no problem-class speedup established.
& \textbf{(1)} needs $4\gamma L\le\lambda$; excludes \textbf{(2)}'s optimal $\gamma=(\sqrt2-1)/L$. Tuned finite-iteration constants are worse for $N\le4$.
\\
\bottomrule
\end{tabularx}
\begin{tablenotes}[flushleft]\footnotesize
\item[] \textbf{Abbr:} $F=f+g$ on $C$; $\iota_C$ is the convex indicator of $C$; $L$ is the smoothness parameter; $R$ bounds the initial distance to $x^\star+\gamma\nabla f(x^\star)$; $\gamma$ is the proximal parameter; $\lambda$ is the relaxation; $N$ is the iteration count; $\varepsilon$ is the objective-gap tolerance. Costs count one evaluation of each proximal map per iteration.
\end{tablenotes}
\end{threeparttable}
\end{table}
\clearpage

\begin{table}[tp]
\centering
\caption{Monotone root finding under mixed Lipschitz--H\"older continuity: \href{https://github.com/k3nfalt/longtails/blob/main/LongTails/Setups/HolderMonotone.lean}{HolderMonotone}.}\label{tab:HolderMonotone}
\parbox{\linewidth}{\footnotesize On a finite-dimensional real inner-product space, $F:E\to E$ is monotone and satisfies $\|F(x)-F(y)\|\le L\|x-y\|+H\|x-y\|^\alpha$, with $L,H\ge0$ and $0<\alpha\le1$. A zero $x_\star$ satisfies $\|x_0-x_\star\|\le R$, and $0<r_0\le R$, $\varepsilon>0$. The exact oracle returns $F(x)$; the target is $\|F(\hat x)\|\le\varepsilon$. Result \textbf{(1)} uses only $x_0,r_0,\varepsilon$ and returns a stopped iterate deterministically; the setup imposes no algorithm restriction.}\par\smallskip
\begin{threeparttable}
\footnotesize
\setlength{\tabcolsep}{4pt}\renewcommand{\arraystretch}{1.25}
\begin{tabularx}{\linewidth}{@{}>{\raggedright\arraybackslash}p{0.20\linewidth}>{\raggedright\arraybackslash}p{0.32\linewidth}>{\raggedright\arraybackslash}p{0.25\linewidth}>{\raggedright\arraybackslash}X@{}}
\toprule
\textbf{Result} & \textbf{Oracle complexity} & \textbf{Novelty} & \textbf{Limitations} \\
\midrule
\textbf{(1)} \href{https://github.com/k3nfalt/longtails/blob/main/LongTails/Theorems/UniversalHolderMonotone.lean}{Universal\allowbreak Holder\allowbreak Monotone} \textbf{(Ours)} $(\star\star)$ \href{https://github.com/k3nfalt/longtails/blob/main/scripts/reviews/UniversalHolderMonotone.md}{(review)}
& $\mathcal{O}\!\left((S+\log\rho)[1+\log(S\rho)]\right)$ & & \\
\midrule[1.2pt]
\textbf{(2)} PolyakEG-LS \citep{hmYoon2026v1}, Thm.~3.4(12), Prop.~3.5.
& $\begin{aligned}\mathcal{O}\bigl(&1+R^2(2M/\varepsilon)^{2/\nu}\\&+(\lambda R/\varepsilon)^2\\&+\log_+\!\tfrac{(2M)^{1/\nu}}{\lambda\varepsilon^{(1-\nu)/\nu}}\bigr)\end{aligned}$.
Pure H\"older: $L=0$, $(M,\nu)=(H,\alpha)$; pure Lipschitz: $H=0$, $(M,\nu)=(L,1)$.
& \textbf{(1)}: For $L=0$, $\alpha=1/2$, $H=R=r_0=1$, fixed $\lambda$, $\varepsilon\downarrow0$: $\mathcal{O}(\varepsilon^{-2}\log(1/\varepsilon))$ instead of $\mathcal{O}(\varepsilon^{-4})$. Scope: mixed continuity. Matches parameter-free updates. At $H=0$, fixed $L,R,r_0,\lambda$, $\varepsilon\downarrow0$: $\mathcal{O}(\varepsilon^{-1}\log(1/\varepsilon))$ instead of $\mathcal{O}(\varepsilon^{-2})$. & \\
\midrule
\textbf{(3)} Fast composite extragradient \citep{hmZhang2026v1}, Thm.~2.
& $\mathcal{O}_{H,R,\alpha}\!\left(1+[\varepsilon^{-2}\ell_\alpha]^{1/(3\alpha-1)}\right)$.
Comparison: $L=0$, $G=0$, one H\"older operator, $\alpha>1/3$. Parameter dependence not explicit.
& \textbf{(1)}: Fixed $H,R,r_0$, $1/3<\alpha<1$, $\varepsilon\downarrow0$: exponent $1/\alpha$ instead of $2/(3\alpha-1)$; no knowledge of $H,\alpha$. Also covers $\alpha\le1/3$.
& \textbf{(1)}: Stopped rather than last iterate; extra logarithm at $\alpha=1$. Equations in finite dimension rather than Hilbert-space inclusions with sums of operators. \\
\midrule
\textbf{(4)} Tseng splitting \citep{hmZhang2026v1}, Thm.~1.
& $\mathcal{O}_{H,R,\alpha}\!\left(1+[\varepsilon^{-2}\ell_\alpha]^{1/\alpha}\right)$.
Comparison: $L=0$, $G=0$, one H\"older operator. Parameter dependence not explicit.
& \textbf{(1)}: Fixed $H,R,r_0,\alpha$, $\varepsilon\downarrow0$: exponent $1/\alpha$ instead of $2/\alpha$, up to logarithms; no knowledge of $H,\alpha$.
& \textbf{(1)}: Equations in finite dimension rather than Hilbert-space inclusions with sums of operators. \\
\midrule
\textbf{(5)} Anchored extragradient \citep{bscYoonRyu2021}, Cor.~2, Sec.~2.3.
& $\mathcal{O}(1+LR/\varepsilon)$. Special case $H=0$; known $L$.
& \textbf{(1)}: No knowledge of $L$; scope extends to mixed Lipschitz--H\"older continuity.
& \textbf{(1)}: Logarithmic overhead at $H=0$ (also radius-guess overhead); stopped rather than last iterate. \\
\bottomrule
\end{tabularx}
\begin{tablenotes}[flushleft]\footnotesize
\item[] \textbf{Abbr:} $S=1+LR/\varepsilon+(8HR^\alpha/\varepsilon)^{1/\alpha}$; $\rho=R/r_0$; $\ell_\alpha=1+(1-\alpha)\log_+(1/\varepsilon)$; $\log_+t=\max\{0,\log t\}$, with $\log_+0=0$. Costs count all evaluations of $F$, including rejected trials. In \textbf{(2)}, $\lambda>0$ is the initial inverse-step estimate; line-search constants are fixed at $A=1/2$, $\beta=3$. In \textbf{(3)}, \textbf{(4)}, $G$ is the set-valued inclusion term; $\mathcal{O}_{H,R,\alpha}$ holds these parameters fixed.
\end{tablenotes}
\end{threeparttable}
\end{table}
\clearpage

\begin{table}[tp]
\centering
\caption{Noisy value-oracle convex optimization under \href{https://github.com/k3nfalt/longtails/blob/main/LongTails/Setups/ZerothOrderSmoothAndConvex.lean}{ZerothOrderSmoothAndConvex}.}\label{tab:ZerothOrderSmoothAndConvex}
\parbox{\linewidth}{\footnotesize A real Hilbert space is linearly isometric to $\mathbb R^d$, $d\ge1$. The differentiable convex function $f$ has $L$-Lipschitz gradient and a minimizer $x_\star$ with $\|x_0-x_\star\|\le R$. Each value query uses a fresh IID seed $\xi\sim p$: jointly measurable $g$ satisfies $\int g(x,\xi)\,dp(\xi)=f(x)$ and $\int(g(x,\xi)-f(x))^2\,dp(\xi)\le\sigma^2$, with these integrands integrable. Known $d,L,R,\sigma,\varepsilon$ satisfy $L,\varepsilon>0$, $R,\sigma\ge0$. Target: measurable output $\widehat x$ with integrable gap and $\int[f(\widehat x(\omega))-f(x_\star)]\,dP(\omega)\le\varepsilon$, where $P$ is the execution law.}\par\smallskip
\begin{threeparttable}
\footnotesize
\setlength{\tabcolsep}{4pt}\renewcommand{\arraystretch}{1.25}
\begin{tabularx}{\linewidth}{@{}>{\raggedright\arraybackslash}p{0.20\linewidth}>{\raggedright\arraybackslash}p{0.32\linewidth}>{\raggedright\arraybackslash}p{0.25\linewidth}>{\raggedright\arraybackslash}X@{}}
\toprule
\textbf{Result} & \textbf{Oracle complexity} & \textbf{Novelty} & \textbf{Limitations} \\
\midrule
\textbf{(1)} \href{https://github.com/k3nfalt/longtails/blob/main/LongTails/Theorems/AcceleratedTwoPointZerothOrderConvex.lean}{Accelerated\allowbreak Two\allowbreak Point\allowbreak Zeroth\allowbreak Order\allowbreak Convex} \textbf{(Ours)} $(\star\star)$ \href{https://github.com/k3nfalt/longtails/blob/main/scripts/reviews/AcceleratedTwoPointZerothOrderConvex.md}{(review)}
& $\mathcal{O}\!\left(\begin{aligned}&1+d\sqrt{LR^2/\varepsilon}\\&+d^2LR^2\sigma/\varepsilon^2\\&+d^2LR^2\sigma^2/\varepsilon^3\end{aligned}\right)$ & & \\
\midrule
\textbf{(2)} \href{https://github.com/k3nfalt/longtails/blob/main/LongTails/Theorems/AcceleratedZerothOrderConvex.lean}{Accelerated\allowbreak Zeroth\allowbreak Order\allowbreak Convex} \textbf{(Ours)} \href{https://github.com/k3nfalt/longtails/blob/main/scripts/reviews/AcceleratedZerothOrderConvex.md}{(review)}
& $\mathcal{O}\!\left(\begin{aligned}&d+d\sqrt{LR^2/\varepsilon}\\&+d^2LR^2\sigma/\varepsilon^2\\&+d^3LR^2\sigma^2/\varepsilon^3\end{aligned}\right)$ & & \\
\midrule[1.2pt]
\textbf{(3)} Derived bound: retuned two-point SGD \citep{zoscBachPerchet2016}, Prop.~7, $\beta=2$
& $\mathcal{O}\!\left(\begin{aligned}&1+dLR^2/\varepsilon\\&+d^2LR^2\sigma^2/\varepsilon^3\\&+d^2LR^2\min\!\left\{\frac{\sigma}{\varepsilon^2},\frac{d\sigma^2}{\varepsilon^3}\right\}\end{aligned}\right)$
& \textbf{(1)}: $\mathcal{O}(d\sqrt T)$ instead of $\mathcal{O}(dT)$ when $T\ge1$, $\sigma\le\varepsilon/(d\sqrt T)$. Matches order for $T\ge1$, $\sigma\ge\varepsilon/d$; unchanged quadratic-noise term.
& \\
\midrule
\textbf{(4)} Derived bound: box smoothing and AC-SA \citep{zoscLan2008}, Cor.~3
& $\mathcal{O}\!\left(\begin{aligned}&d+d\sqrt{LR^2/\varepsilon}\\&+d^2LR^2\sigma/\varepsilon^2\\&+d^3LR^2\sigma^2/\varepsilon^3\end{aligned}\right)$\newline Independent coordinate differences; all $2d$ values per gradient included.
& \textbf{(1)}: $d^2$ instead of $d^3$ in the quadratic-noise term. \textbf{(2)}: matches known rate; Lean formalization.
& \\
\midrule
\textbf{(5)} Two-point SGD \citep{zoscBachPerchet2016}, Prop.~7, published tuning
& $\mathcal{O}\!\left(d^2+\dfrac{d^2(LR^2+\sigma^2/(LR^2))^3}{\varepsilon^3}\right)$\newline Rescaled unconstrained $\beta=2$ guarantee; $T\ge1$.
& \textbf{(2)}: better bound for $\sigma\le LR^2/\sqrt d$; at $\sigma=0$, $\mathcal{O}(d\sqrt T)$ instead of $\mathcal{O}(d^2T^3)$.
& \textbf{(2)}: needs $\sigma$ for tuning; up to $d$ times worse near $\sigma=LR^2$ in the noise-dominated regime $LR^2/\sqrt d<\sigma<d^{1/4}LR^2$.
\\

\bottomrule
\end{tabularx}
\begin{tablenotes}[flushleft]\footnotesize
\item[] \textbf{Abbr:} Costs count individual function-value calls, including all batches. $T=LR^2/\varepsilon$; $\sigma$ bounds noise standard deviation. In \textbf{(7)}, $M$ is the Lipschitz constant of systematic error and $\delta$ bounds adversarial error. In \textbf{(8)}, $G$ is a global Lipschitz constant, $S=\sigma+GR/d$, $\ell=1+\log\max\{1,d,GR,S/\varepsilon\}$, and $\psi_2$ denotes the sub-Gaussian Orlicz norm.
\end{tablenotes}
\end{threeparttable}
\end{table}
\clearpage
\begin{table}[tp]
\centering
\caption{ZerothOrderSmoothAndConvex (continued).}
\begin{threeparttable}
\footnotesize
\setlength{\tabcolsep}{4pt}\renewcommand{\arraystretch}{1.25}
\begin{tabularx}{\linewidth}{@{}>{\raggedright\arraybackslash}p{0.20\linewidth}>{\raggedright\arraybackslash}p{0.32\linewidth}>{\raggedright\arraybackslash}p{0.25\linewidth}>{\raggedright\arraybackslash}X@{}}
\toprule
\textbf{Result} & \textbf{Oracle complexity} & \textbf{Novelty} & \textbf{Limitations} \\
\midrule
\textbf{(6)} Accelerated random search \citep{zoscNesterovSpokoiny2017}, Thm.~9
& $\mathcal{O}(1+d\sqrt{LR^2/\varepsilon})$\newline Exact-value restriction: $\sigma=0$.
& \textbf{(1)}: matches noiseless rate and extends to independent finite-variance noise. \textbf{(2)}: matches for $T\ge1$ and extends to noisy values.
& \\
\midrule
\textbf{(7)} \citep{zoscLoginov2026v1}, Thm.~4(b), $V=L$ (subsequent work)
& $\mathcal{O}\!\left(d+dLR^2/\varepsilon+d^2LR^2\sigma^2/\varepsilon^3\right)$\newline Additional H\"older smoothness $\beta>2$; $R>0$. Allows $M\lesssim\varepsilon/R$, $\delta\lesssim\varepsilon^{3/2}/(R\sqrt{dL})$, with sufficiently small constants.
& \textbf{(1)}: $\mathcal{O}(d\sqrt T)$ instead of $\mathcal{O}(dT)$ when $T\ge1$, $\sigma\le\varepsilon/(d\sqrt T)$, $M=\delta=0$. \textbf{(1)}, \textbf{(2)}: no $\beta>2$ requirement.
& \textbf{(1)}, \textbf{(2)}: require unbiased values without systematic or adversarial errors.
\\
\midrule
\textbf{(8)} \citep{zoscLattimoreGyorgy2023v1}, Thm.~1; high-probability conversion
& $N=\mathcal{O}\!\left(\begin{aligned}&1+d^3S^2\ell^8/\varepsilon^2\\&+d^3S\ell^4/\varepsilon\end{aligned}\right)$\newline Globally $G$-Lipschitz $f$, $GR\ge1$; sub-Gaussian noise of $\psi_2$ norm $\le\sigma$. Averaged-point gap $\le\varepsilon$ with probability $\ge1-7/N$.
& \textbf{(1)}, \textbf{(2)}: permit finite-variance noise and non-Lipschitz objectives, including quadratics. Rates are not directly comparable across expected and high-probability criteria.
& \textbf{(1)}, \textbf{(2)}: need smoothness; give expected accuracy only, rather than a high-probability guarantee.
\\
\bottomrule
\end{tabularx}
\begin{tablenotes}[flushleft]\footnotesize
\item[] \textbf{Abbr:} Costs count individual function-value calls, including all batches. $T=LR^2/\varepsilon$; $\sigma$ bounds noise standard deviation. In \textbf{(7)}, $M$ is the Lipschitz constant of systematic error and $\delta$ bounds adversarial error. In \textbf{(8)}, $G$ is a global Lipschitz constant, $S=\sigma+GR/d$, $\ell=1+\log\max\{1,d,GR,S/\varepsilon\}$, and $\psi_2$ denotes the sub-Gaussian Orlicz norm.
\end{tablenotes}
\end{threeparttable}
\end{table}
\clearpage

\begin{table}[tp]
\centering
\caption{High-probability convex optimization under \href{https://github.com/k3nfalt/longtails/blob/main/LongTails/Setups/StochasticGeneralizedSmoothAndConvex.lean}{StochasticGeneralizedSmoothAndConvex}.}\label{tab:StochasticGeneralizedSmoothAndConvex}
\parbox{\linewidth}{\footnotesize Convex $f$ on a separable real Hilbert space has differentiable gradient and $\|\nabla^2 f(x)\|\le L_0+L_1\|\nabla f(x)\|$. A minimizer $x^\star$ satisfies $\|x_0-x^\star\|\le R$ and $f(x_0)-f(x^\star)\le\Delta$. The jointly measurable oracle is unbiased, with integrable gradient and squared error, $\mathbb E\|g(x,\xi)-\nabla f(x)\|^2\le\sigma^2$, and $\mathbb E\exp(\|g(x,\xi)-\nabla f(x)\|^2/\sigma^2)\le e$ for $\sigma>0$; $\sigma=0$ means exact gradients almost surely. Samples are IID with marginal $p$ under probability law $P$. Parameters $L_0,R,\Delta,\varepsilon>0$, $L_1,\sigma\ge0$, and $\delta\in(0,1)$ are supplied. Target: measurable $\hat x$ with $P\{f(\hat x)-f(x^\star)>\varepsilon\}\le\delta$. Costs count individual stochastic-gradient samples, including batches and trials.}\par\smallskip
\begin{threeparttable}
\footnotesize
\setlength{\tabcolsep}{4pt}\renewcommand{\arraystretch}{1.25}
\begin{tabularx}{\linewidth}{@{}>{\raggedright\arraybackslash}p{0.20\linewidth}>{\raggedright\arraybackslash}p{0.32\linewidth}>{\raggedright\arraybackslash}p{0.25\linewidth}>{\raggedright\arraybackslash}X@{}}
\toprule
\textbf{Result} & \textbf{Oracle complexity} & \textbf{Novelty} & \textbf{Limitations} \\
\midrule
\textbf{(1)} \href{https://github.com/k3nfalt/longtails/blob/main/LongTails/Theorems/StochasticGeneralizedAcceleratedConvex.lean}{Stochastic\allowbreak Generalized\allowbreak Accelerated\allowbreak Convex} \textbf{(Ours)} $(\star)$ \href{https://github.com/k3nfalt/longtails/blob/main/scripts/reviews/StochasticGeneralizedAcceleratedConvex.md}{(review)}
& $\displaystyle\mathcal{O}\!\left(\begin{aligned}&1+\sqrt{\frac{L_0R^2}{\varepsilon}}\\&+\min\{\lceil t^{2/3}\rceil\ell,\,\ell+(t\ell)^{2/3}\}\\&+\frac{\sigma^2R^2}{\varepsilon^2}\log\frac e\delta\end{aligned}\right)$
& & \\
\midrule[1.2pt]
\textbf{(2)} ARC-SG \citep[Thm.~3]{gssgDvinskikh2026}\newline Finite-dimensional case.
& $\displaystyle\widetilde{\mathcal{O}}\!\left(\begin{aligned}&1+t+\sqrt{\frac{(1+t)L_0R^2}{\varepsilon}}\\&+(1+t)^3\frac{\sigma^2R^2}{\varepsilon^2}\\&+(1+t)^3\frac{\sigma^2L_1^2}{L_0^2}\end{aligned}\right)$
& \textbf{(1)} removes $(1+t)^{1/2}$ from optimization and $(1+t)^3$ from noise, plus burn-in: polynomial termwise gains for $t\gg1$. Scope: separable Hilbert spaces.
& \textbf{(1)}'s curvature term $\ell+(t\ell)^{2/3}$ can exceed $t$ through logarithms when $\ell^2\gtrsim t$; no uniform total-cost dominance. \\
\midrule
\textbf{(3)} RCSAG \citep[Thm.~9]{gssgLobanov2026}\newline Derived bound: finite-dimensional case, $p=2$, additive noise, $Q=B(x_0,R)$.
& $\displaystyle\mathcal{O}\!\left(\left[\begin{aligned}&\ell+\sqrt{\frac{L_0R^2}{\varepsilon}}+tJ\\&+\frac{\sigma^2R^2}{\varepsilon^2}\end{aligned}\right]\log\frac{e\ell}{\delta}\right)$
& \textbf{(1)}: $tJ$ becomes $\ell+(t\ell)^{2/3}$. For $\sigma=0$, $\Delta=L_0R^2$, $\varepsilon=L_0R^2/t^{4/3}$, $t\to\infty$, fixed $\delta$: $\mathcal{O}(t^{2/3}\log^{2/3}t)$ versus \textbf{(3)}'s $\mathcal{O}(t\log t\log\log t)$. Otherwise unchanged terms can dominate.
& \textbf{(1)} needs light tails and batching; \textbf{(3)} allows finite $p$th moments, $1<p\le2$, and gradient/gap-dependent noise. Logarithmic losses possible when $L_1^2\Delta\lesssim L_0$ and $\ell^2\gg t$. \\
\midrule
\textbf{(4)} Clipped SGD \citep[Thm.~1]{sgacGaash2025v2}\newline Finite-dimensional case.
& $\displaystyle\mathcal{O}\!\left(\begin{aligned}&1+\left(\frac{L_0R^2}{\varepsilon}+t^2\right)\log\frac{eT}{\delta}\\&+\frac{\sigma^2R^2}{\varepsilon^2}\log^3\frac{eT}{\delta}\end{aligned}\right)$
& \textbf{(1)}: optimization $L_0R^2/\varepsilon\to\sqrt{L_0R^2/\varepsilon}$; startup $t^2\to\ell+(t\ell)^{2/3}$. For $\sigma^2\log(e/\delta)\le L_0\varepsilon$, the optimization-term gain persists; total savings depend on noise/startup. Noise-dominated gain: $\log^3(eT/\delta)\to\log(e/\delta)$. Scope: Hilbert spaces; stopped iterate instead of averaging.
& \textbf{(1)} additionally needs $\Delta$ and batching. \textbf{(4)} uses two samples per step. \\
\bottomrule
\end{tabularx}
\begin{tablenotes}[flushleft]\footnotesize
\item[] \textbf{Abbr:} $t=L_1R$, $\ell=1+\log(1+\Delta/\varepsilon)$, $J=\min\{\ell,1+\log(1+L_1^2\Delta/L_0)\}$. In \textbf{(3)}, $H_0=\mathcal{O}(L_0)$, $H_1=\mathcal{O}(L_1^2)$, $\sigma_0=\sigma$, $\sigma_1=\sigma_2=0$; $p$ there is the moment order. $T$ is the sufficient step count in \textbf{(4)} (two calls per step), defined implicitly by its bound. $\widetilde{\mathcal{O}}$ in \textbf{(2)} suppresses polylogarithms in $1+t$, $1+L_1^2\Delta/L_0$, $1+L_0R^2/\varepsilon$, and $1/\delta$. $\lceil\cdot\rceil$ denotes upward rounding.
\end{tablenotes}
\end{threeparttable}
\end{table}
\clearpage
\endgroup

\let\table\comparisonTable

\clearpage
\section{Additional Comparisons}
\label{sec:add_comparisons}
The following tables complement the results discussed in Section~\ref{gen_inst} and Appendix~\ref{sec:comparions}. Here, we present our additional results and compare them with previous approaches. See also other new results on our \href{https://k3nfalt.github.io/longtails/}{website}. The following tables are AI-generated using the script \href{https://github.com/k3nfalt/longtails/blob/main/scripts/run_results_table_prompt.py}{\texttt{run\_results\_table\_prompt.py}}.
\def\table[#1]{\let\table\comparisonTable\comparisonTable[H]%
  \def\endtable{\end{adjustwidth}\csname float@endH\endcsname}}
\input{results_table_approved.tex}
\let\table\comparisonTable

\appendix

\end{document}